\documentclass[12pt,leqno]{article}

\usepackage{amsmath}
\usepackage{amscd}
\usepackage{amsmath}
\usepackage{amssymb}
\usepackage{latexsym}
\makeatletter

\@addtoreset{equation}{section}
\makeatother

\newtheorem{thm}{Theorem}[subsection]
\newtheorem{pr}[thm]{Proposition}
\newtheorem{df}[thm]{Definition}
\newtheorem{lm}[thm]{Lemma}
\newtheorem{cor}[thm]{Corollary}
\newtheorem{cn}[thm]{Conjecture}
\newtheorem{rmk}[thm]{Remark}

\newcommand{\qed}{\hfill{\rule{7pt}{7pt}}}

\input xy \xyoption{all} \CompileMatrices

\begin{document}
\title{Ramification theory
for varieties over a perfect field}
\author{{\sc Kazuya Kato and Takeshi Saito}}
\maketitle

\begin{abstract}
For an $\ell$-adic sheaf on a variety
of arbitrary dimension over a perfect field,
we define the Swan class
measuring the wild ramification
as a 0-cycle class supported on the ramification locus.
We prove a Lefschetz trace formula
for open varieties and 
a generalization of the 
Grothendieck-Ogg-Shararevich formula
using the Swan class.
\end{abstract}

Let $F$ be a perfect field
and $U$ be a separated and smooth
scheme of finite type
purely of dimension $d$ over $F$.
In this paper,
we study 
ramification of 
a finite \'etale scheme 
$V$ over $U$
along the boundary,
by introducing a map (\ref{eqmap}) below.

We put $CH_0(\overline V\setminus V)=
\varprojlim CH_0(Y\setminus V)$
where $Y$ runs compactifications of $V$
and the transition maps
are proper push-forwards
(Definition \ref{dfCHV}).
The degree maps $CH_0(Y\setminus V)
\to {\mathbb Z}$
induce a map $\deg:CH_0(\overline V\setminus V)
\to {\mathbb Z}$.
The fiber product $V\times_UV$
is smooth purely of dimension $d$ and
the diagonal
$\Delta_V:V\to V\times_UV$
is an open and closed immersion.
Thus the complement
$V\times_UV\setminus \Delta_V$
is also smooth purely of dimension $d$
and the Chow group
$CH_d(V\times_UV\setminus \Delta_V)$
is the free abelian group
generated by the classes
of connected components
of $V\times_UV$ not contained in $\Delta_V$.
If $U$ is connected and if
$V\to U$ is a Galois covering,
the scheme
$V\times_UV$ is the disjoint union
of the graphs $\Gamma_\sigma$
for $\sigma \in G={\rm Gal}(V/U)$ and
the group $CH_d(V\times_UV\setminus \Delta_V)$
is identified with
the free abelian group generated by
$G-\{1\}$. 

The intersection of
a connected component of $V \times_U V \setminus \Delta_V$
with $\Delta_V$ is empty.
However,
we define the intersection product with
the logarithmic diagonal 
\begin{equation}
\begin{CD}
(\ ,\Delta_{\overline V})^{\log}:
CH_d(V\times_UV\setminus \Delta_V)
@>>>
CH_0(\overline V\setminus V)\otimes_{\mathbb Z}{\mathbb Q}
\end{CD}\label{eqmap}
\end{equation}
using log product and alteration
(Theorem \ref{corindep}).
The aim of this paper is
to show that the map (\ref{eqmap})
gives generalizations
to an arbitrary dimension
of the classical 
invariants of wild ramification
of $f:V\to U$.
The image of the map is in fact supported
on the wild ramification locus
(Proposition \ref{prtame}.2).
If we have a strong form of
resolution of singularities,
we do not need $\otimes_{\mathbb Z}{\mathbb Q}$
to define the map (\ref{eqmap}).
We prove
a Lefschetz trace formula
for open varieties
\begin{equation}
\sum_{q=0}^{2d}(-1)^q
{\rm Tr}(\Gamma^*:H^q_c(V_{\bar F},{\mathbb Q}_\ell))
={\rm deg}\ (\Gamma,\Delta_{\overline V})^{\log}
\label{eqGal0}
\end{equation}
in Proposition \ref{prtr}.
If $V\to U$ is a Galois covering
of smooth curves,
the log Lefschetz class 
$(\Gamma_\sigma,\Delta_{\overline V})^{\log}$
for $\sigma\in {\rm Gal}(V/U)\setminus\{1\}$
is an equivalent of
the classical Swan character
(Lemma \ref{lmcurve}).

For a smooth $\ell$-adic sheaf ${\cal F}$
on $U$
where $\ell$ is a prime number
different from the characteristic of $F$,
we define the Swan class
${\rm Sw}({\cal F})\in 
CH_0(\overline U\setminus U)\otimes_{\mathbb Z}{\mathbb Q}$
(Definition \ref{dfSwan2})
also using the map (\ref{eqmap}).
From the trace formula (\ref{eqGal0}),
we deduce a formula
\begin{equation}
\chi_c(U_{\bar F},{\cal F})
={\rm rank}\ {\cal F}\cdot 
\chi_c(U_{\bar F},{\mathbb Q}_\ell)-
{\rm deg}\ {\rm Sw}({\cal F})
\label{eqgos}
\end{equation}
for the Euler characteristic
$\chi_c(U_{\bar F},{\cal F})
=\sum_{q=0}^{2d}(-1)^q
\dim H^q_c(U_{\bar F},{\cal F})$
in Theorem \ref{thmgos}.
If $U$ is a smooth curve,
we have
${\rm Sw}({\cal F})=
\sum_{x\in \overline U\setminus U}
{\rm Sw}_x({\cal F})[x]$
by Lemma \ref{lmdim1}.
Thus the formula
(\ref{eqgos})
is nothing other than 
the Grothendieck-Ogg-Shafarevich formula
\cite{SGA5}, \cite{GOS}.
As a generalization of the Hasse-Arf theorem
(Lemma \ref{lmdim1}),
we state Conjecture \ref{cnint}
asserting that we do not need
$\otimes_{\mathbb Z}{\mathbb Q}$
in the definition of the Swan class.
We prove a part of Conjecture \ref{cnint}
in dimension 2 (Corollary \ref{corint}.1).

The profound insight that the wild ramification
gives rise to invariants as 0-cycle classes
supported on the ramification locus
is due to S.~Bloch \cite{bloch}
and is developed by one of the authors
in \cite{ICM}, \cite{Kato}.
Since a covering ramifies along a divisor
in general,
it is naturally expected that
the invariants defined as 0-cycle classes
should be computable in terms of
the ramification
at the generic points
of irreducible components
of the ramification divisor.
For the log Lefschetz class 
$(\Gamma_\sigma,\Delta_{\overline V})^{\log}$,
we give such a formula
(\ref{eqsgrd})
in Lemma \ref{lmsgr}.
For the Swan class of a sheaf
of rank 1,
we state Conjecture \ref{cncF}
in this direction
and prove it assuming $\dim U\le 2$
in Theorem \ref{thmint}.
We expect that the log
filtration by ramification groups
defined in \cite{AS}
should enable us to compute the Swan classes
of sheaves of arbitrary rank.

In a subsequent paper,
we plan to study ramification of
schemes over a discrete valuation ring
and prove an analogue of
Grothendieck-Ogg-Shafarevich formula
for the Swan conductor of
cohomology (cf. \cite{abbes}, \cite{abbes2}).
In $p$-adic setting,
the relation between the Swan conductor
and the irregularities are studied in 
\cite{Crew}, \cite{Crew2},
\cite{Matsuda} and \cite{Tsuzuki}.
The relation between the Swan classes
defined in this paper
and the characteristic varieties
of ${\cal D}$-modules 
defined in \cite{Berthelot}
should be investigated.

In Section 1,
we recall a log product construction
in \cite{KS}.
In Section 2,
we prove a Lefschetz trace formula
Theorem \ref{thmLTF}
for algebraic correspondences on open varieties,
under a certain assumption.
In Section 3,
we define and study the map
$(\ref{eqmap})$
and prove the trace formula
(\ref{eqGal0}) in
Proposition \ref{prtr}.
In Section 4,
we define
the Swan class of an
$\ell$-adic sheaf and prove 
the formula (\ref{eqgos})
in Theorem \ref{thmgos}.
In Section 5,
we compare the Swan class
in rank 1 case
with an invariant defined in \cite{Kato}.
We also compare the formula
(\ref{eqgos})
with a formula of Laumon
in dimension 2.

\noindent
{\bf Acknowledgement}
The authors are grateful to
Ahmed Abbes
and the referee
for thorough reading
and helpful comments.
They thank
Ahmed Abbes,
H\'el\`ene Esnault
and Luc Illusie
for stimulating discussions
and their interests.
The authors are grateful to
Shigeki Matsuda
for pointing out that
the assumption of
Theorem in \cite{Sa}
is too weak to
deduce the conclusion.
A corrected assumption is given in
Proposition \ref{prrk1}.
The authors are grateful to Bruno Kahn
for showing Lemma \ref{lmCV}.

\tableofcontents

\bigskip 
\noindent
{\Large \bf Notation}

\medskip
In this paper, we fix a base field $F$.
A scheme means a separated scheme
of finite type over $F$
unless otherwise stated explicitly.
For schemes $X$ and $Y$ over $F$,
the fiber product over $F$
will be denoted by $X\times Y$.

The letter $\ell$ denotes
a prime number invertible in $F$.

\section{Log products}

In \S1.1,
we introduce log products
and establish elementary properties.
In \S1.2, we define and study
admissible automorphisms.

\subsection{Log blow-up and log product}

We introduce log blow-ups and log products
with respect to families of Cartier divisors.

\begin{df}\label{dflogbp}
Let $F$ be a field and
let $X$ and $Y$ be separated schemes
of finite type over $F$.
Let 
${\cal D}=(D_i)_{i\in I}$
be a finite family
of Cartier divisors of $X$
and
${\cal E}=(E_i)_{i\in I}$
be a finite family
of Cartier divisors of $Y$
indexed by the same finite set $I$.

For $i\in I$,
let $(X\times Y)_i'\to X\times Y$
be the blow-up at 
$D_i\times E_i
\subset X\times Y$
and let $(X\times Y)_i^\sim
\subset
(X\times Y)_i'$
be the complement of the 
proper transforms of
$D_i\times Y$
and $X\times E_i$.

1.
We define the log blow-up
\begin{equation}
\begin{CD}
p:(X\times Y)'@>>> X\times Y,
\end{CD}
\label{eqblup}
\end{equation}
more precisely denoted by
$((X,{\cal D})\times (Y,{\cal E}))'$,
to be the fiber product 
${\prod_{i\in I}}_X (X\times Y)_i'
\to X\times Y$ of $(X\times Y)_i'\ (i \in I)$
over $X\times Y$.

2. Similarly,
we define the log product
\begin{equation}
(X\times Y)^\sim\subset(X\times Y)',
\label{eqprod}
\end{equation}
or more precisely denoted by
$((X,{\cal D})\times (Y,{\cal E}))^\sim$,
to be the fiber product 
${\prod_{i\in I}}_X (X\times Y)_i^\sim
\to X\times Y$ of $(X\times Y)_i^\sim\ (i \in I)$
over $X\times Y$.

3. If $X=Y$
and ${\cal D}={\cal E}$,
we call 
$(X\times X)^\sim$
the log self product
of $X$ with respect to 
${\cal D}$.
By the universality of blow-up,
the diagonal map
$\Delta:
X\to X\times X$
induces an immersion 
$$X\to (X\times X)^\sim$$
called the log diagonal map.
\end{df}

Locally on $X$ and $Y$,
the log blow-up,
log self-product
and the log diagonal maps 
are described as follows.

\begin{lm}\label{lmlogbp}
Let the notation be as in
Definition \ref{dflogbp}.
Assume that $X={\rm Spec}\ A$ and
$Y={\rm Spec}\ B$
are affine 
and that the Cartier divisors
$D_i$ are defined by $t_i\in A$
and
$E_i$ are defined by $s_i\in B$
respectively.

1. The log product
$(X\times Y)'$
is the union of
\begin{equation}
{\rm Spec}\ \frac{A\otimes_F B
[U_i\ (i\in I_1),V_j\ (j\in I_2)]}
{(t_i\otimes 1-U_i(1\otimes s_i)\ (i\in I_1),
1\otimes s_j-V_j(t_j\otimes 1)\ (j\in I_2))}
\label{eqlogbp}
\end{equation}
for decompositions
$I=I_1\amalg I_2$.

2. The log product
$(X\times Y)^\sim$ 
is given by
\begin{equation}
{\rm Spec}\ A\otimes_F B
[U_i^{\pm1}\ (i\in I)]/
(t_i\otimes 1-U_i(1\otimes s_i)\
(i\in I))
\label{eqlogpr}
\end{equation}

3. Assume further that
$A=B$, $D_i=E_i$
and $t_i=s_i$ for each $i\in I$.
Then in the notation {\rm (\ref{eqlogpr})},
the log diagonal map $\Delta:
X\to (X\times X)^\sim$
is defined by 
the map 
\begin{equation}
A\otimes_F A
[U_i^{\pm1}\
(i\in I)]/
(t_i\otimes 1-U_i(1\otimes t_i)\
(i\in I))\to A
\label{eqlogspr}
\end{equation}
sending
$a\otimes 1$ and $1\otimes a$ to $a\in A$
and $U_i$ to 1 for $i\in I$.
\end{lm}

{\it Proof.}
For each $i\in I$,
the Cartier divisors
$D_i\times Y$ and $X\times E_i$
are locally defined
by a regular sequence.
Thus we obtain 1.
The rest is clear from this and 
the definition.
\qed

For the sake of readers familiar with
log schemes,
we recall an intrinsic definition
using log structures given in \cite{KS}.
The Cartier divisors
$D_1,\ldots,D_m$
define a log structure $M_X$ on $X$.
In the notation in Lemma \ref{lmlogbp},
the log structure $M_X$ is defined
by the chart
${\mathbb N}^m\to A$
sending the standard basis
to $t_1,\ldots,t_m$.
The local charts
${\mathbb N}^m\to A$
induce a map
${\mathbb N}^m
\to 
\Gamma(X,M_X/O_X^\times)$
of monoids.
Similarly,
the Cartier divisors
$E_1,\ldots,E_m$
defines a log structure on $Y$
and a map
${\mathbb N}^m
\to 
\Gamma(Y,M_Y/O_Y^\times)$.
Then, the log product $(X\times Y)^\sim$
represents the functor
attaching to an fs-log scheme $T$ over $F$
the set of pairs $(f,g)$
of morphisms of log schemes 
$f:T\to X$ and $g:T\to Y$ over $F$
such that
the diagram
$$\begin{CD}{\mathbb N}^m
@>>>
\Gamma(X,M_X/O_X^\times)\\
@VVV@VVV\\
\Gamma(Y,M_Y/O_Y^\times)
@>>>
\Gamma(T,M_T/O_T^\times)
\end{CD}$$
is commutative.
The log diagonal 
$\Delta:X\to (X\times X)^\sim$
corresponds to the pair
$({\rm id},{\rm id})$.

The log product satisfies the following
functoriality.
Let $X, X', Y$ and $Y'$ be schemes over $F$
and ${\cal D}=(D_i)_{i\in I}$,
${\cal D}'=(D'_i)_{i\in I}$,
${\cal E}=(E_j)_{j\in J}$,
and ${\cal E}'=(E'_j)_{j\in J}$
be families of Cartier divisors of 
$X, X', Y$ and of $Y'$ respectively.
Let $f:X\to Y$
and $g:X'\to Y'$
be morphisms over $F$
and let $e_{ij}\ge 0,(i,j)\in I\times J$
be integers satisfying
$f^*E_j=\sum_{i\in I}e_{ij}D_i$
and
$f^*E'_j=\sum_{i\in I}e_{ij}D'_i$
for $j\in J$.
Then, the maps $f$ and $g$ induces
a map
$(f\times g)^\sim:
(X\times X')^\sim \to
(Y\times Y')^\sim$.
If $Y=Y'$ and ${\cal E}={\cal E}'$,
we define
$(X\times_YX')^\sim$,
or more precisely
$((X,{\cal D})\times_{(Y,{\cal E})}
(X',{\cal D}'))^\sim$,
to be the fiber product
$(X\times X')^\sim\times_{
(Y\times Y)^\sim}Y$ with the log diagonal
$Y \to
(Y\times Y)^\sim$.

\begin{lm}\label{lmn}
Let $F$ be a field
and $n\ge 1$ be an integer.
Let $Y$ be a separated scheme over $F$.
Let ${\cal L}$ be an invertible $O_Y$-module
and
$\mu:{\cal L}^{\otimes n}\to O_Y$ be
an injection of $O_Y$-modules.
We define an $O_Y$-algebra
${\cal A}=\bigoplus_{i=0}^{n-1}
{\cal L}^{\otimes i}$
with the multiplication defined by
$\mu:{\cal L}^{\otimes n}\to O_Y$
and put
$X={\rm Spec}{\cal A}$.
Let $E$ be the Cartier divisor of $Y$
defined by ${\cal I}_E=
{\rm Im}({\cal L}^{\otimes n}\to O_Y)$
and $D$ be the Cartier divisor of $X$
defined by ${\cal L}O_X$.
Let $(X\times_Y X)^\sim$
be the log self product
defined with respect to $D$ and $E$.

We define an action of
the group scheme $\mu_n
={\rm Spec}F[t]/(t^n-1)$
on $X$ over $Y$ by the multiplication by $t$
on ${\cal L}$.
We consider the action of
$\mu_n$ on $(X\times_Y X)^\sim$
by the action on the first factor $X$.

Then, by the second projection
$(X\times_Y X)^\sim\to X$,
the scheme 
$(X\times_Y X)^\sim$ is a $\mu_n$-torsor on $X$.
Further the log diagonal map
$X\to (X\times_Y X)^\sim$
induces an isomorphism
$\mu_n\times X\to (X\times_Y X)^\sim$.
\end{lm}

\noindent{\it Proof.}
Since the question is local on $Y$,
it is reduced to the
case where 
$Y={\mathbf A}^1={\rm Spec}\ F[T]$
and $\mu$ send a basis $S^n$ to $T$.
Then we have
$X={\mathbf A}^1={\rm Spec}\ F[S]$
and the map $X\to Y$
is given by $T\mapsto S^n$.
Then, by Lemma \ref{lmlogbp}.2,
we have
$(Y\times Y)^\sim=
{\rm Spec}\ F[T,T',U^{\pm1}]/(T'-UT)=
{\rm Spec}\ F[T,U^{\pm1}]$,
$(X\times X)^\sim=
{\rm Spec}\ F[S,S',V^{\pm1}]/(S'-VS)=
{\rm Spec}\ F[S,V^{\pm1}]$,
and the map
$(X\times X)^\sim
\to (Y\times Y)^\sim$
is given by
$T\mapsto S^n$
and
$U\mapsto V^n$.
Since
the log diagonal
$Y\to (Y\times Y)^\sim$
is defined by
$U=1$,
we have
$(X\times_YX)^\sim=
{\rm Spec}\ F[S,V^{\pm1}]/(V^n-1)$.
Thus the assertion is proved.
\qed

Let $F$ be a field and
$X$ be a smooth scheme
purely of dimension $d$ over $F$.
In this paper, we say a
divisor $D$ of $X$ has
simple normal crossings if 
the irreducible components
$D_i\ (i\in I)$
are smooth over $F$
and, for each subset
$J\subset I$,
the intersection
$\bigcap_{i\in J}D_i$
is smooth purely of dimension
$d-|J|$ over $F$.
In other words,
Zariski locally on $X$,
there is an \'etale map
to ${\mathbf A}^d_F
={\rm Spec}\ F[T_1,\ldots,T_d]$
such that $D$ is the pull-back of
the divisor defined by
$T_1\cdots T_r$
for some $0\le r\le d$.
If $D_i$ is an irreducible component,
$D_i$ is smooth and
$\bigcup_{j\neq i}(D_i\cap D_j)$
is a divisor of $D_i$ with simple
normal crossings.

Let $X$ be a smooth scheme over a field $F$
and $D$
be a divisor of $X$ with simple
normal crossings.
Let $D_i\ (i\in I)$ be the irreducible 
components of $D$.
We consider the log blow-up
$p:(X\times X)'\to X\times X$
with respect to the
family $D_i\ (i\in I)$
of irreducible components of $D$,
defined in Definition \ref{dflogbp}.
Let $D^{(1)\prime}\subset (X\times X)'$ and
$D^{(2)\prime}\subset (X\times X)'$ be
the proper transforms
of $D^{(1)}=D\times X$
and of $D^{(2)}=X\times D$
respectively.
Let $E_i=(X\times X)'\times_{X\times X}
(D_i\times D_i)$ be the exceptional divisors
and $E=\bigcup_iE_i
\subset
(X\times X)'$
be the union.

The log blow-up $p:(X\times X)'\to X\times X$
is used in \cite{fal} and in \cite{Pink}
in the study of cohomology of
open varieties.
For an irreducible component $D_i$ of $D$,
the log blow-up
$(D_i\times D_i)'\to D_i\times D_i$
is defined with respect to the family
$D_i\cap D_j, j\neq i$
of Cartier divisors.

\begin{lm}\label{lm'}
Let $X$ be a smooth scheme over $F$,
$D$ be a divisor of $X$ with simple normal crossings
and $U=X\setminus D$ be the complement.
Let
$p:(X\times X)'\to X\times X$
be the log blow-up with respect to the
family of irreducible components of $D$.

1. The scheme
$(X\times X)'$ is smooth over $F$.
The complement
$(X\times X)'\setminus (U\times U)
=D^{(1)^\prime}\cup
D^{(2)^\prime}\cup E$
is a divisor with simple normal crossings.
The log product $(X\times X)^\sim$
is equal to the complement
$$
(X\times X)'\setminus (D^{(1)\prime}\cup D^{(2)\prime}).$$

2. 
Let $D_i$ be an irreducible 
component of $D$.
The projection
$E_i\to D_i\times D_i$
induces a map
$E_i\to (D_i\times D_i)'$
and further a map
$E_i^\circ=E_i\cap (X\times X)^\sim
\to (D_i\times D_i)^\sim$.
We have a canonical isomorphism
\begin{equation}
\begin{CD}
E_i@>>>
{\mathbf P}(N_{D_i\times D_i/X\times X})
\times_{D_i\times D_i}(D_i\times D_i)'
\end{CD}
\label{eqEP}
\end{equation}
to the pull-back of the
${\mathbf P}^1$-bundle
${\bf P}(N_{D_i\times D_i/X\times X})=
{\bf Proj}({\rm S}^\bullet N_{D_i\times D_i/X\times X})$
associated to the conormal sheaf
$N_{D_i\times D_i/X\times X}$.

We identify
$E_i$ with ${\mathbf P}(N_{D_i\times D_i/X\times X})
\times_{D_i\times D_i}(D_i\times D_i)'$
by the isomorphism {\rm (\ref{eqEP})}.
Then the open subscheme 
$E_i^\circ\subset E_i$
is the complement of
the two disjoint sections
$(D_i\times D_i)^\sim
\to
{\mathbf P}(N_{D_i\times D_i/X\times X})
\times_{D_i\times D_i}
(D_i\times D_i)^\sim$
defined by the surjections
$N_{D_i\times D_i/X\times X}\to
N_{D_i\times D_i/D_i\times X}$
and
$N_{D_i\times D_i/X\times X}\to
N_{D_i\times D_i/X\times D_i}$.
\end{lm}

\noindent
{\it Proof.}
1. It follows immediately from the definition and 
the description in
Lemma \ref{lmlogbp}.

2. Clear from the definition.
\qed

\begin{cor}\label{cor'}
Let the notation be as in Lemma \ref{lm'}.
Let $D_i$ be an irreducible component of $D$
and let $D_i\to (D_i\times D_i)^\sim$
be the log diagonal map.
Then the isomorphism {\rm (\ref{eqEP})}
induces an isomorphism
\begin{equation}
\begin{CD}
E_{i,D_i}^\circ=
E_i^\circ \times_{(D_i\times D_i)^\sim}D_i
@>>> {\mathbb G}_{m,D_i}.
\end{CD}
\label{eqEPD}
\end{equation}
The section $D_i\to E_{i,D_i}^\circ$
induced by the log diagonal
$X\to (X\times_YX)^\sim$
is identified with the unit section
$D_i\to {\mathbb G}_{m,D_i}$.
\end{cor}

\noindent{\it Proof.}
The restrictions
of the conormal sheaf
$N_{D_i\times D_i/X\times X}$
to the diagonal $D_i\subset
D_i\times D_i$
is the direct sum of
the restrictions 
$N_{D_i\times D_i/D_i\times X}|_{D_i}$
and 
$N_{D_i\times D_i/X\times D_i}|_{D_i}$.
Further the restrictions
$N_{D_i\times D_i/D_i\times X}|_{D_i}$
and
$N_{D_i\times D_i/X\times D_i}|_{D_i}$
are canonically isomorphic to
$N_{D_i/X}$.
Hence
we have a canonical isomorphism
${\mathbf P}(N_{D_i\times D_i/X\times X})
\times_{D_i\times D_i}D_i\to
{\mathbf P}^1_{D_i}$
and the assertion follows
from Lemma \ref{lm'}.2.
\qed

\begin{pr}\label{prmu}
Let $X$ be a separated smooth scheme purely of dimension $d$
over $F$ and $U=X\setminus D$ be the complement of
a divisor $D=\bigcup_{i\in I}D_i$
with simple normal crossings.
Let $Y$ be a separated scheme over $F$
and $V=Y\setminus B$
be the complement of a Cartier divisor $B$.
We consider a Cartesian diagram
\begin{equation}
\begin{CD}
U@>{\subset}>> X\\
@VfVV @VV{\bar f}V\\
V@>{\subset}>> Y.
\end{CD}\label{eqUV}
\end{equation}
We put $\bar f^* B=\sum_{i\in I}e_iD_i$.

1. 
Let
$(X\times X)^\sim$
be the log product with respect to
the family $(D_i)_{i\in I}$
of irreducible components
and $(Y\times Y)^\sim$
be the log product with respect to $B$.
Let $(X\times_Y X)^\sim=
(X\times X)^\sim\times_{(Y\times Y)^\sim}Y$ be the 
inverse image of the diagonal.
We  keep the notation in Corollary \ref{cor'}.
Let $D_i$ be an irreducible component of $D$.
We identify
$E_{i,D_i}^\circ=
E_i^\circ
\times_{(D_i\times D_i)^\sim}D_i$
with
${\mathbb G}_{m,D_i}$
by the isomorphism
{\rm (\ref{eqEPD})}.

Then the intersection
$E_{i,D_i}^\circ
\cap (X\times_YX)^\sim$
is a closed subscheme of
the subscheme $\mu_{e_i,D_i}\subset 
{\mathbb G}_{m,D_i}$
of $e_i$-th roots of $1$.

2. The closure $\overline {U\times_VU}$
in the log product $(X\times X)'$
satisfies the equality
\begin{equation}
\overline {U\times_VU} \cap D^{(1)\prime}
=
\overline {U\times_VU} \cap D^{(2)\prime}
\label{eqP'}
\end{equation}
of the underlying sets.
\end{pr}

\noindent{\it Proof.}
1. 
The assertion is local on $D_i\subset (D_i\times D_i)^\sim$.
Hence, we may assume
that $X={\rm Spec}\ A$
is affine
and that the divisor $D_k$
is defined by $t_k\in A$ for $k\in I$.
We may also assume
that the Cartier divisor $B$ of $Y$
is defined by 
a function $s$.
Then,
we have $f^*s=v\prod_{k\in I}
t_k^{e_k}$
for a unit $v\in A^\times$.
We identify
$(X\times X)^\sim
={\rm Spec}\ A\otimes_F A
[U_k^{\pm1}\ (k\in I)]/
(t_k\otimes 1-U_k(1\otimes t_k)\
\ (k\in I))$
as in (\ref{eqlogspr}).
Then on the closed subscheme
$(X\times_YX)^\sim\subset
(X\times X)^\sim$,
we have an equation
$$\frac {v\otimes 1}
{1\otimes v}
\prod_{k\in I}U_k^{e_k}
=1.$$
On the log diagonal
$D_i\subset (D_i\times D_i)^\sim$,
we have $v\otimes 1=1\otimes v$ and
$U_k=1$
for $k\in I\setminus \{i\}$.
Since the coordinate of the ${\mathbb G}_m$-bundle
$E_{i,D_i}$ is given by 
$U_i$,
the assertion follows.

2. It suffices to show the equality
$\overline {\Gamma} \cap D^{(1)\prime}
=
\overline {\Gamma} \cap D^{(2)\prime}$
for any integral closed subscheme
$\Gamma\subset U\times_VU$.
We regard
$\overline \Gamma$ as a closed subscheme
of $(X\times X)'$
with an integral scheme
structure
and let $p_1,p_2:
\overline \Gamma\to X$ denote the compositions
with the projections.
We consider the Cartier divisors
$p_1^*D_i$ and $p_2^*D_i$ of 
$\overline \Gamma$.
We also consider the Cartier divisors
$(D_i\times X)'
\cap \overline \Gamma$
and
$(X\times D_i)'
\cap \overline \Gamma$.

By the Cartesian diagram
(\ref{eqUV}),
we have $e_i>0$
in $X\times_YB=\sum_{i\in I}e_iD_i$
for all $i$.
Since
$\Gamma\subset U\times_VU$,
the closure $\overline \Gamma$
is a closed subscheme of
the pull-back
$(X\times X)'\times_{Y\times Y}Y$
of the diagonal.
Hence, we have
an equality 
$\sum_ie_ip_1^*D_i=\sum_ie_ip_2^*D_i$
of Cartier divisors of $\overline \Gamma$.
Thus, we have
an equality 
$\sum_ie_i
(D_i\times X)'
\cap \overline \Gamma
=\sum_ie_i
(X\times D_i)'
\cap \overline \Gamma$.
Since $e_i>0$ for all $i$,
we obtain
$$\overline \Gamma \cap D^{(1)\prime}
=
\bigcup_i
(D_i\times X)'
\cap \overline \Gamma
=
\bigcup_i
(X\times D_i)'
\cap \overline \Gamma=
\overline \Gamma \cap D^{(2)\prime}.$$
\qed

We consider tamely ramified coverings.

\begin{df}\label{dftame}

1. Let $K$
be a complete discrete
valuation field.
We say a finite separable extension
$L$ of $K$ is
tamely ramified if
the ramification index $e_{L/K}$
is invertible in the residue field
and if the extension
of the residue field is separable.

2. Let 
$$\begin{CD}
U@>{\subset}>> X\\
@VfVV @VV{\bar f}V\\
V@>{\subset}>> Y
\end{CD}$$
be a Cartesian diagram of
locally noetherian normal schemes.
We assume that 
$Y$ is regular,
$V$ is the complement of
a divisor with simple normal crossings
and that 
$U$ is a dense open subscheme of $X$.
We also assume that
the map $f:U\to V$ is finite \'etale
and $\bar f:X\to Y$ is quasi-finite.

We say $\bar f:X\to Y$
is tamely ramified
if, for each point $\xi\in X\setminus U$
such that $O_{X,\xi}$
is a discrete valuation ring,
the extension of
the complete discrete valuation fields
${\rm Frac}(\hat O_{X,\xi})$
over
${\rm Frac}(\hat O_{Y,\bar f(\xi)})$
is tamely ramified.
\end{df}

\begin{lm}\label{lmtame}
Let 
$$\begin{CD}
U@>{\subset}>> X\\
@VhVV @VV{\bar h}V\\
V'@>{\subset}>> Y'\\
@VgVV @VV{\bar g}V\\
V@>{\subset}>> Y
\end{CD}$$
be a Cartesian diagram of
separated normal schemes
of finite type over $F$.
We assume that 
$X$ and $Y$ are smooth over $F$,
$U\subset X$ and $V\subset Y$ are the complements of
divisors with simple normal crossings
and 
$V'$ is a dense open subscheme of $Y'$.
We also assume that
$g:V'\to V$ is finite \'etale
and $\bar g:Y'\to Y$ is quasi-finite
and tamely ramified.

Then, 
in $(X\times X)^\sim$,
the intersection of the closure
$\overline{U\times_VU\setminus U\times_{V'}U}$
with the log diagonal
$X\subset (X\times X)^\sim$
is empty.
\end{lm}

\noindent{\it Proof.}
The assertion is \'etale local on $X$ and on $Y$.
We put $f=g\circ h$ and $\bar f=\bar g\circ \bar h$.
Let $\bar x$ be a geometric point of $X$
and $\bar y=\bar f(x)$ be its image. 
We take \'etale maps
$Y\to {\mathbf A}^d_F=
{\rm Spec}\ F[T_1,\ldots, T_d]$
and
$X\to {\mathbf A}^n_F=
{\rm Spec}\ F[S_1,\ldots, S_n]$
such that
$V=Y\times_{{\mathbf A}^d_F}
{\rm Spec}\ F[T_1,\ldots, T_d]
[(T_1\cdots T_r)^{-1}]$
and
$U=X\times_{{\mathbf A}^n_F}
{\rm Spec}\ F[S_1,\ldots, S_n]
[(S_1\cdots S_q)^{-1}]$.
Since the assertion is \'etale local on $Y$,
we may assume 
that there exist an integer $e\ge 1$
invertible in $F$
and a surjection
$Y_e=
Y\times_{{\mathbf A}^d_F}
{\rm Spec}\ F[T_1,\ldots, T_d]
[T_1^{1/e},\ldots,T_r^{1/e}]
\to Y'$ over $Y$
by Abhyankar's lemma.
Further we may assume 
that there exists a surjection
$X_e=
X\times_{{\mathbf A}^n_F}
{\rm Spec}\ F[S_1,\ldots, S_n]
[S_1^{1/e},\ldots,S_r^{1/e}]
\to X\times_{Y'}Y_e$ over $X$.

We put
$V_e=V\times_YY_e$
and
$U_e=U\times_XX_e$.
Then, $(X_e\times X_e)^\sim
\to (X\times X)^\sim$
is finite,
$X_e\to X$ is surjective
and the inverse image of
$U\times_VU\setminus U\times_{V'}U$
is a subset of
$U_e\times_VU_e\setminus U_e\times_{V_e}U_e$.
Hence,
it is reduced to the case
where $X\to Y'$
is $X_e\to Y_e$
and further 
to the case $X_e=Y_e$.
Since $(Y_e\times_YY_e)^\sim
\to Y_e$ is finite \'etale
as in Lemma \ref{lmn},
the assertion is proved.
\qed

\subsection{Admissible automorphisms}

Let $X$ be a smooth scheme over $F$,
$D$ be a divisor of $X$ with simple normal crossings
and $U=X\setminus D$ be the complement.
We study an automorphism of $X$
stabilizing $U$.

\begin{df}\label{dfadm}
Let $X$ be a smooth scheme over $F$,
$D$ be a divisor of $X$ with simple normal crossings
and $U=X\setminus D$ be the complement.
Let $D_1,\ldots,D_m$
be the irreducible components of $D$.

Let $\sigma$ be an automorphism
of $X$ over $F$
satisfying $\sigma(U)=U$.
We say $\sigma$ is admissible
if, for each $i=1,\ldots,m$,
we have either 
$\sigma(D_i)=D_i$
or $\sigma(D_i)\cap D_i=\emptyset$.
\end{df}

We define the blow-up $X_\Sigma\to X$
associated to the subdivision by baricenters
and show that the induced action
on $X_\Sigma$ is admissible.

\begin{df}\label{dfbari}
Let $X$ be a smooth scheme 
purely of dimension $d$ over $F$,
$D$ be a divisor of $X$ with simple normal crossings
and let $D_1,\ldots,D_m$
be the irreducible components of $D$.
For a subset $I\subset \{1,\ldots,m\}$,
we put $D_I=\bigcap_{i\in I}D_i$.
We put $X=X_0$
and, for $0\le i<d$,
we define $X_{i+1}\to X_i$
to be the blow-up at the proper transforms
of $D_I$ for $|I|=d-i$ inductively.
We call $X_\Sigma=X_d \to X$
the blow-up
associated to the subdivision by baricenters.
\end{df}

\begin{lm}\label{lmdec}
Let $X$ be a smooth scheme over $F$,
$D$ be a divisor of $X$ with simple normal crossings
and let $D_1,\ldots,D_m$
be the irreducible components of $D$.
Let $U=X\setminus D$ be the complement
and let $p:X_\Sigma\to X$ be the blow-up 
associated to the subdivision by baricenters.

1.
The scheme $X_\Sigma$ is smooth over $F$
and the
complement $D'=X_\Sigma\setminus U$
is a divisor with simple normal crossings.
For an irreducible component
$D'_j$ of $D'$,
we put $I=\{i|D'_j\subset p^{-1}(D_i),1\le i\le m\}$
and $k=|I|$.
Then there exists
an irreducible component $Z$ of $D_I$
satisfying the following condition.
Let $Z'\subset X_k$ be the proper transform of $Z$
in $X_k$ and $E_Z\subset X_{k+1}$
be the inverse image of $Z'$.
Then $D'_j$ is the proper transform of $E_Z$.

2.
For an automorphism
$\sigma$ of $X$ over $F$
satisfying $\sigma(U)=U$,
the induced action of $\sigma$
on $X_\Sigma$ is admissible.
\end{lm}

\noindent{\it Proof.}
1.
It suffices to study \'etale locally on $X$.
Hence, it suffices to consider
the case where
$X={\mathbf A}^d={\rm Spec}\ F[T_1,\ldots,T_d]$
and $D$ is defined by $T_1\cdots T_m=0$.
Then $X_\Sigma$ is obtained by patching
${\rm Spec}\ A_\varphi$
where 
$$A_\varphi
=
F\left[T_{\varphi(1)},\frac{T_{\varphi(2)}}{T_{\varphi(1)}},
\ldots,\frac{T_{\varphi(m)}}{T_{\varphi(m-1)}},
T_{m+1},\ldots,T_d\right]$$
for bijections
$\varphi:\{1,\ldots,m\}\to 
\{1,\ldots,m\}$.
The assertion follows easily from this.

2.
Let $D'_1,\ldots,D'_{m'}$
be the irreducible components of $D'$
and $\Sigma=\{I\subset\{1,\ldots,m\}\}$
be the power set of $\{1,\ldots,m\}$.
We define a map
$\psi:\{1,\ldots,m'\}
\to \Sigma$
by putting
$\psi(j)=
\{i|D'_j\subset p^{-1}(D_i),1\le i\le m\}$.
Then by 1, for irreducible components
$D'_j\neq D'_{j'}$
such that
$D'_j\cap D'_{j'}\neq \emptyset$,
we have either
$\psi(j)\subsetneqq \psi(j')$
or
$\psi(j)\supsetneqq \psi(j')$.
The map
$\psi:\{1,\ldots,m'\}
\to \Sigma$
is compatible with the
natural actions of $\sigma$.
Therefore, if $\sigma(D'_j)=D'_{\sigma(j)}
\neq D'_j$,
we have $|\psi(\sigma(j))|=
|\sigma(\psi(j))|=
|\psi(j)|$
and $\sigma(D'_j)\cap D'_j=\emptyset$.
\qed

We define the log fixed part
for an admissible automorphism.

\begin{lm}\label{lmadm}
Let $X$ be a separated and smooth scheme 
of finite type over $F$,
$D$ be a divisor of $X$ with simple normal crossings
and $U=X\setminus D$ be the complement.
Let $\sigma$ be an admissible automorphism
of $X$ over $F$
satisfying $\sigma(U)=U$.
Then, the closed immersion 
$(1,\sigma):U\to U\times U$
is extended to a closed immersion
\begin{equation}
\begin{CD}
\tilde \Gamma_\sigma:
X\setminus 
\bigcup_{i:\sigma(D_i)\neq D_i}
D_i
@>>>
(X\times X)^\sim.
\end{CD}
\label{eqsigma}
\end{equation}
\end{lm}

\noindent{\it Proof.}
By the assumption that
$\sigma$ is admissible,
the closed immersion $(1,\sigma):
X\to X\times X$ induces
a closed immersion
$X\to (X\times X)'$.
Let $\Gamma'_\sigma$
denote $X$
regarded as 
a closed subscheme of $(X\times X)'$
by this immersion.
Then, it induces an isomorphism
$X\setminus 
\bigcup_{i:\sigma(D_i)\neq D_i}
D_i
\to
\Gamma_\sigma'\cap (X\times X)^\sim.$
\qed

\begin{df}\label{dfsig}
Let $X$ be a separated and smooth scheme 
of finite type over $F$,
$D$ be a divisor of $X$ with simple normal crossings
and $U=X\setminus D$ be the complement.
Let $\sigma$ be an admissible automorphism
of $X$ over $F$
satisfying $\sigma(U)=U$
and let $\tilde \Gamma_\sigma\subset
(X\times X)^\sim$ denote
the image of the closed immersion
$\tilde \Gamma_\sigma:
X\setminus 
\bigcup_{i:\sigma(D_i)\neq D_i}D_i\to
(X\times X)^\sim.$
We call the closed subscheme
\begin{equation}
X_{\log}^\sigma=
\Delta_X\cap 
\tilde \Gamma_\sigma=
X\times_{(X\times X)^\sim}
\tilde \Gamma_\sigma
\label{eqsig}
\end{equation}
of $X$ the log $\sigma$-fixed part.
\end{df}

\begin{lm}\label{lmsig}
Let $X$ be a separated and smooth scheme 
of finite type over $F$,
$D$ be a divisor of $X$ with simple normal crossings
and $U=X\setminus D$ be the complement.
Let $\sigma$ be an admissible automorphism
of $X$ over $F$
satisfying $\sigma(U)=U$.

1. The closed subscheme $X_{\log}^\sigma\subset X$
is a closed subscheme of 
the $\sigma$-fixed part $X^\sigma=
X\times_{X\times X
\swarrow\Gamma_\sigma}X$.

2. Let $k\in {\mathbb Z}$
be an integer and assume
$\sigma^k$
is also admissible.
Then, we have an inclusion
$$X_{\log}^\sigma
\subset X_{\log}^{\sigma^k}$$
of closed subschemes.

3. Assume $U^\sigma=\emptyset$
and $\sigma$
is of finite order invertible in $F$.
Then, we have
$$X_{\log}^\sigma
=\emptyset.$$
\end{lm}

\noindent{\it Proof.}
1. Clear from the commutative diagram
$$\begin{CD}
X\setminus 
\bigcup_{i:\sigma(D_i)\neq D_i}D_i
@>{\Gamma_\sigma}>>
(X\times X)^\sim\\
@V{\cap}VV @VVV\\
X
@>{(1,\sigma)}>>X\times X.
\end{CD}$$

2. Since 
$X^\sigma_{\log}=
X^{\sigma^{-1}}_{\log}$
and $X^{\rm id}_{\log}=X$,
we may assume $k\ge 1$.
Let $J_\sigma$
and $J_{\sigma^k}$
be the ideals of $O_X$
defining the closed subschemes
$X^\sigma_{\log}$
and 
$X^{\sigma^k}_{\log}$ respectively.
By 1,
it is sufficient to show
the inclusion 
$J_{\sigma^k,x}\subset J_{\sigma,x}$
of the ideals of $O_{X,x}$
for each $x\in X^\sigma$.
Let $x$ be a point of $X^\sigma$.
The ideal
$J_{\sigma,x}$ is generated by
$\sigma(a)-a$
and
$\sigma(b)/b-1$
for $a\in O_{X,x}$ and
$b\in O_{X,x}\cap j_*O_{U,x}^\times$
where $j:U\to X$ is the open immersion.
Similarly, 
$J_{\sigma^k,x}$ is generated by
$\sigma^k(a)-a$
and
$\sigma^k(b)/b-1$
for $a\in O_{X,x}$ and
$b\in O_{X,x}\cap j_*O_{U,x}^\times$.
Since $\sigma$ is admissible,
we have
$\sigma(b)/b\in O_{X,x}^\times$ for
$b\in O_{X,x}\cap j_*O_{U,x}^\times$.
We have
$\sigma^k(a)-a=
\sum_{i=0}^{k-1}(\sigma(\sigma^i(a))-\sigma^i(a))
\in J_{\sigma,x}$
and 
$\sigma^k(b)/b-1=
\sum_{i=0}^{k-1}(\sigma(\sigma^i(b))/\sigma^i(b)-1)(\sigma^i(b)/b)
\in J_{\sigma,x}$
for $a\in O_{X,x}$ and
$b\in O_{X,x}\cap j_*O_{U,x}^\times$.
Hence, we have
$J_{\sigma^k,x}\subset J_{\sigma,x}$.

3. By 1,
it is sufficient to show
$J_{\sigma,x}=O_{X,x}$
for each closed point $x\in X^\sigma$.
Let $x$ be a closed point in 
$X^\sigma$
and $e$ be the order of $\sigma$.
Since the question is \'etale local,
we may assume $F$ contains a primitive
$e$-th root of unity.
We take a regular system 
$t_1,\ldots, t_d$
of parameters
of $O_{X,x}$
such that $t_1\cdots t_r$
defines $D$ at $x$.
By replacing $t_i$'s if necessary,
we may assume there is
a unique $e$-th root $\zeta_i$
of unity such that
$\sigma(t_i)\equiv \zeta_i t_i\bmod m_x^2$
for each $t_i$.
Replacing $t_i$
by $\sum_{k=1}^e \zeta_i^{-k}\sigma^k(t_i)/e$,
we may assume
$\sigma(t_i)=\zeta_it_i$.
Then, the ideal
$J_{\sigma,x}$
is generated by
$\zeta_i-1$ for $1\le i\le r$
and $(\zeta_i-1)t_i$
for $r<i\le d$.
Since $\zeta_i-1$
is invertible unless $\zeta_i=1$, we have
either $J_{\sigma,x}=O_{X,x}$
or $J_{\sigma,x}=
((\zeta_i-1)t_i,(r<i\le d))$.
By the assumption that
$U^\sigma=\emptyset$,
we have
$J_{\sigma,x}=O_{X,x}$
and the assertion follows.
\qed

\begin{cor}\label{corsig}
Let the notation be as in Lemma \ref{lmsig}.
Assume $\sigma$ is of finite order $e$
and $\sigma^j$ is admissible
for each $j\in{\mathbb Z}$.

1. If $j$ is prime to $e$,
we have
$$X_{\log}^\sigma=
X_{\log}^{\sigma^j}.$$

2. If $U^\sigma=\emptyset$
and if $e$ is not a power of characteristic of $F$,
then we have
$$X_{\log}^\sigma
=\emptyset.$$
\end{cor}

\noindent{\it Proof.}
Clear from Lemma \ref{lmsig}.2 and 3.
\qed

\section{A Lefschetz trace formula for open varieties}

In preliminary subsections
\S\S2.1 and 2.2, we recall some facts on the cycle class map
and
a lemma of Faltings
on the cohomology of the log self product respectively.
In \S2.3, we prove
a Lefschetz trace formula,
Theorem \ref{thmLTF}, for open varieties.

In this section, 
we keep the notation that $F$ denotes a field
and $\ell$ denotes a prime number 
invertible in $F$.

\subsection{Complements on cycle maps}

We recall some facts on cycle maps.
Let $X$ be a smooth scheme over $F$
and $i:Y\to X$
be a closed immersion
of codimension $d$.
Then, 
the cycle class $[Y]\in H^{2d}_Y(X,{\mathbb Z}_\ell(d))$
and the corresponding map
${\mathbb Z}_\ell\to Ri^!{\mathbb Z}_\ell(d)[2d]$
are defined in \cite{cycle}.

\begin{lm}\label{lmcup}
Let $X$ be a smooth scheme over $F$
and $j:U\to X$ be an open immersion.
Let $i:Y\to U$ be a closed immersion
and assume that the composition
$i'=j\circ i:Y\to X$ is also a closed immersion.
Assume that $Y$ is of codimension $d$ in $X$.
Then,
for an integer $q\in {\mathbb Z}$,
the composition
$$\begin{CD}
H^q_c(X,{\mathbb Z}_\ell)@>{i^{\prime *}}>>
H^q_c(Y,{\mathbb Z}_\ell)
@>{i_*}>>
H^{q+2d}_c(U,{\mathbb Z}_\ell(d))
\end{CD}$$
is the cup-product with the image of 
the cycle class 
$[Y]\in H^{2d}_Y(X,{\mathbb Z}_\ell(d))$
by the map
$H^{2d}_Y(X,{\mathbb Z}_\ell(d))=
H^{2d}_Y(X,j_!{\mathbb Z}_\ell(d))\to
H^{2d}(X,j_!{\mathbb Z}_\ell(d))$.
\end{lm}

\noindent{\it Proof.}
The cycle class $[Y]\in H^{2d}_Y(X,{\mathbb Z}_\ell(d))$
defines a map
${\mathbb Z}_\ell\to Ri^{\prime!}{\mathbb Z}_\ell(d)[2d]$.
The push-forward map
$i_*:H^q_c(Y,{\mathbb Z}_\ell)\to
H^{q+2d}_c(U,{\mathbb Z}_\ell(d))$
is the composition 
of the map
$H^q_c(Y,{\mathbb Z}_\ell)\to
H^{q+2d}_{Y!}(X,j_!{\mathbb Z}_\ell(d))$
induced by 
${\mathbb Z}_\ell\to Ri^{\prime!}{\mathbb Z}_\ell(d)[2d]$
with the canonical map
$H^{q+2d}_{Y!}(X,j_!{\mathbb Z}_\ell(d))
\to
H^{q+2d}_c(U,{\mathbb Z}_\ell(d))$
in the notation of \cite{cycle} 1.2.5, 2.3.1.
Hence the assertion follows.
\qed

\begin{lm}\label{lmcompa}
Let $X$ and $Y$ be smooth schemes
purely of dimensions $n$ and $m$ over $F$
and $f:X\to Y$
be a morphism over $F$.
Let $Z$ be a closed subscheme of $Y$
of codimension $d$
and put $W=Z\times_YX$.
Then,
the image of
the cycle class 
$[Z]\in H^{2d}_Z(Y,{\mathbb Z}_\ell(d))$
by the pull-back map
$f^*:
H^{2d}_Z(Y,{\mathbb Z}_\ell(d))\to
H^{2d}_W(X,{\mathbb Z}_\ell(d))$
is equal to the cycle class
$[f^!(Z)]$
of the image $f^!(Z)$ of the Gysin map
$f^!:CH_{m-d}(Z)
\to CH_{n-d}(W)$.
\end{lm}

\noindent
{\it Proof.}
In the case $f:X\to Y$ is smooth,
the assertion is in \cite{cycle} Th\'eor\`eme 2.3.8 (ii).
By decomposing $f:X\to Y$
as the composition of the graph map $X\to X\times Y$
with the projection $X\times Y\to Y$,
we may assume $f:X\to Y$
is a closed immersion.
We prove this case using
the deformation to normal cone.

Let $(Y\times {\mathbf A}^1)'\to Y\times {\mathbf A}^1$
be the blow-up at $X\times \{0\}$ and
let $Y'$ be the complement of the proper
transform of $Y\times \{0\}$ in
$(Y\times {\mathbf A}^1)'$.
Let $Z'$ be 
the proper transform of $Z\times {\mathbf A}^1$
in $Y'$.
The fiber $Y'\times_{{\mathbf A}^1}\{0\}$
at $0$ is naturally identified with
the normal bundle $N=N_{X/Y}$ of $X$ in $Y$
and $Z'\times_{{\mathbf A}^1}\{0\}$
is also identified with the normal cone
$C=C_WZ$ of $W=X\times_YZ$ in $Z$ \cite{fulton} Chapter 5.1.
Let $f':X\times {\mathbf A}^1\to Y'$
denote the immersion
and $g:X\to N$ be the 0-section.
We consider the commutative diagram
$$\begin{CD}
H^{2d}_Z(Y,{\mathbb Z}_\ell(d))@<{1^*}<<
H^{2d}_{Z'}(Y',{\mathbb Z}_\ell(d))@>{0^*}>>
H^{2d}_C(N,{\mathbb Z}_\ell(d))\\
@V{f^*}VV @V{f^{\prime *}}VV @VV{g^*}V\\
H^{2d}_W(X,{\mathbb Z}_\ell(d))@<{1^*}<<
H^{2d}_{W\times {\mathbf A}^1}
(X\times {\mathbf A}^1,{\mathbb Z}_\ell(d))@>{0^*}>>
H^{2d}_W(X,{\mathbb Z}_\ell(d)).
\end{CD}$$
The lower horizontal arrows are the same
and are isomorphisms.
In the upper line,
the images of the cycle class $[Z']$
in the middle are
the cycle classes
$[Z]$ and $[C]$
respectively by \cite{cycle} Th\'eor\`eme 2.3.8 (ii).
Since $f^!(Z)$
is defined as $g^!(C)$ \cite{fulton}
Chapter 6.1 (1),
it is reduced to showing
the equality 
$g^*([C])=[g^!(C)]$.

We put $N_W=N\times_XW$.
Since $C\subset N_W$,
the pull-back $g^*:
H^{2d}_C(N,{\mathbb Z}_\ell(d))\to
H^{2d}_W(X,{\mathbb Z}_\ell(d))$
is the composition
$H^{2d}_C(N,{\mathbb Z}_\ell(d))\to
H^{2d}_{N_W}(N,{\mathbb Z}_\ell(d))
\overset{g^*}\to
H^{2d}_W(X,{\mathbb Z}_\ell(d))$.
Thus it is reduced to showing that the diagram
$$\begin{CD}
CH_{m-d}(N_W)@>{cl}>>
H^{2d}_{N_W}(N,{\mathbb Z}_\ell(d))\\
@V{g^!}VV @VV{g^*}V\\
CH_{n-d}(W)@>>{cl}>
H^{2d}_W(X,{\mathbb Z}_\ell(d))
\end{CD}$$
is commutative.
Let $p:N\to X$ be the projection.
Then the maps $g^!$ and $g^*$
are the inverse of the pull-back map $p^*$.
Hence it is reduced to the case
where $f=p$ is smooth.
\qed

\subsection{Cohomology of the log self products}

We recall a lemma of Faltings on
the cohomology of 
the log self products.
To state it,
we introduce a notation.
Let $Y$ be a smooth scheme over $F$
and $D_1,D_2$ be 
relatively prime divisors of $Y$
such that
the sum $D_1\cup D_2$
has simple normal crossings.
Let
$$
\begin{CD}
Y\setminus (D_1\cup D_2)@>{k_2}>>
Y\setminus D_1\\
@V{k_1}VV@V{j_1}VV\\
Y\setminus D_2@>{j_2}>>Y
\end{CD}
$$
be open immersions.
Let $\ell$ be a prime number invertible in $F$.
Then, the base change map
\begin{equation}
\begin{CD}
j_{1!}Rk_{2*}{\mathbb Z}_\ell
@>>>
Rj_{2*}k_{1!}{\mathbb Z}_\ell
\end{CD}
\label{eqjkY}
\end{equation}
is an isomorphism.
We will identify 
$j_{1!}Rk_{2*}{\mathbb Z}_\ell
=Rj_{2*}k_{1!}{\mathbb Z}_\ell$
by the isomorphism
(\ref{eqjkY}).
We define
$$H^q(Y,D_{1!},D_{2*},
{\mathbb Z}_\ell)=
H^q(Y,D_{2*},D_{1!},
{\mathbb Z}_\ell)$$
to be
$H^q(Y,j_{1!}Rk_{2*}{\mathbb Z}_\ell)
=
H^q(Y,Rj_{2*}k_{1!}{\mathbb Z}_\ell)$.
If $D_1$ or $D_2$ is empty,
we write simply
$H^q(Y,D_{1!},\emptyset_*,
{\mathbb Z}_\ell)
=H^q(Y,D_{1!},
{\mathbb Z}_\ell)$
or
$H^q(Y,D_{2*},\emptyset_!,
{\mathbb Z}_\ell)=
H^q(Y,D_{2*},
{\mathbb Z}_\ell)$
respectively.
With this convention,
we have
$H^q(Y,D_{1!},
{\mathbb Z}_\ell)
=H^q_c(Y-D_1,
{\mathbb Z}_\ell)$,
if $Y$ is proper,
and
$H^q(Y,D_{2*},
{\mathbb Z}_\ell)=
H^q(Y-D_2,
{\mathbb Z}_\ell)$.

Let $X$ be a smooth scheme 
of finite type over a field $F$,
$D$ be a divisor of $X$ with simple normal crossings
and $U=X\setminus D$ be the complement.
We consider a commutative diagram
\begin{equation}
\xymatrix{
 &(X\times X)'\setminus D^{(2)'}\ar[dl]_{j'_2}&&
(X\times X)^\sim\ar[dl]^{k'_2}\ar[ll]_{k'_1}\\
(X\times X)'\ar[dd]_p&&
(X\times X)'\setminus D^{(1)'}
\ar[ll]_{\qquad j'_1}&\\
 &X\times U\ar[dl]_{j_2}\ar[uu]_{e_1}&&
U\times U\ar[dl]^{k_2}
\ar[ll]_{\!\!\!\!\!\!\!\!\!\!\!\!\!\!\!\! k_1}
\ar[uu]_{\tilde j}\\
X\times X&&U\times X.\ar[ll]_{j_1}\ar[uu]_{e_2}&}
\end{equation}
All the arrows
except the log blow-up
$p:(X\times X)'\to X\times X$
are open immersions.
The four faces consisting of open immersions
are Cartesian.
Let $\ell$ be a prime number invertible in $F$.
The canonical maps
$\tilde j_!{\mathbb Z}_\ell
\to 
{\mathbb Z}_\ell
\to 
R\tilde j_*{\mathbb Z}_\ell$
induce maps
\begin{equation}
\begin{CD}
(j'_1 e_2)_!Rk_{2*}{\mathbb Z}_\ell @. @.\\
=j'_{1!}Rk'_{2*}\tilde j_!{\mathbb Z}_\ell
@>>> 
j'_{1!}Rk'_{2*}{\mathbb Z}_\ell @.\\
@.=
Rj'_{2*}k'_{1!}{\mathbb Z}_\ell
@>>>
Rj'_{2*}k'_{1!}R\tilde j_*{\mathbb Z}_\ell\\
@. @. =
R(j'_2e_1)_*k_{1!}{\mathbb Z}_\ell.
\end{CD}
\label{eqj'!}
\end{equation}
The equalities
refer to the identification by
(\ref{eqjkY}).

\begin{lm}\label{lmfal}
{\rm (Faltings)}
Let $X$ be a smooth scheme over $F$, 
$D$ be a divisor of $X$ with simple normal crossings
and $p:(X\times X)'\to X\times X$
be the blow-up {\rm (\ref{eqblup})}.
The maps {\rm (\ref{eqj'!})}
induce isomorphisms
\begin{equation}
\begin{CD}
j_{1!}Rk_{2*}{\mathbb Z}_\ell @. @. \\
=
Rp_*(j'_1 e_2)_!Rk_{2*}{\mathbb Z}_\ell
@>>>
Rp_* j'_{1!}Rk'_{2*}{\mathbb Z}_\ell
@>>>
Rp_* R(j'_2e_1)_*k_{1!}{\mathbb Z}_\ell\\
@. @. =
Rj_{2*}k_{1!}{\mathbb Z}_\ell
\end{CD}\label{eqfal}
\end{equation}
and the composition
is the isomorphism
{\rm (\ref{eqjkY})}.
\end{lm}

For the sake of completeness,
we recall the proof in \cite{fal}.

\noindent{\it Proof.}
Since $j_1=p\circ j'_1\circ e_2$,
$j_2=p\circ j'_2\circ e_1$
and $p$ is proper,
we have
$j_{1!}Rk_{2*}{\mathbb Z}_\ell
=
Rp_*(j'_1 e_2)_!Rk_{2*}{\mathbb Z}_\ell$
and
$Rp_* R(j'_2e_1)_*k_{1!}{\mathbb Z}_\ell
=
Rj_{2*}k_{1!}{\mathbb Z}_\ell$.
It is clear that the composition
is the isomorphism
{\rm (\ref{eqjkY})}.
Thus, it is sufficient to show that the 
first arrow
\begin{equation}
\begin{CD}
Rp_*(j'_1 e_2)_!Rk_{2*}{\mathbb Z}_\ell
@>>>
Rp_* j'_{1!}Rk'_{2*}{\mathbb Z}_\ell
\end{CD}\label{eq!*}
\end{equation}
is an isomorphism.
Since the question is \'etale local on $X\times X$,
it is reduced to the case where
$X={\rm Spec}\ F[T_1,\ldots,T_d]$
and $D$ is defined by
$T_1\cdots T_r=0$.
Further by the K\"unneth formula,
it is reduced to the case where
$X={\mathbf A}^1={\rm Spec}\ F[T]$
and $D$ is defined by
$T=0$.
In this case,
by the proper base change theorem,
the assertion follows from
$H^q({\mathbf A}^1_{\bar F},j_!{\mathbb Z}_\ell)=0$
for $q\in {\mathbb Z}$
where $j:{\mathbf A}^1\setminus\{0\}
\to {\mathbf A}^1$ is the open immersion.
\qed

\begin{cor}\label{corfal}
{\rm (Faltings)}
Let the notation be as in
Lemma \ref{lmfal}.
If $X$ is proper over $F$,
the maps
\begin{eqnarray*}
&&H^q(X_{\bar F}\times X_{\bar F},D^{(1)}_{\bar F!},
D^{(2)}_{\bar F*}, {\mathbb Z}_\ell(d))=
H^q(
(X\times X)'_{\bar F},
(D^{(1)\prime}\cup E)_{\bar F!},D^{(2)\prime}_{\bar F*},
{\mathbb Z}_\ell(d))\\
&\longrightarrow&
H^q(
(X\times X)'_{\bar F},D^{(1)\prime}_{\bar F!},
D^{(2)\prime}_{\bar F*},
{\mathbb Z}_\ell(d))\\
&\longrightarrow&
H^q(
(X\times X)'_{\bar F},D^{(1)\prime}_{\bar F!},
(D^{(2)\prime}\cup E)_{\bar F*},
{\mathbb Z}_\ell(d))
\end{eqnarray*}
are isomorphisms
for $q\in {\mathbb Z}$.
\end{cor}

{\it Proof.}
Clear from Lemma \ref{lmfal}.
\qed

\subsection{A Lefschetz trace formula for open varieties}

Let $F$ be a field,
$X$ be a proper scheme over $F$
and $U$ be a dense open subscheme of $X$.
Let $\Gamma\subset U\times U$ be
a closed subscheme.
Let $p_1,p_2:\Gamma \to U$
denote the compositions of
the closed immersion
$i:\Gamma\to U\times U$ with
the projections $pr_1,pr_2:U\times U\to U$.

\begin{lm}\label{lmproper}
Let $X$ be a proper scheme over $F$.
Let $D$ be a closed subscheme
and $U=X\setminus D\subset X$ be the complement.
Let $\Gamma\subset U\times U$ be
a closed subscheme and
$\overline \Gamma$ be the closure of
$\Gamma$ in $X\times X$.
We put 
$D^{(1)}=D\times X$ and
$D^{(2)}=X\times D$.
Then, the second projection $p_2:\Gamma \to U$
is proper
if and only if
we have the inclusion 
\begin{equation}
\overline \Gamma \cap D^{(1)}
\subset 
\overline \Gamma \cap D^{(2)}
\label{eqproper}
\end{equation}
of the underlying sets.
\end{lm}

\noindent{\it Proof.}
The projection 
$p_2:\Gamma\to U$ is proper
if and only if
$\overline \Gamma\cap (X\times U)=
\overline \Gamma\cap (U\times U)=\Gamma$.
Taking the complement,
it is equivalent to
$\overline \Gamma\cap D^{(2)}=
\overline \Gamma\cap (D^{(1)}\cup D^{(2)})$.
It is further equivalent to
$\bar\Gamma\cap D^{(1)}\subset
\overline \Gamma\cap D^{(2)}$.
\qed

In the following,
we assume that 
$U$ is smooth purely of dimension $d$,
that $\Gamma$ is purely of dimension $d$
and that $p_2:\Gamma\to U$ is proper.
For a prime number $\ell$
invertible in $F$,
we define an endomorphism $\Gamma^*$
of $H^q_c(V_{\bar F},{\mathbb Z}_\ell)$
to be $p_{1*}\circ p_2^*$
and consider
the alternating sum
$${\rm Tr}(\Gamma^*:
H^*_c(U_{\bar F},{\mathbb Z}_\ell))=
\sum_{q=0}^{2d}(-1)^q{\rm Tr}(\Gamma^*:
H^q_c(U_{\bar F},{\mathbb Z}_\ell)).$$
Since $p_2$ is assumed proper,
the pull-back
$p_2^*:
H^q_c(U_{\bar F},{\mathbb Z}_\ell)
\to
H^q_c(\Gamma_{\bar F},{\mathbb Z}_\ell)$
is defined.
We briefly recall the definition of the push-forward map 
$p_{1*}:
H^q_c(\Gamma_{\bar F},{\mathbb Z}_\ell)\to
H^q_c(U_{\bar F},{\mathbb Z}_\ell)$.
Let $f:U\to  {\rm Spec}\ F$ and
$g:\Gamma\to {\rm Spec}\ F$
denote the structural maps.
Then the trace map
$Rg_!{\mathbb Z}_\ell(d)[2d]\to {\mathbb Z}_\ell$
induces the cycle class map
${\mathbb Z}_\ell(d)[2d]\to Rg^!{\mathbb Z}_\ell$.
Since $U$ is smooth of dimension $d$,
the cycle class map for $U$ induces an isomorphism
$Rp_1^!{\mathbb Z}_\ell(d)[2d]\to 
Rp_1^!Rf^!{\mathbb Z}_\ell\to 
Rg^!{\mathbb Z}_\ell$.
Thus, we obtain a canonical map
${\mathbb Z}_\ell
\to
Rp_1^!{\mathbb Z}_\ell$
and hence 
$Rp_{1!}{\mathbb Z}_\ell
\to
{\mathbb Z}_\ell$
by adjunction.
The map $Rp_{1!}{\mathbb Z}_\ell
\to
{\mathbb Z}_\ell$
induces the push-forward map 
$p_{1*}:
H^q_c(\Gamma_{\bar F},{\mathbb Z}_\ell)\to
H^q_c(U_{\bar F},{\mathbb Z}_\ell)$.

We give another description of the map
$\Gamma^*=p_{1*}\circ p_2^*$
using the cycle class of $\Gamma$.
We put $H^q_{!,*}(U_{\bar F}\times U_{\bar F},{\mathbb Z}_\ell(d))=
H^q(X_{\bar F}\times U_{\bar F},
(j\times{\rm id})_!{\mathbb Z}_\ell(d))$.
By the assumption that
$p_2:\Gamma\to U$ is proper,
$\Gamma$ is closed in 
$X\times U$ and hence
the canonical maps
$$
H^{2d}_{\Gamma}(X\times U,
(j\times{\rm id})_!{\mathbb Z}_\ell(d))
\to
H^{2d}_{\Gamma}(X\times U,{\mathbb Z}_\ell(d))
\to
H^{2d}_{\Gamma}(U\times U,{\mathbb Z}_\ell(d))$$
are isomorphisms.
Thus the cycle class
$[\Gamma]\in
H^{2d}_{\Gamma}(U\times U,{\mathbb Z}_\ell(d))$
defines a class
$[\Gamma]\in
H^{2d}(X_{\bar F}\times U_{\bar F},(j\times{\rm id})_!{\mathbb Z}_\ell(d))=
H^{2d}_{!,*}(U_{\bar F}\times U_{\bar F},{\mathbb Z}_\ell(d))$.
By the K\"unneth formula
and Poincar\'e duality,
we have
canonical isomorphisms
$$\begin{CD}
\bigoplus_q
H^q_c(U_{\bar F},{\mathbb Q}_\ell)\otimes
H^{2d-q}(U_{\bar F},{\mathbb Q}_\ell(d))
@>>>
H^{2d}_{!,*}(U_{\bar F}\times U_{\bar F},
{\mathbb Q}_\ell(d))
\\
@VVV @.\\
\bigoplus_q
H^q_c(U_{\bar F},{\mathbb Q}_\ell)\otimes
Hom(H^q_c(U_{\bar F},{\mathbb Q}_\ell),
{\mathbb Q}_\ell)
@>>>\prod_{q=0}^{2d}
{\rm End}\ H^q_c(U_{\bar F},{\mathbb Q}_\ell).
\end{CD}$$

\begin{lm}\label{lmcyccor}
Let $\Gamma\subset U\times U$ be 
a closed subscheme of dimension $d$.
Assume that $p_2:\Gamma\to U$ is proper.
Then, by the canonical isomorphism
$H^{2d}_{!,*}(U_{\bar F}\times U_{\bar F},{\mathbb Q}_\ell(d))\to 
\prod_{q=0}^{2d}
{\rm End}\ H^q_c(U_{\bar F},{\mathbb Q}_\ell)$,
the image of the cycle class $[\Gamma]$ 
is $\Gamma^*$.
\end{lm}

\noindent
{\it Proof.}
It is sufficient to show the equality
$$
p_{1*}p_2^*\alpha=
pr_{1*}([\Gamma]\cup pr_2^*\alpha)$$
in $H^q_c(U_{\bar F},{\mathbb Q}_\ell)$
for an arbitrary integer $q\in {\mathbb Z}$
and $\alpha\in H^q_c(U_{\bar F},{\mathbb Q}_\ell)$.
Let $i:\Gamma\to U\times U$ be the immersion
and $i':\Gamma\to U\times U\to X\times U$ be
the composition.
Since 
$p_{1*}p_2^*\alpha=
pr_{1*}(i_*i^{\prime*}pr_2^*\alpha)$,
it is reduced to showing the equality
$$i_*i^{\prime*}\beta =
[\Gamma]\cup \beta$$
in $H^q_c(U_{\bar F}\times U_{\bar F},{\mathbb Q}_\ell)$
for $\beta \in H^q_c(X_{\bar F}\times U_{\bar F},{\mathbb Q}_\ell)$.
By Lemma \ref{lmcup},
the class 
$i_*i^{\prime*}\beta$
is the product with the class of $\Gamma$.
Thus the assertion follows.
\qed

\begin{lm}\label{lmtr}
Let $U$ and $V$ be connected
separated smooth schemes
of finite type purely of dimension $d$
over $F$.
Let $g:U\to V$
be a proper and generically finite morphism 
of constant degree $[U:V]$
over $F$.
Then, for
$\Gamma^*\in 
H^{2d}_{!,*}(V_{\bar F}\times V_{\bar F},{\mathbb Q}_\ell)=
\prod_{q=0}^{2d}
{\rm End}\ H^q_c(V_{\bar F},{\mathbb Q}_\ell)$,
we have
\begin{equation}
{\rm Tr}(\Gamma^*:H^*_c(V_{\bar F},{\mathbb Q}_\ell))=
\frac 1{[U:V]}
{\rm Tr}(((g\times g)^*\Gamma)^*:H^*_c(U_{\bar F},{\mathbb Q}_\ell)).
\label{eqWV}
\end{equation}
\end{lm}

\noindent{\it Proof.}
Since 
$g^*:H^*_c(V_{\bar F},{\mathbb Q}_\ell)
\to
H^*_c(U_{\bar F},{\mathbb Q}_\ell)$
is injective
and $g_*\circ g^*$
is the multiplication by $[U:V]$,
it is sufficient to show that
$((g\times g)^*\Gamma)^*$
is the composition
$g^*\circ \Gamma^*\circ g_*$.
In other words,
it suffices to show the equality
$$
pr_{1*}(((g\times g)^*[\Gamma]\cup pr_2^*\alpha)=
g^* (pr_{1*}([\Gamma]\cup pr_2^*g_*\alpha))$$
for $q\in {\mathbb Z}$ and
$\alpha \in H^q_c(U_{\bar F},{\mathbb Q}_\ell)$.
In the commutative diagram
$$\begin{CD}
U\times U&@>{pr_2}>> &U\\
@V{1\times g}VV @. @VVgV\\
U\times V @>{g\times 1}>>  V\times V@>>{pr_2}> V\\
@V{pr_1}VV @VV{pr_1}V @.\\
U @>g>> V, @.
\end{CD}
$$
$\alpha$ lives on $U$ in the northeast
and $\Gamma$ lives on $V\times V$.
Thus, by the projection formula,
we compute
\begin{eqnarray*}
&&pr_{1*}((g\times g)^*[\Gamma]\cup pr_2^*\alpha)=
pr_{1*}((1\times g)^*(g\times 1)^*[\Gamma]\cup pr_2^*\alpha)\\
&=&
pr_{1*}((g\times 1)^*[\Gamma]\cup (1\times g)_*pr_2^*\alpha)=
pr_{1*}((g\times 1)^*([\Gamma]\cup pr_2^* g_*\alpha))\\
&=&
g^* pr_{1*}([\Gamma]\cup pr_2^*g_*\alpha).
\end{eqnarray*}
\qed

We prove a Lefschetz trace formula
for open varieties.

\begin{thm}\label{thmLTF}
Let $X$ be a proper and smooth scheme
purely of dimension $d$
over a field $F$
and $U$ be the complement
of a divisor $D$ with simple normal crossings.
Let $\Gamma\subset U\times U$ be 
a closed subscheme 
purely of dimension $d$.
Let $D^{(1)\prime },D^{(2)\prime }\subset (X\times X)'$
denote the proper transforms of
$D^{(1)},D^{(2)}$
respectively
and let $\overline \Gamma'$ be the closure of $\Gamma$
in $(X\times X)'$.
We assume that 
we have an inclusion 
\begin{equation}
\overline \Gamma' \cap D^{(1)\prime}
\subset 
\overline \Gamma' \cap D^{(2)\prime}
\label{eqthm}
\end{equation}
of the underlying sets.

Then, 
the map
$p_2:\Gamma\to U$ is proper and
we have an equality
\begin{equation}
{\rm Tr}(\Gamma^*:
H^*_c(U_{\bar F},{\mathbb Q}_\ell))={\rm deg}\
(\overline \Gamma',\Delta'_X)_{(X\times X)'}.
\label{eqLTF}
\end{equation}
The right hand side is the
intersection product 
in $(X\times X)'$ 
of the closure $\overline \Gamma'$ 
with the image $\Delta'_X$
of the log diagonal closed immersion
$\Delta':
X\to (X\times X)'$.
\end{thm}

\noindent
{\it Proof.}
First, we show
the map
$p_2:\Gamma\to U$ is proper.
By the assumption (\ref{eqthm})
$\overline \Gamma' \cap D^{(1)\prime}
\subset 
\overline \Gamma' \cap D^{(2)\prime}$,
we have
$\overline \Gamma' \cap (D^{(1)\prime}\cup E)
\subset 
\overline \Gamma' \cap (D^{(2)\prime}\cup E)$.
Hence we have (\ref{eqproper})
$\overline \Gamma \cap D^{(1)}
\subset 
\overline \Gamma \cap D^{(2)}$
and the assertion follows by Lemma \ref{lmproper}.

Since the restriction of
$j_{1!}Rk_{2*}{\mathbb Z}_\ell(d)$
on the diagonal $X\subset X\times X$
is
$j_!{\mathbb Z}_\ell(d)$,
the pull-back map 
$$\begin{CD}
\Delta^*:&
H^{2d}_{!,*}(U_{\bar F}\times U_{\bar F},{\mathbb Z}_\ell(d))
@>>>
H^{2d}_c(U_{\bar F},{\mathbb Z}_\ell(d))
\\
&=
H^{2d}(X_{\bar F}\times X_{\bar F},
j_{1!}Rk_{2*}{\mathbb Z}_\ell(d))
@.=H^{2d}(X_{\bar F},j_!{\mathbb Z}_\ell(d))
\end{CD}$$
by the diagonal is defined.
Then, 
by Lemma \ref{lmcyccor}
and by the standard argument (cf. Proposition 3.3 \cite{cycle})
in the proof of Lefschetz trace formula,
we have
\begin{equation}
{\rm Tr}(\Gamma^*:
H^*_c(U_{\bar F},{\mathbb Q}_\ell))=
{\rm Tr}(\Delta^*([\Gamma])).
\label{eqtr}
\end{equation}

In the notation introduced in the beginning of \S2.2,
we have
$H^{2d}_{!,*}(U_{\bar F}\times U_{\bar F},{\mathbb Z}_\ell(d))
=
H^{2d}(X_{\bar F}\times X_{\bar F},D^{(1)}_{\bar F!},
D^{(2)}_{\bar F*}, {\mathbb Z}_\ell(d))$
and
$H^{2d}_c(U_{\bar F},{\mathbb Z}_\ell(d))=
H^{2d}(X_{\bar F},D_{\bar F!},{\mathbb Z}_\ell)$.
The canonical map
$(X\times X)'\to X\times X$
induces an isomorphism
$H^q(X_{\bar F}\times X_{\bar F},D^{(1)}_{\bar F!},
D^{(2)}_{\bar F*}, {\mathbb Z}_\ell(d))
\to
H^q((X\times X)'_{\bar F},(D^{(1)\prime}\cup E)_{\bar F!},
D^{(2)\prime}_{\bar F*}, {\mathbb Z}_\ell(d))$.
Thus the composition
\begin{eqnarray}
&&H^{2d}_{!,*}(U_{\bar F}\times U_{\bar F},{\mathbb Z}_\ell(d))
=H^q(X_{\bar F}\times X_{\bar F},D^{(1)}_{\bar F!},
D^{(2)}_{\bar F*}, {\mathbb Z}_\ell(d))
\nonumber\\
&\longrightarrow&
H^q((X\times X)'_{\bar F},(D^{(1)\prime}\cup E)_{\bar F!},
D^{(2)\prime}_{\bar F*}, {\mathbb Z}_\ell(d))
\label{eq!}\\
&\longrightarrow&
H^q((X\times X)'_{\bar F},D^{(1)\prime}_{\bar F!},
D^{(2)\prime}_{\bar F*}, {\mathbb Z}_\ell(d))
\nonumber
\end{eqnarray}
is an isomorphism by Corollary \ref{corfal}.

We put $\Gamma'=
\overline \Gamma'\setminus
\overline \Gamma'\cap D^{(2)\prime}$.
By the assumption (\ref{eqthm}),
we have $\Gamma'\cap D^{(1)\prime}=\emptyset$.
Thus the cycle class $[\Gamma']\in 
H^{2d}((X\times X)'_{\bar F},D^{(1)\prime}_{\bar F!},
D^{(2)\prime}_{\bar F*},{\mathbb Z}_\ell(d))$
is defined.
We show that the arrow (\ref{eq!})
sends $[\Gamma]$ to $[\Gamma']$.
By Corollary \ref{corfal},
the map
$$H^{2d}((X\times X)'_{\bar F},D^{(1)\prime}_{\bar F!},
D^{(2)\prime}_{\bar F*},{\mathbb Z}_\ell(d))
\to
H^{2d}((X\times X)'_{\bar F},D^{(1)\prime}_{\bar F!},
(E\cup D^{(2)\prime})_{\bar F*},{\mathbb Z}_\ell(d))$$
is an isomorphism.
By this isomorphism,
both
$[\Gamma']$ and the image of $[\Gamma]$
are sent to $[\Gamma]$.
Hence the arrow (\ref{eq!})
sends $[\Gamma]$ to $[\Gamma']$.

Since 
$\Delta'_X
\cap D^{(2)\prime}=\emptyset$,
the map
$\Delta^{\prime *}:
H^q((X\times X)'_{\bar F},
D^{(2)\prime}_{\bar F*},
{\mathbb Z}_\ell(d))
\to
H^q(X_{\bar F},{\mathbb Z}_\ell(d))$
is defined.
We consider the commutative diagram
\begin{equation}
\begin{CD}
[\Gamma]&\in &
H^{2d}_{!,*}(U_{\bar F}\times U_{\bar F},
{\mathbb Z}_\ell(d))
@>{\Delta^*}>>
H^{2d}_c(U_{\bar F},{\mathbb Z}_\ell(d))\\
\downarrow&&@V{(\ref{eq!})}VV @VVV\\
[\Gamma']&\in &
H^{2d}((X\times X)'_{\bar F},
D^{(1)\prime}_{\bar F!},D^{(2)\prime}_{\bar F*},{\mathbb Z}_\ell(d))
@>{\Delta^{\prime*}}>>
H^{2d}(X_{\bar F},{\mathbb Z}_\ell(d))\\
\downarrow&&@VVV @|\\
[\Gamma']&\in &
H^{2d}((X\times X)'_{\bar F},D^{(2)\prime}_{\bar F*},{\mathbb Z}_\ell(d))
@>{\Delta^{\prime*}}>>
H^{2d}(X_{\bar F},{\mathbb Z}_\ell(d))\\
\uparrow&&@AAA @|\\
[\overline \Gamma']&\in &
H^{2d}((X\times X)'_{\bar F},{\mathbb Z}_\ell(d))
@>{\Delta^{\prime*}}>>
H^{2d}(X_{\bar F},{\mathbb Z}_\ell(d)).
\end{CD}
\end{equation}
As we have shown above,
the arrow (\ref{eq!})
sends $[\Gamma]$ to $[\Gamma']$.
Since the middle and the lower
left vertical arrows send
$[\Gamma']$ and $[\overline \Gamma']$
to $[\Gamma']$ respectively,
we have
\begin{equation}
{\rm Tr}(\Delta^*([\Gamma]))=
{\rm Tr}(\Delta^{\prime *}([\overline \Gamma'])).
\label{eqtr'}
\end{equation}
Since
$${\rm Tr}(\Delta^{\prime *}([\overline \Gamma']))
={\rm deg}\
(\overline \Gamma',\Delta'_X)_{(X\times X)'},$$
the assertion follows from
the equalities (\ref{eqtr}) and 
(\ref{eqtr'}).
\qed

\begin{rmk}\label{rmkLTF}
In Theorem \ref{thmLTF},
we can not replace the assumption {\rm (\ref{eqthm})}
$\overline \Gamma' \cap D^{(1)\prime}
\subset 
\overline \Gamma' \cap D^{(2)\prime}$
by a weaker assumption {\rm (\ref{eqproper})}
$\overline \Gamma \cap D^{(1)}
\subset
\overline \Gamma \cap D^{(2)}$
as the following example shows.
Let $X={\mathbf P}^1$,
$U={\mathbf A}^1$,
and $n\ge1$ be an integer.
Let $f:U\to U$ be
the $n$-th power map
and $\Gamma\subset U\times U$
be the transpose 
$\Gamma=\{(x,y)\in U\times U|
x=y^n\}$
of the graph of $f$.
Then, we have
${\rm Tr}(\Gamma^*:H^*_c(U_{\bar F},{\mathbb Z}_\ell))=
{\rm Tr}(f_*:H^2_c(U_{\bar F},{\mathbb Z}_\ell))=1$
while
$(\Gamma,\Delta)_{(X\times X)'}=n$.
\end{rmk}

One can deduce a part of
a conjecture of Deligne
from Theorem \ref{thmLTF} as follows.
The conjecture of Deligne itself is 
proved assuming
resolution of singularities
by Pink in \cite{Pink}
and proved unconditionally 
by Fujiwara in \cite{Fuji}
using rigid geometry.
In the proof below,
we will not
use rigid geometry or
assume resolution of singularities.

We introduce some notations
assuming
$F$ is a finite field.
For a scheme over $F$,
let $Fr$ denote the Frobenius endomorphism
over $F$.
Let $U$ be a separated smooth scheme
of finite type of pure dimension $d$
over $F$.
Let $\Gamma\subset U\times U$
be a closed subscheme of dimension $d$
and assume the composition
$p_2:\Gamma\to U$
with the projection is proper.
For an integer $n\ge 0$
and a prime number $\ell$
different from the characteristic of $F$,
we consider the alternating sum
${\rm Tr}(Fr_F^{*n}\Gamma^*:
H^*_c(U_{\bar F},{\mathbb Q}_\ell))$.
Let $i_n:\Gamma\to U\times U$
be the composition of
the immersion $i:\Gamma\to U\times U$
with the endomorphism
$1\times Fr^n$ of $U\times U$.
Let $\Gamma_n$ denote the scheme $\Gamma$
regarded as a scheme over
$U\times U$ by $i_n$.
If the fiber product
$\Gamma_n\times_{U\times U}\Delta_U$
is proper over $F$,
the degree of the intersection product
$(\Gamma_n,\Delta_U)_{U\times U}\in 
CH_0(\Gamma_n\times_{U\times U}\Delta_U)$
is defined.

\begin{pr}\label{prpink}
{\rm (cf. \cite{Fuji}, \cite{Pink})}
Let $U$ be a separated smooth scheme
of finite type of pure dimension $d$
over a field $F$
and $\ell$ be a prime number
different from the characteristic of $F$.
Let $\Gamma\subset U\times U$
be a closed subscheme of dimension $d$.
Assume the composition
$p_2:\Gamma\to U$
with the projection 
is proper.
Then, we have the following.

1. The alternating sum
${\rm Tr}(\Gamma^*:
H^*_c(U_{\bar F},{\mathbb Q}_\ell))$
is in ${\mathbb Z}[\frac 1p]$
and is independent of $\ell$
invertible in $F$.

2. Assume $F$ is a finite field.
Then, there exists an integer $n_0\ge 0$
satisfying the following property.

For an integer $n\ge n_0$, the fiber product
$\Gamma_n\times_{U\times U}\Delta_U$
is proper over $F$
and we have
\begin{equation}
{\rm Tr}(Fr_F^{*n}\Gamma^*:
H^*_c(U_{\bar F},{\mathbb Q}_\ell))
=\deg(\Gamma_n,\Delta_U)_{U\times U}.
\label{eqpink}
\end{equation}
\end{pr}

{\it Proof.}
1. It is reduced to 2 by a standard argument
using specialization.

2. By the main result of de Jong \cite{dJ}
and Lemma \ref{lmtr},
we may assume that there exists
a proper smooth scheme $X$
containing 
$U$ as the complement of a divisor with 
simple normal crossings.
We will derive Proposition 
from Theorem \ref{thmLTF}
using the following Lemma.

\begin{lm}\label{lmpink}
Let $X$ be a proper smooth scheme
over a finite field $F$ of order $q$
and $D\subset X$
be a divisor with 
simple normal crossings.
Let $U=X\setminus D$ be the complement and
let $\Gamma\subset U\times U$
be an integral closed subscheme.
Assume $p_2:\Gamma\to U$ is proper.

Then, there exists an integer
$n_0\ge 0$ such that,
for all $n\ge n_0$,
the closure $\overline{i_n(\Gamma_n)}
\subset
(X\times X)'$ of the
image $i_n(\Gamma_n)\subset
U\times U$ satisfies the inclusion
\begin{equation}
\overline{i_n(\Gamma_n)}
\setminus i_n(\Gamma_n)
\subset D^{(2)\prime}.
\label{eqthmp}
\end{equation}
\end{lm}

{\it Proof.}
Let $\overline \Gamma
\subset (X\times X)'$
be the closure of $\Gamma$.
By the main result of de Jong \cite{dJ},
there exist
a proper smooth integral scheme
$Z$ of dimension $d$,
a proper map $Z\to \overline \Gamma$
over $F$
such that the inverse image
$W=Z\times_{\overline \Gamma}\Gamma$
is the complement of
a divisor $B$ with simple normal crossings.
Let $Z'\to Z$ be the blow-up
associated to the subdivision by baricenters
and $B'=Z'\setminus W$ be the complement.

Let $\bar r_1,\bar r_2:Z'\to X$
be the compositions with the projections.
Let $D_i,\ (i\in I)$
be the irreducible components of $D$
and
$B'_j,\ (j\in J)$
be the irreducible components of $B'$.
We put $\bar r_1^*D_i=\sum_{j\in J}e_{ij}^{(1)}B'_j$
and $\bar r_2^*D_i=\sum_{j\in J}e_{ij}^{(2)}B'_j$
for $i\in I$.
By the assumption $p_2:\Gamma\to U$
is proper,
the composition $r_2:W\to U$
is proper and hence
the support of $\bar r_2^*D
=\sum_{j\in J}(\sum_{i\in I}e_{ij}^{(2)})B'_j$
equals $B'$.
In other words,
for every $j\in J$,
there exists
an index $i\in I$
such that $e_{ij}^{(2)}>0$.

Let $J_0\subset J$
be the subset 
$\{j\in J|B'_j$
is the proper transform of
an irreducible component of $B\}$.
Then, if 
$B'_j\cap B'_{j'}\neq\emptyset$
and if $j\in J_0$,
we have
$e_{ij}^{(2)}\le
e_{ij'}^{(2)}$.
Hence, if $e_{ij}^{(2)}=0$
and $e_{ij'}^{(2)}>0$
for $B'_j\cap B'_{j'}\neq\emptyset$,
then we have $j\in J_0$.

We show that,
for every $z\in B'$,
there exists
an index $i\in I$
such that $e_{ij}^{(2)}>0$
for all $B'_j\ni z$.
We prove this by contradiction.
Assume there exists
$z\in B'$
such that, for every $i\in I$,
there exists a component
$B'_j\ni z$ such that 
$e_{ij}^{(2)}=0$.
First, we show that
there exists
an element $j_0\in J_0$
such that $z\in B'_{j_0}$.
Let $B'_j$ be a component
containing $z$.
Then, as we have seen above,
there exists an index $i\in I$
such that $e^{(2)}_{ij}>0$.
By the hypothesis,
we also have an index $j_0\in J$
such that $z\in B'_{j_0}$ and
$e^{(2)}_{ij_0}=0$.
Since $z\in B'_{j_0}\cap B'_j$,
we have $j_0\in J_0$.
We show 
$e^{(2)}_{ij_0}=0$
for every $i\in I$,
to get a contradiction.
For $i\in I$,
by the hypothesis,
there exists 
$B'_j\ni z$ such that
$e^{(2)}_{ij}=0$.
Since $z\in B'_{j_0}\cap B'_j\neq
\emptyset$,
we have
$0=e^{(2)}_{ij}
\ge
e^{(2)}_{ij_0}\ge 0$.
Thus we get a contradiction.

We take $n_0\ge 0$
such that $q^{n_0}> \max_{i\in I,j\in J}
e_{ij}^{(1)}$.
Then, for every $z\in B'$,
there exists an index $i\in I$
such that $q^{n_0}e_{ij}^{(2)}>
e_{ij}^{(1)}$
for all $B'_j\ni z$.
Namely,
we have a strict inequality
\begin{equation}
q^{n_0}\bar r_2^*D_i>
\bar r_1^*D_i
\label{eqD>}
\end{equation}
of germs of Cartier divisors
at $z$.

We show the inclusion
(\ref{eqthmp}) for $n\ge n_0$.
We consider the product
$\bar i_n:W \to (X\times X)'\times Z'$
of the composition
$W\to \Gamma$ with $i_n:\Gamma_n\to 
U\times U\subset (X\times X)'$
and the inclusion $W\to Z'$.
Let $Z_n$ be the closure of the image of
the immersion
$\bar i_n:W \to (X\times X)'\times Z'$
with the reduced scheme structure.
Let $\bar r_n:Z_n\to (X\times X)'$ and
$f_n:Z_n\to Z'$ be the projections.
Further, let $\bar r_{1,n},\bar r_{2,n}:
Z_n\to X$ be the compositions
of $\bar r_n$ with the projections.
Then, since
$W\subset Z_n$ is dense,
the diagram
$$\begin{CD}
Z_n @>{f_n}>> Z'\\
@V{\bar r_{1,n}\times\bar r_{2,n}}VV 
@VV{\bar r_1\times\bar r_2}V\\
X\times X@<{1\times Fr^n}<<
X\times X
\end{CD}$$
is commutative. 
Thus, we have equalities
$\bar r_{1,n}^*D_i=
f_n^*\bar r_1^*D_i$
and
$\bar r_{2,n}^*D_i=
q^nf_n^*\bar r_2^*D_i$
of Cartier divisors on $Z_n$
for each $i\in I$.

Since $W\to \Gamma$ is proper and surjective,
we have
$\overline{i_n(\Gamma_n)}
\setminus i_n(\Gamma_n)
=\bar r_n(Z_n\setminus W)$.
For every point
$z\in Z_n\setminus W$,
there exists an index $i\in I$
satisfying a strict inequality 
$$\bar r_{2,n}^*D_i=
q^nf_n^*\bar r_1^*D_i>
f_n^*\bar r_1^*D_i=
\bar r_{1,n}^*D_i$$
of germs of Cartier divisors
at $z$
by (\ref{eqD>}).
Namely, we have
$z\in \bar r_n^{-1}(X\times D_i)'$.
Thus, we have 
$\bar r_n(Z_n\setminus W)
\subset D^{(2)\prime}
=\bigcup_{i\in I}
(X\times D_i)'$
and the assertion follows.
\qed

We complete the proof of Proposition
\ref{prpink}.
Take a proper scheme $\overline{\Gamma_n}$
over $F$
containing $\Gamma$
as a dense open subscheme
and a map $\bar i_n:\overline{\Gamma_n}
\to (X\times X)'$
extending the map
$i_n:\Gamma_n
\to U\times U$.
The intersection of the
log diagonal $\Delta'_X\subset (X\times X)'$
with $D^{(2)\prime}$ is empty.
Hence by the inclusion (\ref{eqthmp})
in Lemma \ref{lmpink},
the intersection
$\overline{i_n(\Gamma_n)}
\cap \Delta'_X$
with the log diagonal
equals 
$i_n(\Gamma_n)
\cap \Delta_U$.
Hence the fiber product
$\Gamma_n\times_{U\times U}\Delta_U
=
\overline{\Gamma_n}\times_{(X\times X)'}\Delta'_X$
is proper over $F$
and we have
$(\overline{\Gamma_n},
\Delta'_X)_{(X\times X)'}
=(\Gamma_n,\Delta_U)_{U\times U}.$

Also by the inclusion (\ref{eqthmp})
in Lemma \ref{lmpink},
the assumption (\ref{eqthm})
of Theorem \ref{thmLTF}
is satisfied for the support of the cycle
$\bar i_{n*}(\overline \Gamma_n)$.
Thus, by Theorem \ref{thmLTF},
we have
${\rm Tr}(Fr_F^{*n}\Gamma^*:
H^*_c(U_{\bar F},{\mathbb Q}_\ell))
=\deg(\overline{\Gamma_n},
\Delta'_X)_{(X\times X)'}=
\deg(\Gamma_n,\Delta_U)_{U\times U}.$
\qed

\section{Intersection product
with the log diagonal and
a trace formula}

We introduce the target group
$CH_0(\overline V\setminus V)$
of the map (\ref{eqmap}) in \S3.1.
We define 
the map (\ref{eqmap})
and prove the trace formula
(\ref{eqGal0})
in \S3.2.
We establish elementary properties
of the map (\ref{eqmap})
in \S3.3.
We define and compute the wild different
of a covering and
the log Lefschetz class of an automorphism
using the map (\ref{eqmap})
in \S3.4.

In this section,
$F$ denotes a perfect field
and $f:V\to U$
is a finite \'etale morphism
of separated and smooth
schemes of finite type
purely of dimension $d$ over $F$.

\subsection{Chow group of 0-cycles
on the boundary}

In this subsection,
we introduce the target group
$CH_0(\overline V\setminus V)$
of the map (\ref{eqmap}).

\begin{df}\label{dfCHV}
Let $V$ be a separated smooth scheme
of finite type
over a field $F$.

1.
Let ${\cal C}_V$
be the following category.
An object of ${\cal C}_V$
is a proper scheme $Y$ over $F$
containing $V$ as a dense open subscheme.
A morphism $Y'\to Y$ in ${\cal C}_V$
is a morphism $Y'\to Y$ over $F$
inducing the identity on $V$.

Let ${\cal C}_V^{\rm sm}$
be the full subcategory of ${\cal C}_V$
consisting of smooth objects.
Let
${\cal C}_V^{\rm sm,0}$
be the full subcategory of ${\cal C}_V$
consisting of smooth objects $Y$
such that $V$ is the complement
of a divisor with simple normal crossings.

2.
We put
\begin{equation}
\begin{CD}
CH_0(\overline V\setminus V)=
\varprojlim_{{\cal C}_V}CH_0(Y\setminus V).
\end{CD}
\label{eqchv}
\end{equation}
The transitions maps
are proper push-forwards.
Let 
\begin{equation}
\begin{CD}
{\rm deg}:
CH_0(\overline V\setminus V)
@>>> {\mathbb Z}
\end{CD}
\label{eqdeg}
\end{equation}
be the limit of the degree maps
$CH_0(Y\setminus V)\to {\mathbb Z}$.
\end{df}

Recall that we assume $F$ is perfect.
The resolution of singularities
means that the full subcategory
${\cal C}_V^{\rm sm}$ is cofinal in 
${\cal C}_V$.
A strong form of the resolution of singularities
means that 
${\cal C}_V^{\rm sm,0}$ is cofinal in 
${\cal C}_V$.
Thus, it is known that 
${\cal C}_V^{\rm sm,0}$ is cofinal in 
${\cal C}_V$
if dimension $V$ is at most $2$.
More precisely,
if dimension is at most $2$,
we have a strong form of equivariant
resolution of singularities as follows.

\begin{lm}\label{lmalter2}
Let $V$ be a separated smooth scheme
of finite type of dimension $\le 2$
over a perfect field $F$
and $G$ be a finite group
of automorphisms of $V$ over $F$.

Then the full subcategory of
${\cal C}_V^{\rm sm,0}$
consisting of $Y$
with an admissible action of $G$
extending that on $V$
is cofinal in 
${\cal C}_V$.
\end{lm}

\noindent{\it Proof.}
Let $Y_0$ be an object of
${\cal C}_V$.
Let $Y_1$ be the closure of
the image of the map
$V\to \prod_{\sigma\in G}V
\subset \prod_{\sigma\in G}Y_0$
sending $v$ to $(\sigma(v))_{\sigma\in G}$.
Let $Y_2$ be the minimal resolution
of the normalization of $Y_1$.
By blowing-up $Y_2$ successively at the closed points
where the complement $Y_2\setminus V$
does not have simple normal crossing,
we obtain $Y_3$ in $C_V^{\rm sm,0}$
with an action of $G$.
The action of $G$ on the blow-up $Y$ of $Y_3$
associated to the subdivision by baricenters
is admissible by Lemma \ref{lmdec}.2
\qed

Let $Y$ be a separated scheme
of finite type over $F$
containing $V$ as a dense open subscheme.
Then there exists a unique map
$CH_0(\overline V\setminus V)
\to
CH_0(Y\setminus V)$
satisfying the following property.
Let $Y'$ be an object of ${\cal C}_V$
containing $Y$ as a dense open subscheme.
Then it is the same as the composition
of the projection
$CH_0(\overline V\setminus V)
\to
CH_0(Y'\setminus V)$
and the restriction
$CH_0(Y'\setminus V)
\to
CH_0(Y\setminus V)$.

Let $f:V\to U$ be a finite flat morphism
of smooth schemes over $F$.
The push-forward maps induces
a map
$f_*:CH_0(\overline V\setminus V)
\to
CH_0(\overline U\setminus U)$.
The flat pull-back map
$f^*:CH_0(\overline U\setminus U)
\to
CH_0(\overline V\setminus V)$
is defined as follows.

\begin{lm}\label{lmGR}
Let $f:V\to U$ be a finite flat morphism
of smooth schemes over $F$.
Then, the following holds.

1. 
Let $X$ be an object of ${\cal C}_U$
and $Y$ be an object of ${\cal C}_V$.
Then there exist
a morphism $X'\to X$ in ${\cal C}_U$,
a morphism $Y'\to Y$ in ${\cal C}_V$
and a finite flat morphism $\bar f':Y'\to X'$
over $F$
extending $f:V\to U$.

If $f:V\to U$ is a Galois covering of group $G$,
there exists
$\bar f':Y'\to X'$ as above
such that the action of $G$ 
is extended to an action on $Y'$.

2. 
Let
$g:X'\to X$ be a morphism in ${\cal C}_U$ and
$h:Y'\to Y$ be a morphism in ${\cal C}_V$.
Let 
\begin{equation}\label{eqGR}
\begin{CD}
Y'@>h>> Y\\
@V{\bar f'}VV @VV{\bar f}V\\
X'@>g>> X
\end{CD}
\end{equation}
be a commutative diagram of morphisms
over $F$ where the vertical arrows
extend $f:V\to U$.
Then, 
the diagram
$$\begin{CD}
CH_0(Y'\setminus V)@>{h_*}>> CH_0(Y\setminus V)\\
@A{\bar f^{\prime*}}AA @AA{\bar f^*}A \\
CH_0(X'\setminus U)@>{g_*}>> CH_0(X\setminus U)
\end{CD}$$
is commutative.
\end{lm}

{\it Proof.}
1. By replacing $Y$ by the closure
of the graph $\Gamma_f\subset V\times U
\subset Y\times X$,
we may assume that there exists
a proper map
$\bar f:Y\to X$
extending $f:V\to U$.
Then, 
we obtain
a finite flat morphism $\bar f':Y'\to X'$
by applying
Th\'eor\`eme (5.2.2)
of \cite{GR}.

Assume $V\to U$
is a Galois covering.
Then, by replacing $Y$
by the closure of the image
$V\to \prod_{\sigma\in G}V
\subset \prod_{\sigma\in G}Y$
sending $v$ to $(\sigma(v))_{\sigma\in G}$,
we may assume the action of 
$G$ on $V$ is extended
to an action on $Y$.
It suffices to apply the construction above.

2. 
Since the assertion is clear if the diagram
(\ref{eqGR}) is Cartesian,
we may assume $X'=X$.
Let $x\in X\setminus U$
be a closed point
and put $A=\hat O_{X,x}$.
For $y\in \bar f^{-1}(x)$,
we put $B_y=\hat O_{Y,y}$.
For $y'\in \bar f^{\prime -1}(x)$,
we put $B'_{y'}=\hat O_{Y',y'}$.
Then, we have
$\bar f^*([x])=\sum_{y\in \bar f^{-1}(x)}
{\rm rank}_AB_y/
[\kappa(y):\kappa(x)]\cdot [y]$
and
$\bar f^{\prime*}([x])=\sum_{y'\in \bar f^{\prime -1}(x)}
{\rm rank}_AB'_{y'}/
[\kappa(y'):\kappa(x)]\cdot [y']$.
For each $y\in \bar f^{-1}(x)$,
we have
${\rm rank}_AB_y
=\sum_{y'\in h^{-1}(y)}
{\rm rank}_AB'_{y'}$.
Thus the assertion follows.
\qed

By Lemma \ref{lmGR},
the flat pull-back maps
$\bar f^*:CH_0(X\setminus U)
\to
CH_0(Y\setminus V)$
induce
$f^*:CH_0(\overline U\setminus U)
\to
CH_0(\overline V\setminus V)$.

\begin{cor}\label{f*f!}
Let $f:V\to U$
be a finite flat morphism 
of smooth schemes
of constant degree $N$.

1. Then, the composition 
$f_*\circ f^*:
CH_0(\overline U\setminus U)
\to
CH_0(\overline U\setminus U)$
is the multiplication by $N$.

2. Assume further that
$V\to U$ is a Galois covering of Galois group $G$.
Then, the composition 
$f^*\circ f_*:
CH_0(\overline V\setminus V)
\to
CH_0(\overline V\setminus V)$
is equal to
$\sum_{\sigma\in G}\sigma^*$.

The pull-back map
$f^*$
induces an isomorphism
$f^*:
CH_0(\overline U\setminus U)
\otimes_{\mathbb Z}{\mathbb Q}
\to
(CH_0(\overline V\setminus V)
\otimes_{\mathbb Z}{\mathbb Q})^G$
to the $G$-fixed part.
The inverse is given by
$\frac1{|G|}f_*$.
\end{cor}

\noindent{\it Proof.}
Clear from Lemma \ref{lmGR}.
\qed

If we admit resolution of singularities,
the projective limit
$CH_0(\overline V\setminus V)$
is computed by a smooth object
in ${\cal C}_V$
as we see in Corollary \ref{corCV} below.

\begin{lm}\label{lmCV}
Let $V$ be a separated smooth scheme 
of finite type over $F$.
Let $Y$ and $Y'$
be separated smooth schemes over $F$
containing $V$ as dense open subschemes
and $g:Y'\to Y$
be a morphism over $F$
inducing the identity on $V$.

Then, the Gysin map
$g^!:CH_0(Y\setminus V)\to
CH_0(Y'\setminus V)$
is a surjection.
Further if $g:Y'\to Y$
is proper,
the map
$g^!:CH_0(Y\setminus V)\to
CH_0(Y'\setminus V)$
is an isomorphism and is the inverse of
$g_*:CH_0(Y'\setminus V)\to
CH_0(Y\setminus V)$.
\end{lm}

\noindent{\it Proof.}
Let ${\cal K}_d$
denote the Zariski sheaf of Quillen's K-theory.
Then, by the Gersten resolution,
the Chow group
$CH_0(Y\setminus V)$
is identified with the cohomology
$H^d_{Y\setminus V}(Y,{\cal K}_d)$
with support
and the Gysin map
$g^!:CH_0(Y\setminus V)\to
CH_0(Y'\setminus V)$
is identified with the pull-back map
$g^*:
H^d_{Y\setminus V}(Y,{\cal K}_d)\to
H^d_{Y'\setminus V}(Y',{\cal K}_d)$.
Thus,
we have a commutative diagram
of exact sequences
$$\begin{CD}
H^{d-1}(V,{\cal K}_d)@>>>
CH_0(Y\setminus V)
@>>>
CH_0(Y)
@>>>
CH_0(V)\\
@| @VV{g^!}V @VV{g^!}V @|\\
H^{d-1}(V,{\cal K}_d)@>>>
CH_0(Y'\setminus V)
@>>>
CH_0(Y')
@>>>
CH_0(V).
\end{CD}$$
Since $CH_0(Y')$
is generated by the 0-cycles on 
the dense open $V\subset Y'$,
the map $g^!:
CH_0(Y)\to CH_0(Y')$
is surjective.
Thus a diagram chasing shows the surjectivity of
$g^!:CH_0(Y\setminus V)
\to 
CH_0(Y'\setminus V)$.

If $g$ is proper,
we have $g_*\circ g^!={\rm id}$
by the projection formula.
Hence $g^!$ is an isomorphism and
is the inverse of $g_*$.
\qed

\begin{cor}\label{corCV}
Let $V$ be a separated smooth scheme of finite
finite type over $F$.
Assume the full subcategory ${\cal C}_V^{\rm sm}$
consisting of smooth objects is cofinal
in ${\cal C}_V$.

1. Then, 
the projection
$CH_0(\overline V\setminus V)\to
CH_0(Y\setminus V)$
is an isomorphism for an object $Y$
of ${\cal C}_V^{\rm sm}$.
Their inverses induce an isomorphism
$\varinjlim_{{\cal C}_V^{\rm sm,opp}}
CH_0(Y\setminus V)\to 
CH_0(\overline V\setminus V)$
where the transition maps are
Gysin maps.

2. Let $f:V\to U$
be a finite flat morphism 
of smooth schemes.
Assume the full subcategory ${\cal C}_U^{\rm sm}$
is also cofinal
in ${\cal C}_U$.

Then,
the pull-back map
$f^*:
CH_0(\overline U\setminus U)\to
CH_0(\overline V\setminus V)$
is the same as the map
$\varinjlim_{{\cal C}_U^{\rm sm,opp}}
CH_0(X\setminus U)\to 
\varinjlim_{{\cal C}_V^{\rm sm,opp}}
CH_0(Y\setminus V)$
induced by the Gysin maps.
\end{cor}

{\it Proof.}
1. Clear from Lemma \ref{lmCV}.

2. 
Let $X$ and $Y$ be objects
of ${\cal C}_U^{\rm sm}$
and of ${\cal C}_V^{\rm sm}$
respectively and let
$\bar f:Y\to X$
be a morphism over $F$
extending $f:V\to U$.
It is sufficient to show that
$\bar f^!([x])=f^*([x])$
for an arbitrary closed point $x\in X\setminus U$.
Let $X'\to X$ be the blow-up at $x$
and $Y'$ be an object of
${\cal C}_V^{\rm sm}$
dominating $Y\times_XX'$.
Replacing $Y\to X$
by $Y'\to X'$,
we may assume that the map $\bar f:Y\to X$
is finite flat
on a neighborhood of $x$.
Then, we have
$\bar f^!([x])=[\bar f^{-1}(x)]$.
By applying 
Th\'eor\`eme (5.2.2)
of \cite{GR},
we also get
$f^*([x])=[\bar f^{-1}(x)]$.
\qed

\subsection{Definition  of the
intersection product with the log diagonal}

First, we recall the existence of alteration.

\begin{lm}\label{lmalter}
Let $f:V\to U$
be a finite \'etale morphism
of separated and smooth
schemes of finite type purely of dimension $d$ over 
a perfect field $F$.
Let $Y$ be a separated scheme of finite
type over $F$
containing $V$
as a dense open subscheme.

Then, there exists
a commutative diagram
\begin{equation}
\xymatrix{
W\ar[d]_g\ar[r]^{\subset}&Z\ar[d]_{\bar g}\ar[ddr]^{\bar h}&\\
V\ar[d]_f\ar[r]^{\subset}&Y
\\
U\ar[rr]^{\subset}&&X
}\label{eqalt}
\end{equation}
satisfying the following conditions:

{\rm (\ref{eqalt}.1)} 
$U$ is the complement of a Cartier divisor $B$ of $X$.

{\rm (\ref{eqalt}.2)} 
$Z$ is
smooth purely of dimension $d$
over $F$
and
$W$ is the complement of a
divisor $D$
of $Z$ with simple normal crossings.

{\rm (\ref{eqalt}.3)}
The two quadrangles are Cartesian.

{\rm (\ref{eqalt}.4)}
$\bar g:Z\to Y$ is proper.
The map $g:W\to V$ is
a generically finite surjection
of constant degree $[W:V]$.
\end{lm}

\noindent{\it Proof.}
By \cite{nagata},
there exists a proper scheme $X$
over $F$ containing $U$
as a dense open subscheme.
By replacing $X$
by its blow-up at a closed subscheme
whose support is the complement of $U$,
the condition (\ref{eqalt}.1) is satisfied.
By replacing $Y$
by the closure of the graph
of $f:V\to U$
in $Y\times X$,
we may assume there is a 
commutative diagram
\begin{equation}
\begin{CD}
V@>{\subset}>>Y\\
@VfVV @VV{\bar f}V\\
U@>{\subset}>>X.
\end{CD}\label{eqXY}
\end{equation}
Since $V$ is proper over $U$
and is dense in $U\times_XY$,
the diagram
(\ref{eqXY})
is Cartesian.
Now, it is sufficient to
apply the main result of de Jong \cite{dJ}
to $V\subset Y$
to find $W\subset Z$.
\qed

Next, we study
the intersection product
with the log diagonal
on the level of alteration.
We consider a Cartesian diagram
\begin{equation}
\begin{CD}
W@>{\subset}>> Z  \\
@VhVV @VV{\bar h}V\\
U@>{\subset}>>  X
\end{CD}\label{eqalter}
\end{equation}
of separated schemes of finite type over $F$
satisfying the conditions:

(\ref{eqalt}.1)
$U$ is the complement of a Cartier divisor 
$B$ of $X$.

{\rm (\ref{eqalt}.2)} 
$Z$ is
smooth purely of dimension $d$
over $F$
and
$W$ is the complement of a
divisor $D$
of $Z$ with simple normal crossings.

\noindent
Let $D_1,\ldots,D_m$
be the irreducible components of $D$
and let $(Z\times Z)^\sim$
be the log product
with respect to
the divisors $D_1,\ldots,D_m$.
The scheme 
$(Z\times Z)^\sim$
is smooth over $F$
and contains 
$W\times W$ as the complement
of a divisor with simple normal crossings
by Lemma \ref{lm'}.
The log diagonal map
$\Delta_Z:Z\to (Z\times Z)^\sim$
is a regular closed immersion of codimension $d$.
Let
$(Z\times_XZ)^\sim=
(Z\times Z)^\sim\times_{(X\times X)^\sim}X$
be the relative log product
defined with respect to the Cartier divisor
$B$
and the family
$D_1,\ldots,D_m\subset Z$ of Cartier divisors.

Let $T$ be an open neighbourhood of 
$\Delta_W$ in $W\times_UW$.
Then the closure
$\overline{W\times_UW
\setminus T}$
in
$(Z\times Z)^\sim$
satisfies
$\overline{W\times_UW
\setminus T}\cap \Delta_Z
\subset Z\setminus W$
since
$\overline{W\times_UW
\setminus T}
\cap W\times W=
W\times_UW
\setminus T$.
Thus the intersection product
in $(Z\times Z)^\sim$
defines a map
\begin{equation}
\begin{CD}
(\ ,\Delta_Z)_{(Z\times Z)^\sim}:
CH_d(\overline{W\times_UW
\setminus T})
@>>>
CH_0(Z \setminus W).
\end{CD}\label{edDY}
\end{equation}

\begin{pr}\label{prcut}
Let $Z$ be a smooth scheme
purely of dimension $d$ over $F$
and $W\subset Z$ be
the complement of a divisor $D$ with
simple normal crossings.
Let $W\to U$ be a morphisms
of schemes of finite type over $F$
and $T\subset W\times_UW$
be an open neighborhood of
the diagonal $\Delta_W$.
Assume there exists a Cartesian diagram
{\rm (\ref{eqalter})}
satisfying the conditions
{\rm (\ref{eqalt}.1)} and
{\rm (\ref{eqalt}.2)}.

1. 
Let $\overline{W\times_UW\setminus T}$
be the closure in $(Z\times Z)^\sim$.
Then,
there exists
a unique map
\begin{equation}
\begin{CD}
(\ ,\Delta_Z)^{\log}:
CH_d(W\times_UW\setminus T)
@>>>
CH_0(Z \setminus W)
\end{CD}\label{eqmapZ}
\end{equation}
making the diagram
$$\xymatrix{
CH_d(\overline{W\times_UW\setminus T})
\ar[d]_{\text {\rm restriction}}
\ar[rrd]^{\qquad (\ ,\Delta_Z)_{(Z\times Z)^\sim}}&&
\\
CH_d(W\times_UW\setminus T)
\ar[rr]_{\qquad (\ ,\Delta_Z)^{\log}}&&
CH_0(Z\setminus W)
}$$
commutative.

2. 
Further, let 
$$\begin{CD}
W'@>{\subset }>>Z'\\
@V{k}VV
@VV{\bar k}V\\
W@>{\subset}>>Z
\end{CD}$$
be a Cartesian diagram of schemes over $F$.
We assume that $Z'$ is smooth over $F$ and
that $W'$ is the complement of a divisor of $Z'$
with simple normal crossings.
Then, we have a commutative diagram
\begin{equation}
\begin{CD}
CH_d(W\times_UW\setminus T)
@>{(\ ,\Delta_Z)^{\log}}>>
CH_0(Z \setminus W)\\
@V{(k\times k)^{!}}VV @VV{\bar k^!}V\\
CH_d(W'\times_UW'\setminus (k\times k)^{-1}(T))
@>{(\ ,\Delta_{Z'})^{\log}}>>
CH_0(Z' \setminus W')
\end{CD}
\end{equation}
where the left vertical arrow is the Gysin map
for $k\times k:W'\times W'
\to W\times W$.
\end{pr}

\noindent{\it Proof.}
1. Take a Cartesian diagram
{\rm (\ref{eqalter})}
satisfying the conditions
{\rm (\ref{eqalt}.1)} and
{\rm (\ref{eqalt}.2)}.
Then $(Z\times_XZ)^\sim$
is a closed subscheme
of $(Z\times Z)^\sim$
containing $W\times_UW$
as an open subscheme.
Hence, 
$\overline{W\times_UW\setminus T}$
is closed in $(Z\times_XZ)^\sim
\setminus T$
and $W\times_UW\setminus T$
is open in $(Z\times_XZ)^\sim
\setminus T$.
Thus,
it suffices to show that the map
$(\ ,\Delta_Z)_{(Z\times Z)^\sim}:
CH_d((Z\times_XZ)^\sim
\setminus T)
\to CH_0(Z\setminus W)$
factors through
the restriction map
$CH_d((Z\times_XZ)^\sim
\setminus T)
\to
CH_d(W\times_UW\setminus T)$.
The kernel of the surjection
$CH_d((Z\times_XZ)^\sim
\setminus T)
\to
CH_d(W\times_UW\setminus T)$
is generated by the image of
$CH_d((Z\times_XZ)^\sim
\setminus W\times_UW)$.

We use the notation 
in Lemma \ref{lm'}
replacing $X\to Y$ by $Z\to X$.
Then, the complement
$(Z\times Z)^\sim\setminus (W\times W)$
is the union of divisors
$E_i^\circ$.
Hence 
the complement $(Z\times_XZ)^\sim
\setminus (W\times_UW)
=(Z\times_XZ)^\sim
\cap ((Z\times Z)^\sim\setminus (W\times W))$
is the union of
$(Z\times_XZ)^\sim\cap E_i^\circ$.
Thus the kernel of 
the restriction map
$CH_d((Z\times_XZ)^\sim
\setminus T)
\to CH_d(W\times_UW\setminus T)$
is generated by the images of
$CH_d((Z\times_XZ)^\sim\cap E_i^\circ)$.

The pull-back of the Cartier divisor
$E_i^\circ \subset (Z\times Z)^\sim$
by the log diagonal map
$\Delta_Z\to (Z\times Z)^\sim$
is the Cartier divisor $\Delta_{D_i}\subset \Delta_Z$.
Hence we have 
$(C,\Delta_Z)_{(Z\times Z)^\sim}
=
(C,\Delta_{D_i})_{E_i^\circ}$
for a cycle $C$ in $E_i^\circ$.
Thus, it is sufficient to show that
the map
\begin{equation}
\begin{CD}
(\ ,\Delta_{D_i})_{E_i^\circ}
:CH_d((Z\times_XZ)^\sim\cap E_i^\circ)
@>>> CH_0(D_i)
\end{CD}
\label{eqEi}
\end{equation}
is the 0-map.

The log diagonal map
$D_i\to (D_i\times D_i)^\sim$
is a regular immersion of codimension $d-1$.
The restriction $E_{i,D_i}^\circ$ of
the ${\mathbb G}_m$-bundle
$E_i^\circ\to (D_i\times D_i)^\sim$
to the log diagonal
$D_i\subset (D_i\times D_i)^\sim$
has a canonical isomorphism
$E_{i,D_i}^\circ\to {\mathbb G}_{m,D_i}$
(\ref{eqEPD}).
The immersion
$\Delta_{D_i}=\Delta_Z\cap E_i^\circ \to E_i^\circ$
gives the unit section
$D_i\to E_{i,D_i}^\circ\to {\mathbb G}_{m,D_i}$.
Hence the map (\ref{eqEi})
is the composition of the maps
\begin{equation}
\begin{CD}
CH_d((Z\times_XZ)^\sim\cap E_i^\circ )
@>{(\ ,D_i)_{(D_i\times D_i)^\sim}}>>
CH_1((Z\times_XZ)^\sim\cap 
E_{i,D_i}^\circ )\\
@>{(\ ,D_i)_{E_{i,D_i}^\circ}}>>
CH_0(D_i).
\end{CD}
\label{eqsec}
\end{equation}

By Proposition \ref{prmu}.1,
the intersection
$(Z\times_XZ)^\sim\cap 
E_{i,D_i}^\circ$
is a closed subscheme of
$\mu_{e_i,D_i}
\subset E_{i,D_i}^\circ={\mathbb G}_{m,D_i}$.
Hence the second map in (\ref{eqsec})
is the composition
$$\begin{CD}
CH_1((Z\times_XZ)^\sim\cap 
E_{i,D_i}^\circ )
\to
CH_1(\mu_{e_i,D_i})
\to 
CH_1({\mathbb G}_{m,D_i})
@>{(\ ,D_i)_{{\mathbb G}_{m,D_i}}}>>
CH_0(D_i).
\end{CD}$$
Since 
the composition of the last two maps is the 0-map,
the map
$(\ ,\Delta_Z)_{(Z\times Z)^\sim}:
CH_d((Z\times_XZ)^\sim
\setminus T)
\to CH_0(Z\setminus W)$
induces a map
$CH_d(W\times_UW\setminus T)
\to
CH_0(Z \setminus W)$.
Thus the assertion follows.

2. 
We consider the commutative diagram
$$\begin{CD}
(Z'\times Z')^\sim@<<< 
(Z'\times_{X'}Z')^\sim@<<< Z'\\
@V{(\bar k\times \bar k)^\sim}VV 
@V{(\bar k\times \bar k)^\sim}VV 
@VV{\bar k}V \\
(Z\times Z)^\sim@<<< 
(Z\times_{X'}Z)^\sim@<<< Z
\end{CD}$$
where the right horizontal arrows
are the log diagonal maps.
Then, we have a commutative diagram
$$\begin{CD}
CH_d((Z\times_{X'}Z)^\sim\setminus T)
@>{(\ ,\Delta_Z)_{(Z\times Z)^\sim}}>>
CH_0(Z \setminus W)\\
@V{(\bar k\times \bar k)^{\sim!}}VV @VV{\bar k^!}V\\
CH_d((Z'\times_{X'}Z')^\sim\setminus (k\times k)^{-1}(T))
@>{(\ ,\Delta_{Z'})_{(Z'\times Z')^\sim}}>>
CH_0(Z' \setminus W')
\end{CD}$$
and the assertion follows.
\qed

\begin{thm}\label{corindep}
Let $f:V\to U$
be a finite \'etale morphism
of separated and smooth
schemes of finite type
purely of dimension $d$ over 
a perfect field $F$.

1. There exists a unique map
$$\begin{CD}
(\ ,\Delta_{\overline V})^{\log}:
CH_d(V\times_UV\setminus \Delta_V)
@>>>
CH_0(\overline V \setminus V)\otimes_{\mathbb Z}{\mathbb Q}
\end{CD}\leqno{(\ref{eqmap})}
$$
that makes the diagram
\begin{equation}
\xymatrix{
CH_d(V\times_UV\setminus \Delta_V)
\ar[r]^{(\ ,\Delta_{\overline V})^{\log}}
\ar[d]_{(g\times g)^!}&
CH_0(\overline V \setminus V)\otimes_{\mathbb Z}{\mathbb Q}
\ar[r]&
CH_0(Y \setminus V)\otimes_{\mathbb Z}{\mathbb Q}\\
CH_d(W\times_UW\setminus W\times_VW)
\ar[r]^{\quad \quad \quad (\ ,\Delta_Z)^{\log}}&
CH_0(Z \setminus W)
\ar[ru]_{\frac 1{[W:V]} \bar g_*}
&}\label{eqindep}
\end{equation}
commutative
for an arbitrary commutative diagram 
{\rm (\ref{eqalt})} 
satisfying the condition

{\rm (\ref{eqalt}.0)}
$Y$ contains $V$ as a dense open subscheme.

\noindent
and the conditions
{\rm (\ref{eqalt}.1)-\rm (\ref{eqalt}.4)}.

2. Assume the full subcategory ${\cal C}_V^{\rm sm,0}$
is cofinal in ${\cal C}_V$.
Then, there exists a unique map
\begin{equation}
\begin{CD}
(\ ,\Delta_{\overline V})^{\log}_{\mathbb Z}:
CH_d(V\times_UV\setminus \Delta_V)
@>>>
CH_0(\overline V \setminus V)
\end{CD}\label{eqmapVZ}
\end{equation}
satisfying the following property.

Let $Y$ be an arbitrary smooth separated scheme of finite type
containing $V$ as the complement
of a divisor with simple normal crossings
and let
$(\ ,\Delta_Y)^{\log}:
CH_d(V\times_UV\setminus \Delta_V)
\to
CH_0(Y \setminus V)$ be
the map {\rm (\ref{eqmapZ})}
for $Z=Y$.
Then the diagram
\begin{equation}
\xymatrix{
CH_d(V\times_UV\setminus \Delta_V)
\ar[rr]^{(\ ,\Delta_{\overline V})^{\log}_{\mathbb Z}}
\ar[rrd]_{(\ ,\Delta_Y)^{\log}}&&
CH_0(\overline V \setminus V)
\ar[d]&\\
&&
CH_0(Y \setminus V)
}\label{eqindep0}
\end{equation}
is commutative if there exists a Cartesian diagram
$$\begin{CD}
V@>>> Y\\
@VfVV @VV{\bar f}V\\
U@>>> X
\end{CD}$$
of separated scheme of finite type
where $X$ contains $U$ as the complement of
a Cartier divisor.
\end{thm}

\noindent{\it Proof.}
1. We consider 
an arbitrary commutative diagram 
{\rm (\ref{eqalt})} 
satisfying the conditions
{\rm (\ref{eqalt}.0)}-{\rm (\ref{eqalt}.4)}.
By the assumption that
$V\to U$ is \'etale,
the fiber product
$T=W\times_VW$ is an open 
neighborhood of $\Delta_W$
in $W\times_UW$
and the map
$(\ ,\Delta_Z)^{\log}:
CH_d(W\times_UW\setminus W\times_VW)
\to
CH_0(Z \setminus W)$ 
is defined by
Proposition \ref{prcut}.1.

For an object $Y$ of ${\cal C}_V$,
there exists
a  commutative diagram 
{\rm (\ref{eqalt})} 
satisfying the conditions
{\rm (\ref{eqalt}.1)}-{\rm (\ref{eqalt}.4)}
by Lemma \ref{lmalter}.
The composition
$CH_d(V\times_UV\setminus \Delta_V)
\to
CH_0(Y \setminus V)\otimes_{\mathbb Z}{\mathbb Q}$
via the lower line in (\ref{eqindep})
is independent of the choice of
diagram {\rm (\ref{eqalt})} 
by Proposition \ref{prcut}.2.
We define the map
$(\ ,\Delta_{\overline V})^{\log}:
CH_d(V\times_UV\setminus \Delta_V)
\to
CH_0(\overline V \setminus V)\otimes_{\mathbb Z}{\mathbb Q}$
as the limit.
Then it is clear that
the map
$(\ ,\Delta_{\overline V})^{\log}$ satisfies the 
condition.

2. 
By the assumption and Corollary \ref{corCV}.1,
the group
$CH_0(\overline V\setminus V)$
is identified with the inductive limit
$\varinjlim_{{\cal C}_V^{\rm sm,0,opp}}
CH_0(Y\setminus V)$
with respect to the Gysin maps.
Hence
it follows from 
Proposition \ref{prcut}.
\qed

If ${\cal C}_V^{\rm sm,0}$
is cofinal in ${\cal C}_V$,
the map
$(\ ,\Delta_{\overline V})^{\log}$
is induced by
$(\ ,\Delta_{\overline V})^{\log}_{\mathbb Z}$.

We prove the trace formula (\ref{eqGal0})
in Proposition \ref{prtr}.
Let $V\to U$ be a finite \'etale morphism
of separated smooth schemes of dimension $d$ over $F$.
Let $\ell$ be a prime number invertible in $F$
and $\bar F$ be an algebraic closure of $F$.
For an open and closed subscheme
$\Gamma$  of 
$V\times_UV\setminus \Delta_V$,
we define 
an endomorphism $\Gamma^*$
of $H^q_c(V_{\bar F},{\mathbb Q}_\ell)$
to be $p_{1*}\circ p_2^*$.
We put
$${\rm Tr}(\Gamma^*:H^*_c(V_{\bar F},{\mathbb Q}_\ell))
=
\sum_{q=0}^{2d}(-1)^q
{\rm Tr}(\Gamma^*:H^q_c(V_{\bar F},{\mathbb Q}_\ell)).$$

\begin{pr}\label{prtr}
Let $f:V\to U$
be a finite \'etale morphism
of separated and smooth
schemes of finite type
purely of dimension $d$ over 
a perfect field $F$.
Let $\ell$ be a prime number invertible in $F$.

Then, for an open and closed
subscheme $\Gamma$ of 
$V\times_UV\setminus \Delta_V$,
we have
\begin{equation}
{\rm Tr}(\Gamma^*:H^*_c(V_{\bar F},{\mathbb Q}_\ell))
={\rm deg}(\Gamma ,\Delta_{\overline V})^{\log}.
\label{eqLef}
\end{equation}
\end{pr}

\noindent
{\it Proof.}
Take a diagram (\ref{eqalt})
with $X,Y$ and $Z$ proper over $F$
satisfying the conditions
(\ref{eqalt}.0)-(\ref{eqalt}.4).
By Lemma \ref{lmtr},
we have
$${\rm Tr}(\Gamma^*:H^*_c(V_{\bar F},{\mathbb Q}_\ell))=
\frac 1{[W:V]}
{\rm Tr}(((g\times g)^*\Gamma)^*:H^*_c(W_{\bar F},{\mathbb Q}_\ell)).$$
By Lemma \ref{lmcompa},
we have
$(g\times g)^*[\Gamma]=
[(g\times g)^!(\Gamma)]$.
Take an element 
$\tilde \Gamma=\sum_in_i[C_i]\in
Z_d(W\times_UW\setminus W\times_VW)$
representing
$[(g\times g)^!(\Gamma)]\in
CH_d(W\times_UW\setminus W\times_VW)$.
By Proposition \ref{prmu}.2,
the closures $\overline C_i
\subset (Z\times Z)'$ satisfy
the condition (\ref{eqthm}).
Hence by Theorem \ref{thmLTF},
we have
\begin{eqnarray*}
{\rm Tr}(((g\times g)^*\Gamma)^*:H^*_c(W_{\bar F},{\mathbb Q}_\ell))
&=&{\rm deg}\ (\textstyle{\sum_i}n_i
[\overline C_i],\Delta_Z)_{(Z\times Z)'}\\
&=&{\rm deg}\ ((g\times g)^!\Gamma,\Delta_Z)^{\log}.
\end{eqnarray*}
By the definition of
$(\Gamma,\Delta_{\overline V})^{\log}$,
we have
$${\rm deg}\ (\Gamma,\Delta_{\overline V})^{\log}=
\frac1{[W:V]}\cdot 
{\rm deg}\ ((g\times g)^!\Gamma,\Delta_Z)^{\log}.$$
Thus the equality (\ref{eqLef}) is proved.
\qed

\subsection{Properties of the
intersection product with the log diagonal}

We keep the notation that
$f:U\to V$ denotes
a finite \'etale morphism of
separated smooth schemes
of finite type purely of dimension $d$
over a perfect field $F$.

The maps 
$(\ ,\Delta_{\overline V})^{\log}:
CH_d(V\times_UV\setminus \Delta_V)\to
CH_0(\overline V\setminus V)\otimes_{\mathbb Z}{\mathbb Q}$
satisfies the following functoriality.

\begin{lm}\label{lmfun}
Let $U$ be a
separated smooth scheme
of finite type purely of dimension $d$
over a perfect field $F$.

1.
Let $V\to U'$ be a morphism 
of finite and \'etale schemes over $U$.
Then the map $(\ ,\Delta_{\overline V})^{\log}:
CH_d(V\times_{U'}V\setminus \Delta_V)
\to
CH_0(\overline V\setminus V)\otimes_{\mathbb Z}{\mathbb Q}$
is equal to the restriction of
$(\ ,\Delta_{\overline V})^{\log}:
CH_d(V\times_UV\setminus \Delta_V)
\to
CH_0(\overline V\setminus V)\otimes_{\mathbb Z}{\mathbb Q}$.

2. 
Let $g:V\to V'$ be a morphism 
of finite and \'etale schemes over $U$.
Then, the diagram
\begin{equation}
\begin{CD}
CH_d(V'\times_UV'\setminus \Delta_{V'})
@>{(\ ,\Delta_{\overline {V'}})^{\log}}>>
CH_0(\overline {V'}\setminus V')\otimes_{\mathbb Z}{\mathbb Q}\\
@V{(g\times g)^*}VV @VV{g^*}V\\
CH_d(V\times_UV\setminus \Delta_V)
@>{(\ ,\Delta_{\overline{V}})^{\log}}>>
CH_0(\overline {V}\setminus V)\otimes_{\mathbb Z}{\mathbb Q}
\end{CD}
\label{eqpull1}
\end{equation}
is commutative.
\end{lm}

\noindent{\it Proof.}
1. Clear from the definition
and Proposition \ref{prcut}.2.

2. 
We may assume 
$U$, $V$ and $V'$ are connected.
Then by Corollary \ref{f*f!}.1,
the right vertical arrow $g^*$ in (\ref{eqpull1})
is injective.
Hence, we may replace $V$ by its
Galois closure over $U$
and may assume $V\to U$
is a Galois covering.
Let $G$ be the Galois group 
and $H\subset G$ be the
subgroup corresponding
to $V'$.
Then, the images of the
both compositions
are in the $H$-fixed part of
$CH_0(\overline V\setminus V)\otimes_{\mathbb Z}{\mathbb Q}$.
Hence,
by Corollary \ref{f*f!}.2,
it suffices to show
the diagram
$$\begin{CD}
CH_d(V'\times_UV'\setminus \Delta_{V'})
@>{(\ ,\Delta_{\overline {V'}})^{\log}}>>
CH_0(\overline {V'}\setminus V')\otimes_{\mathbb Z}{\mathbb Q}\\
@V{(g\times g)^*}VV @AA{\frac 1{|H|} g_*}A\\
CH_d(V\times_UV\setminus \Delta_V)
@>{(\ ,\Delta_{\overline{V}})^{\log}}>>
CH_0(\overline {V}\setminus V)\otimes_{\mathbb Z}{\mathbb Q}
\end{CD}
$$
is commutative.
This is clear from the definition
and Proposition \ref{prcut}.2.
\qed

If the subcategory
${\cal C}_V^{\rm sm,0}$
is cofinal in 
${\cal C}_V$,
we can remove
$\otimes_{\mathbb Z}{\mathbb Q}$
in 1.
Further if 
${\cal C}_{V'}^{\rm sm,0}$
is cofinal in 
${\cal C}_{V'}$,
we can remove
$\otimes_{\mathbb Z}{\mathbb Q}$
in 2.
We will omit to state remarks
on integrality of this type in the sequel.

\begin{lm}\label{lmkum}
Let $V\to U'\to U$
be finite \'etale morphisms
of separated and smooth
schemes of finite type
purely of dimension $d$ over 
a perfect field $F$.
Let $n\ge 1$ be an integer
invertible in $F$
and assume
$g:U'\to U$ is a  ${\mathbb Z}/n{\mathbb Z}$-torsor
over $U$.

Then, the restriction
$$\begin{CD}
CH_d(V\times_UV\setminus 
V\times_{U'}V)\subset
CH_d(V\times_UV\setminus \Delta_V)
@>{(\ ,\Delta_{\overline V})^{\log}}>>
CH_0(\overline V\setminus V)
\otimes_{\mathbb Z}{\mathbb Q}
\end{CD}$$
is the $0$-map.
\end{lm}

\noindent{\it Proof.}
By enlarging $F$,
we may assume $F$ contains a primitive 
$n$-th root of 1.
Let $\chi:{\mathbb Z}/n{\mathbb Z}
\to F^\times$
be a character of order $n$.
Then the $\chi$-part
${\cal L}_U$ of $g_*O_{U'}$ is 
an invertible $O_U$-module.
The multiplication defines an isomorphism
$\mu_U:{\cal L}_U^{\otimes n}\to O_U$
of $O_U$-modules.
The $O_U$-algebra $g_*O_{U'}$ is isomorphic to 
$\bigoplus_{i=0}^{n-1}
{\cal L}_U^{\otimes i}$
with the multiplication defined by
$\mu_U$.
We take a proper scheme $X$ over $F$
containing $U$ as
a dense open subscheme.
Replacing $X$ by a blow-up,
we may assume ${\cal L}_U$
is extended to an invertible $O_X$-module ${\cal L}$
and the map
$\mu_U:{\cal L}_U^{\otimes n}\to O_U$
is extended to an injection
$\mu:{\cal L}^{\otimes n}\to O_X$.
We define a finite flat scheme $\bar g:X'\to X$
over $X$ by the $O_X$-algebra
$\bigoplus_{i=0}^{n-1}{\cal L}^{\otimes i}$
with the multiplication defined by $\mu$.
By Lemma \ref{lmn},
the diagonal 
$X'\to (X'\times_XX')^\sim$
is an open immersion.

We take a proper scheme $Y$
containing $V$ as a dense open subscheme
such that the map $V\to U'$
is extended to $Y\to X'$
and an alteration $Z\to Y$ as in
Lemma \ref{lmalter}.
Then, the inverse image of
$V\times_UV\setminus V\times_{U'}V$
in $(Z\times Z)^\sim$
is contained in the inverse image of
$(X'\times_XX')^\sim\setminus X'$.
Thus the assertion follows from
the definition of the map
(\ref{eqmap}).
\qed

For a separated scheme $Y$ 
of finite type over $F$
containing $V$ as a dense open subscheme,
let
\begin{equation}
\begin{CD}
(\ ,\Delta_Y)^{\log}:
CH_d(V\times_UV\setminus \Delta_V)
@>>>
CH_0(Y\setminus V)\otimes_{\mathbb Z}{\mathbb Q}
\end{CD}
\label{eqmapYQ}
\end{equation}
denote the composition of the maps
in the upper line of the diagram
(\ref{eqindep}).
The map
$(\ ,\Delta_Y)^{\log}:
CH_d(V\times_UV\setminus \Delta_V)
\to
CH_0(Y\setminus V)\otimes_{\mathbb Z}{\mathbb Q}$
is characterized by the commutativity of
the diagram 
\begin{equation}
\begin{CD}
CH_d(V\times_UV\setminus \Delta_V)
@>{(\ ,\Delta_Y)^{\log}}>>
CH_0(Y \setminus V)\otimes_{\mathbb Z}{\mathbb Q}\\
@V{(g\times g)^!}VV
@AA{\frac 1{[W:V]} \bar g_*}A\\
CH_d(W\times_UW\setminus W\times_VW)
@>{(\ ,\Delta_Z)^{\log}}>>
CH_0(Z \setminus W)
\end{CD}
\label{eqindepY}
\end{equation}
for an arbitrary commutative diagram 
{\rm (\ref{eqalt})} 
satisfying the conditions
{\rm (\ref{eqalt}.1)-\rm (\ref{eqalt}.4)}.
If $Y$ is smooth and
$V$ is the complement of
a divisor with simple normal crossings,
the map
$(\ ,\Delta_Y)^{\log}$
is given by
the map (\ref{eqmapZ})
for $Z=Y$.

We give sufficient conditions for
the vanishing of the map
$(\ ,\Delta_Y)^{\log}:
CH_d(V\times_UV\setminus \Delta_V)
\to
CH_0(Y\setminus V)\otimes_{\mathbb Z}{\mathbb Q}$.

\begin{lm}\label{lmet}
Let
$$\begin{CD}
V@>{\subset}>> Y\\
@VfVV @VV{\bar f}V\\
U@>{\subset}>> X
\end{CD}\leqno{\rm (\ref{eqXY})}
$$
be a Cartesian diagram
of separated schemes of finite type over a perfect field $F$.
We assume
$U\subset X$ and $V\subset Y$
are dense open subschemes,
$U$ is smooth purely of dimension $d$
over $F$
and 
$f:V\to U$ is finite and \'etale. 

Let $\Gamma\subset V\times_UV$
be an open and closed subscheme.
If the intersection $\overline \Gamma\cap \Delta_Y$
of the closure $\overline \Gamma
\subset Y\times_XY$ of $\Gamma$
and the diagonal $Y \subset Y\times_XY$
is empty,
we have
$(\Gamma,\Delta_Y)^{\log}=0$
in $CH_0(Y\setminus V)\otimes_{\mathbb Z}{\mathbb Q}$.
\end{lm}

\noindent{\it Proof.}
By replacing $X$ by a blow-up,
we may assume $U\subset X$
is the complement of a Cartier divisor.
We take a Cartesian diagram
$$\begin{CD}
W@>{\subset}>> Z\\
@VgVV @VV{\bar g}V\\
V@>{\subset}>> Y
\end{CD}$$
satisfying the conditions
{\rm (\ref{eqalt}.2)} and {\rm (\ref{eqalt}.4)}.
We consider the natural map
$\bar g\times \bar g:
(Z\times_XZ)^\sim\to Y\times_XY$
induced by $\bar g:Z\to Y$.
The closure of
$(g\times g)^{-1}(\Gamma)$
is in
$(\bar g\times \bar g)^{-1}(\overline \Gamma)$
and does not meet
the log diagonal
$\Delta_Z
\subset 
(\bar g\times \bar g)^{-1}(\Delta_Y)$
by the assumption.
Hence we have
$((g\times g)^!(\Gamma),\Delta_Z)_{(Z\times Z)^\sim}=0$
and 
the assertion follows.
\qed

\begin{cor}\label{coret}
Let the notation be as in Lemma \ref{lmet}.

1. If $\bar f:Y\to X$ is \'etale,
the map
$(\ ,\Delta_Y)^{\log}:
CH_d(V\times_UV\setminus \Delta_V)
\to
CH_0(Y\setminus V)\otimes_{\mathbb Z}{\mathbb Q}$
is the 0-map.

2. Let $\bar \sigma$ be an automorphism of $Y$ over $X$
and $\sigma$ be the restriction on $V$.
Let $\Gamma_\sigma \subset V\times_UV$
and $\Gamma_\sigma \subset Y\times_XY$
be the graphs.
If $Y^\sigma=\Gamma_\sigma\cap \Delta_Y$
is empty,
we have
$(\Gamma_\sigma ,\Delta_Y)^{\log}=0$
in $CH_0(Y\setminus V)\otimes_{\mathbb Z}{\mathbb Q}$.
\end{cor}

\noindent{\it Proof.}
Clear from Lemma \ref{lmet}.
\qed

We show that the image of the map 
$(\ ,\Delta_Y)^{\log}:
CH_d(V\times_UV\setminus \Delta_V)\to
CH_0(Y\setminus V)\otimes_{\mathbb Z}{\mathbb Q}$
is supported on the wild ramification locus.

\begin{pr}\label{prtame}
Let
$$\begin{CD}
V@>{\subset}>> Y\\
@VfVV @VV{\bar f}V\\
U@>{\subset}>> X
\end{CD}\leqno{\rm (\ref{eqXY})}
$$
be a Cartesian diagram
of separated schemes of finite type over 
a perfect field $F$.
We assume
$X$ is smooth purely of dimension $d$
over $F$,
$U$ is the complement of
a divisor $B$ with simple normal crossings,
$V\subset Y$
is a dense open subscheme,
and 
$f:V\to U$ is finite and \'etale. 

1. Let $V\subset V'\subset Y$ be 
an open normal subscheme.
If $V'$ is 
tamely ramified over $X$,
then the map
$(\ ,\Delta_Y)^{\log}:
CH_d(V\times_UV\setminus \Delta_V)\to
CH_0(Y\setminus V)\otimes_{\mathbb Z}{\mathbb Q}$
is decomposed as the composition
$$
CH_d(V\times_UV\setminus \Delta_V)\to
CH_0(Y\setminus V')\otimes_{\mathbb Z}{\mathbb Q}
\to
CH_0(Y\setminus V)\otimes_{\mathbb Z}{\mathbb Q}.$$

2. Suppose there exists
a commutative diagram
$$\begin{CD}
V@>{\subset}>> Y\\
@VVV @VVV\\
U'@>{\subset}>> X'\\
@VgVV @VV{\bar g}V\\
U@>{\subset}>> X
\end{CD}
$$
of separated normal schemes
of finite type over $F$,
$g:U'\to U$ is finite \'etale
and $\bar g:X'\to X$ is tamely ramified.
Then, the restriction
$$\begin{CD}
CH_d(V\times_UV\setminus 
V\times_{U'}V)\subset
CH_d(V\times_UV\setminus \Delta_V)
@>{(\ ,\Delta_Y)^{\log}}>>
CH_0(Y\setminus V)\otimes_{\mathbb Z}{\mathbb Q}
\end{CD}$$
is the 0-map.
\end{pr}

\noindent{\it Proof.}
It follows from 
the characterization of the map
(\ref{eqmapYQ}) and
Lemma \ref{lmtame}.
\qed

\subsection{Wild differents
and log Lefschetz classes}

\begin{df}\label{dfdis}
Let $f:V\to U$
be a finite \'etale morphism
of separated and smooth
schemes of finite type
purely of dimension $d$ over 
a perfect field $F$.

1. We call the {\rm 0}-cycle class
\begin{equation}
D^{\log}_{V/U}=(V\times_UV\setminus \Delta_V,
\Delta_{\overline V})^{\log}
\in 
CH_0(\overline V\setminus V)\otimes_{\mathbb Z}{\mathbb Q}
\label{eqwd}
\end{equation}
the wild different of $V$ over $U$.

2.
Let $\sigma$ be an automorphism of $V$
over $U$
that is not the identity on
any component of $V$.
Let $\Gamma_\sigma\subset V\times_UV$
be the graph of $\sigma$.
Then, we call
the {\rm 0}-cycle class
\begin{equation}
(\Gamma_\sigma,\Delta_{\overline V})^{\log}\in 
CH_0(\overline V\setminus V)\otimes_{\mathbb Z}{\mathbb Q}
\label{eqLc}
\end{equation}
the log Lefschetz class
of $\sigma$.
\end{df}

If the subcategory
${\cal C}_V^{\rm sm,0}$
is cofinal in
${\cal C}_V$,
the wild different
$D^{\log}_{V/U}$
and the log Lefschetz class
$(\Gamma_\sigma,\Delta_{\overline V})^{\log}$
are defined in 
$CH_0(\overline V\setminus V)$.

\begin{lm}\label{lmchain}
For a morphism $g:V\to V'$
of finite and \'etale schemes over $U$, we have
\begin{equation}
D^{\log}_{V/U}=
D^{\log}_{V/V'}+
g^*D^{\log}_{V'/U}.
\label{eqD}
\end{equation}
\end{lm}

\noindent{\it Proof.}
We have $V\times_UV\setminus \Delta_V=
(V\times_{V'}V\setminus \Delta_V)\amalg
(g\times g)^{-1}
(V'\times_UV'\setminus \Delta_{V'})$.
Hence, the equalities follow from 
Lemma \ref{lmfun}.2.
\qed

\begin{pr}\label{prord}
Let $f:V\to U$
be a finite and \'etale morphism
of connected separated and smooth
schemes of finite type 
purely of dimension $d$ over 
a perfect field $F$
and let $\sigma$ be an automorphism of $V$
over $U$.

If the order of $\sigma$
is not a power of the characteristic $p$ of $F$,
we have
$$(\Gamma_\sigma,\Delta_{\overline V})^{\log}=0.$$
\end{pr}

\noindent{\it Proof.}
Let $n$ be the prime-to-$p$
part of the order $e$ of $\sigma$.
By Lemma \ref{lmfun}.1,
we may replace $U$ by 
the quotient $V/\langle\sigma\rangle$.
Then it suffices to apply Lemma
\ref{lmkum}
to $V\to U'=V/\langle \sigma^n\rangle
\to U=V/\langle\sigma\rangle$.
\qed

We expect the following holds.

\begin{cn}\label{cnord}
Let $f:V\to U$
be a finite and \'etale morphism
of connected separated and smooth
schemes of finite type 
purely of dimension $d$ over 
a perfect field $F$
and let $\sigma$ be a non-trivial automorphism of $V$
over $U$.

If $j$ is an integer prime to the order of $\sigma$,
we have $$(\Gamma_\sigma, \Delta_{\overline V})^{\log}
=(\Gamma_{\sigma^j}, \Delta_{\overline V})^{\log}.$$
\end{cn}

We will prove Conjecture \ref{cnord}
assuming $\dim \le 2$ in Lemma \ref{lmord}.

\begin{lm}\label{lmtr2}
Let the notation be as in Definition \ref{dfdis}
and let $\ell$ be a prime number invertible in 
a perfect field $F$.

1. If $f:V\to U$ is of constant degree $[V:U]$,
we have
\begin{equation}
{\rm deg}\ D_{V/U}^{\log}=
[V:U]\chi_c(U_{\bar F},{\mathbb Q}_\ell)
-\chi_c(V_{\bar F},{\mathbb Q}_\ell)
\label{eqd}
\end{equation}

2. Let $\sigma$ be an automorphism of $V$
over $U$
that is not the identity on
any component of $V$.
Then, we have
\begin{equation}
{\rm deg}\ 
(\Gamma_\sigma,\Delta_{\overline V})^{\log}
={\rm Tr}(\sigma^*:H^*_c(V_{\bar F},{\mathbb Q}_\ell)).
\label{eqauto}
\end{equation}
\end{lm}

\noindent
{\it Proof.}
1. By Proposition \ref{prtr},
we have
$${\rm deg}\ D_{V/U}^{\log}
={\rm Tr}((V\times_UV)^*:H^*_c(V_{\bar F},{\mathbb Q}_\ell))-
{\rm Tr}(\Delta_V^*:H^*_c(V_{\bar F},{\mathbb Q}_\ell)).$$
By Lemma \ref{lmtr},
we have
${\rm Tr}((V\times_UV)^*:H^*_c(V_{\bar F},{\mathbb Q}_\ell))=
[V:U]\cdot \chi_c(U_{\bar F},{\mathbb Q}_\ell).$
Hence the assertion follows.

2. It suffices to apply Proposition \ref{prtr}
to the graph $\Gamma_\sigma$.
\qed

\begin{cor}\label{corvan}{\rm (\cite{Il} Lemma 2.5)}
Let $U$ be a separated scheme of finite
type over a field $F$
and $V\to U$ be an \'etale ${\mathbb Z}/n{\mathbb Z}$-torsor.
Let $\sigma$ be the automorphism
defined by the generator $1\in {\mathbb Z}/n{\mathbb Z}$
and assume $n$ is not
a power of $p$.
Then, we have
$${\rm Tr}(\sigma^*:H^*_c(V_{\bar F},{\mathbb Q}_\ell))=0.$$
\end{cor}

\noindent{\it Proof.}
We may assume $F$ is perfect.
If the assertion holds for
the base changes to
a closed subscheme $Z\subset U$
and to the complement $U\setminus Z$,
it holds for $U$.
Hence, by induction on dimension,
it is reduced to the case where $U$ is smooth.
Then it follows from Lemma \ref{lmtr2}.2
and Proposition \ref{prord}.
\qed

In the rest of this subsection,
we give some computations of
wild differents
and log Lefschetz classes.

In the classical case where
$U$ is a smooth curve over $F$,
Definition \ref{dfdis}
gives the classical invariants
of wild ramifications as follows.
Let $A$ be a complete discrete valuation ring
and $B$ be the integral closure of $A$
in a finite separable extension $L$
of the fraction field $K$.
Let $e_{L/K}$
be the ramification index
of $L$ over $K$.
Then the wild different
$D^{\log}_{B/A}
\in {\mathbb N}$
is defined by
$$D^{\log}_{B/A}
={\rm length}_B
\Omega^1_{B/A}-
(e_{L/K}-1).$$
For a non-trivial automorphism $\sigma$
of $L$ over $K$,
we put
$$j_B(\sigma)=
{\rm length}_B
B/\left(\frac{\sigma(b)}b-1;b\in B\setminus\{0\}\right).
$$

\begin{lm}\label{lmcurve}
Let $U$ be a smooth connected curve over 
a perfect field $F$
and 
$f:V\to U$ 
be a finite \'etale morphism over $F$.
Let $X$ be the proper smooth curve
containing $U$ as a dense open subscheme
and $\bar f:Y\to X$ be the 
normalization in $V$.
We put $B=X\setminus U$
and $D=Y\setminus V$
and identify
$CH_0(\overline V\setminus V)=
\bigoplus_{y\in D}{\mathbb Z}$.

1. We have
$$D^{\log}_{V/U}
=
[{\rm Coker}(\bar f^*\Omega^1_{X/F}(\log B)
\to
\Omega^1_{Y/F}(\log D))]
=\sum_{y\in D}D^{\log}_{\hat O_{Y,y}/\hat O_{X,\bar f(y)}}
\cdot [y].
$$

2. Let $\sigma$ be a non-trivial automorphism
of $V$ over $U$.
Then, we have
$$(\Gamma_\sigma,\Delta_{\overline V})^{\log}_{\mathbb Z}=
\sum_{y\in D,\sigma(y)=y}j_{\hat O_{Y,y}}(\sigma)\cdot [y].$$
\end{lm}

\noindent
{\it Proof.}
Follows from 
Proposition \ref{prdis}
and
Lemma \ref{lmsgr}
below.
\qed

We compute the wild different
assuming a strong form of resolution.
Before doing it,
we recall some general facts
on intersection theory and
localized Chern classes.

Let $X$ be a scheme of finite type
over $F$ and $Z\subset X$
be a closed subscheme.
Let ${\cal E}$
and ${\cal F}$ be
locally free $O_X$-modules
of rank $d$
and $f:{\cal E}\to {\cal F}$
be an $O_X$-linear map.
We assume that
$f:{\cal E}\to {\cal F}$
is an isomorphism on $X\setminus Z$.
We consider the complex
${\cal K}=[{\cal E}\to {\cal F}]$
of $O_X$-modules
by putting ${\cal F}$ on degree 0.
Then, the localized Chern class
$c^X_Z({\cal K})-1$
is defined as an element of
$CH^*(Z\to X)$ in
\cite{fulton} Chapter 18.1.
We define an element
$c({\cal F}-{\cal E})^X_Z=
(c_i({\cal F}-{\cal E})^X_Z)_{i>0}$
of $CH^*(Z\to X)$ by
\begin{equation}
c({\cal F}-{\cal E})^X_Z=
c({\cal E})\cap 
(c^X_Z({\cal K})-1).
\label{eqchern}
\end{equation}
In other words,
we put
$c_i({\cal F}-{\cal E})^X_Z
=\sum_{j=0}^{\min(d,i-1)}{c_j}({\cal E})
\cap {c_{i-j}}^X_Z({\cal K})$
for $i>0$.
The image of
$c({\cal F}-{\cal E})^X_Z$
in $CH^*(X)$
is the difference $c({\cal F})-c({\cal E})$
of Chern classes.

\begin{lm}\label{lmchern}
Let $X$ be a scheme of finite type
over $F$ and $Z\subset X$
be a closed subscheme.
Let ${\cal E}$
and ${\cal F}$ be
locally free $O_X$-modules
of rank $d$
and $f:{\cal E}\to {\cal F}$
be an $O_X$-linear map
such that
$f:{\cal E}\to {\cal F}$
is an isomorphism on $X\setminus Z$.
Then,

1.  We have
$c_i({\cal F}-{\cal E})^X_Z=0$
for $i>d$.

2. Let
$$\begin{CD}
0@>>>{\cal E}'
@>>>{\cal E}
@>>>{\cal E}''
@>>>0\\
@.@V{f'}VV @VfVV @V{f''}VV @.\\
0@>>>{\cal F}'
@>>>{\cal F}
@>>>{\cal F}''
@>>>0
\end{CD}$$
be a commutative diagram of
exact sequences of locally 
free $O_X$-modules.
We assume that
the maps $f'$ and $f''$
are isomorphism on $X\setminus Z$.
We assume
${\cal E}'$ and ${\cal F}'$
are of rank $d'$
and
${\cal E}''$ and ${\cal F}''$
are of rank $d''$.
Then, we have
\begin{equation}
c({\cal F}-{\cal E})^X_Z=
c({\cal F}'-{\cal E}')^X_Z
\cap c({\cal F}'')+
c({\cal F}''-{\cal E}'')^X_Z
\cap c({\cal E}').
\label{eqchern2}
\end{equation}
\end{lm}

\noindent{\it Proof.}
1. 
The localized Chern classes
${c_i}^X_Z({\cal F})$
and ${c_i}^X_Z({\cal E})$
are defined for $i>d$
in \cite{bloch} \S1.
Further they are equal to 0
since ${\cal F}$ and ${\cal E}$
are locally free of rank $d$.
Hence, by the 
distinguished triangle
$\to {\cal E}\to {\cal F}\to {\cal K}\to $,
we have an equality
$0=
{c_i}^X_Z({\cal F})=
\sum_{j=0}^d{c_j}({\cal E})
\cap {c_{i-j}}^X_Z({\cal K})$
as in Proposition 1.1 (iii) loc.\ cit.
Since the right hand side
is $c_i({\cal F}-{\cal E})^X_Z$,
the assertion follows.

2. We put
${\cal K}'=[{\cal E}'\to {\cal F}']$
and
${\cal K}''=[{\cal E}''\to {\cal F}'']$
as above.
Then, by the assumption,
we have
${c_i}^X_Z({\cal K})=
\sum_{j=0}^{i-1}{c_j}({\cal K}'')
\cap {c_{i-j}}^X_Z({\cal K}')
+
{c_i}^X_Z({\cal K}'')$
for $i>0$
(cf. \cite{fulton} Example 18.1.3, Proposition 18.1 (b) and
\cite{bloch} Proposition 1.1 (iii)).
In other words,
we have
$c^X_Z({\cal K})-1=
(c^X_Z({\cal K}')-1)\cap c({\cal K}'')+
(c^X_Z({\cal K}'')-1)$.
Multiplying
$c({\cal E})=
c({\cal E}')
\cap c({\cal E}'')$
and substituting
$c({\cal E}'')
\cap c({\cal K}'')=
c({\cal F}'')$,
we obtain
$c({\cal E})
\cap (c^X_Z({\cal K})-1)=
c({\cal E}')
\cap (c^X_Z({\cal K}')-1)\cap c({\cal F}'')+
c({\cal E}')
\cap c({\cal E}'')
\cap (c^X_Z({\cal K}'')-1)$
and the assertion follows.
\qed

\begin{lm}\label{lmdis}
Let 
$$\begin{CD}
W@>>> Y\\
@VgVV @VVfV\\
V@>>> X
\end{CD}$$
be a commutative
diagram of separated schemes
of finite type over $F$.
We assume that 
$Y$ is purely of dimension $n$ and
the horizontal arrows
$V\to X$ and $W\to Y$
are regular closed
immersions of codimension $d$.
Let $N_{V/X}$
and $N_{W/Y}$
be the conormal sheaves.

Let $U$ be a
dense open subscheme of $Y$.
We assume that $W\cap U$ is dense in $W$
and that the closed immersion 
$W\cap U\to V\times_XU$
is an open immersion.
We put $Z=W\setminus (W\cap U)$
and $Z'=(V\times_XY)\setminus (W\cap U)$.
Then, we have the following.

1. The canonical map
$g^*N_{V/X}\to N_{W/Y}$
is an isomorphism on $W\cap U=
W\setminus Z$
and 
$c(g^*N_{V/X}-N_{W/Y})^W_Z
\in CH^*(Z\to W)$
is defined.

2. The canonical map $Z_{n-d}(W)\oplus 
CH_{n-d}(Z')\to 
CH_{n-d}(V\times_XY)$
is an isomorphism.
The projection
$CH_{n-d}(V\times_XY)
\to Z_{n-d}(W)$
is given by the restriction map
$CH_{n-d}(V\times_XY)
\to 
CH_{n-d}(W\cap U)=
Z_{n-d}(W\cap U)\simeq
Z_{n-d}(W)$.

3. There exists
a unique element
$[f^!V-W]\in CH_{n-d}(Z')$
satisfying
$[W]+[f^!V-W]=[f^!V]$
in $CH_{n-d}(V\times_XY)$.
Further, we have an equality
\begin{equation}
(W,[f^!V-W])_Y
=
(-1)^{d-1}
c_d(N_{W/Y}-g^*N_{V/X})^W_Z\cap [W]
\label{eqdis2}
\end{equation}
in $CH_{n-2d}(Z)$.
\end{lm}

\noindent{\it Proof.}
1. By the assumption
that the closed immersion 
$W\cap U\to V\times_XU$
is an open immersion,
the canonical map
$g^*N_{V/X}\to N_{W/Y}$
is an isomorphism on 
$W\cap U=W\setminus Z$.
Hence 
$c(g^*N_{V/X}-N_{W/Y})^W_Z
\in CH^*(Z\to W)$
is defined.

2. By the assumption,
the canonical maps $Z_{n-d}(W)\oplus 
Z_{n-d}(Z')\to 
Z_{n-d}(V\times_XY)$
and $Z_{n-d}(W)\to 
CH_{n-d}(W)$
are isomorphisms.
Thus the assertion follows.

3. 
By the assumption,
the restriction of $[f^!V]$
to the open subscheme
$W\cap U\subset V\times_XY$
is $[W\cap U]$.
Hence, by 2,
there exists
a unique element
$[f^!V-W]\in CH_{n-d}(Z')$
satisfying
$[W]+[f^!V-W]=[f^!V]$.

Let $p:Y'\to Y$
be the blow-up at
$V\times_XY\subset Y$ and at $W\subset Y$.
Let
$D=V\times_XY'$
and $D'=W\times_YY'$
be the exceptional divisors.
We compute
$(W,[f^!V-W])_Y$
using $p:Y'\to Y$.
Let $h:D\to V$
and $h':D'\to W$
be the canonical maps
and let 
$N_1={\rm Ker}(h^*N_{V/X}
\to N_{D/Y'})$
and
$N'_1={\rm Ker}(h^{\prime *}N_{W/Y}
\to N_{D'/Y'})$
be the excess conormal sheaves.
The $O_D$-module $N_1$
and the $O_{D'}$-module $N'_1$
are locally free of rank $d-1$.
By the excess intersection formula,
we have
$f^!V=
p_*(V,Y')_X=
(-1)^{d-1}p_*c_{d-1}(N_1)\cap [D]$ and
$W=
p_*(W,Y')_Y=
(-1)^{d-1}p_*c_{d-1}(N'_1)\cap [D']$.

Let $i:D'\to D$
be the immersion.
Since $[f^!V-W]\in CH_{n-d}(Z')$
is characterized by the property
that $[f^!V-W]+[W]=[f^!V]$
in $CH_{n-d}(V\times_XY)$,
we obtain
\begin{align*}
[f^!V-W]
=
(-1)^{d-1}p_*
\left(c_{d-1}(N_1)\cap ([D]-[D'])+
c_{d-1}(N'_1-i^*N_1)^{D'}_{Z_{D'}}\cap [D']\right)
\end{align*}
in $CH_{n-d}(Z)$.
Further by
the excess intersection formula,
we have
\begin{eqnarray}
& &(W,[f^!V-W])_Y\label{eqexc}
\\
&=&
p_*\left(
c_{d-1}(N'_1)\cap [D']\cap \left(c_{d-1}(N_1)\cap  ([D]-[D'])+
c_{d-1}(N'_1-i^*N_1)^{D'}_{Z_{D'}}\cap [D']\right)\right)\nonumber
\end{eqnarray}
in $CH_{n-2d}(Z)$.

Since
$$[D']\cdot ([D]-[D'])
=
([D]-[D'])\cdot [D']
=c_1(N_{D'/Y'}-i^*N_{D/Y'})^{D'}_{Z_{D'}}
\cap [D'],
$$
the right hand side of
(\ref{eqexc})
is equal to
\begin{align*}
p_*&
\Bigl(\left(c_{d-1}(N_1)\cap 
c_1(N_{D'/Y'}-i^*N_{D/Y'})^{D'}_{Z_{D'}}
+
c_{d-1}(N'_1-i^*N_1)^{D'}_{Z_{D'}}\cap c_1(N_{D'/Y'})\right)\cap\\
&\qquad c_{d-1}(N'_1)
\cap [D']\Bigr).
\end{align*}
By the commutative diagram of exact sequences
$$\begin{CD}
0@>>>
i^*N_1@>>>
i^*h^*N_{V/X}@>>>
i^*N_{D/Y'}@>>>0\\
@.@VVV@VVV@VVV @.\\
0@>>>
N'_1@>>>
h^{\prime*}N_{W/Y}@>>>
N_{D'/Y'}@>>>0,
\end{CD}$$
and by Lemma \ref{lmchern}.2,
it is further equal to
\begin{align*}
&p_*(c_d(h^{\prime*}N_{W/Y}-i^*h^*N_{V/X})^{D'}_{Z_{D'}}
\cap c_{d-1}(N'_1)\cap [D'])\\
=&\
c_d(N_{W/Y}-g^*N_{V/X})^W_Z
\cap p_*(c_{d-1}(N'_1)\cap [D']).
\end{align*}
Since $(-1)^{d-1}p_*(
c_{d-1}(N'_1)\cap  [D'])=[W]$,
the assertion follows.
\qed

Let $f:V\to U$ be a finite \'etale morphism
of smooth separated schemes of finite type over $F$
and 
$Y$ be a separated smooth scheme of finite type
containing $V$ as the complement
of a divisor with simple normal crossings.
We put $D^{\log}_{V/U,Y}=
(V\times_UV\setminus \Delta_V,
\Delta_Y)^{\log}_{\mathbb Z}
\in 
CH_0(Y\setminus V)$.
Its image in
$CH_0(Y\setminus V)\otimes_{\mathbb Z}{\mathbb Q}$
is the same as the image of
$D^{\log}_{V/U}$.

\begin{pr}\label{prdis}
Let
$$\begin{CD}
V@>{\subset}>> Y\\
@VfVV @VV{\bar f}V\\
U@>{\subset}>> X
\end{CD}\leqno{\rm (\ref{eqXY})}
$$
be a Cartesian diagram
of separated schemes of finite type over $F$.
We assume
$X$ and $Y$ smooth purely of dimension $d$
over $F$,
$U\subset X$ and $V\subset Y$
are the complements of
divisors $B$ and $D$ with simple normal crossings
respectively and 
$f:V\to U$ is finite and \'etale. 

Then, 
the canonical map
$f^*\Omega^1_{X/F}(\log B)\to
\Omega^1_{Y/F}(\log D)$ is
an isomorphism on $V=Y\setminus D$
and we have
\begin{equation}
D^{\log}_{V/U,Y}
=(-1)^{d-1}
c_d\left(\Omega^1_{Y/F}(\log D)-f^*\Omega^1_{X/F}(\log B)\right)^Y_D\cap [Y].
\label{eqdis}
\end{equation}
\end{pr}

\noindent{\it Proof.}
We consider the commutative diagram
$$\begin{CD}
Y@>>> (Y\times Y)^\sim\\
@VfVV @VV{(f\times f)^\sim}V\\
X@>>> (X\times X)^\sim.
\end{CD}$$
As in \cite{KS} Corollary 4.2.8,
the conormal sheaves
$N_{X/(X\times X)^\sim}$
and 
$N_{Y/(Y\times Y)^\sim}$
are naturally identified with
$\Omega^1_{X/F}(\log B)$ and
$\Omega^1_{Y/F}(\log D)$ respectively.
Hence, it is sufficient to apply Lemma \ref{lmdis}
to the diagram by taking
$V\times V\subset (Y\times Y)^\sim$
as the open subscheme 
$U\subset Y$ in Lemma \ref{lmdis}.
\qed

We compute the log Lefschetz class
assuming an equivariant resolution.
For a closed immersion $Z\to Y$,
let $s(Z/Y)\in \bigoplus_iCH_i(Y)$
be the Segre class.
For a locally free $O_Y$-module ${\cal E}$,
let $c({\cal E})^* =
c({\cal E}^*) =
\sum_i(-1)^ic_i({\cal E})
\in \bigoplus_iCH^i(Y\to Y)$
be the bivariant Chern class
\cite{fulton} Chapter 17.3
of the dual ${\cal E}^*={\cal H}om({\cal E},O_Y)$,
loc.~cit.~Remark 3.2.3 (a).

\begin{lm}\label{lmsgr}
Let $Y$ be a separated 
and smooth scheme
of finite type purely of dimension $d$
over a perfect field $F$
and $V\subset Y$ be the complement of
a divisor $D$ with simple normal crossings.
Let $\sigma$ be an
automorphism of $Y$ over $F$.
We assume that $\sigma$
induces an automorphism of $V$,
$\sigma$ is admissible
and that $V^\sigma=\emptyset$.
Then, we have
\begin{equation}
(\Gamma_\sigma, \Delta_Y)^{\log}_{\mathbb Z}
=\{c(\Omega^1_{Y/F}(\log D))^*\cap
s(Y_{\log}^\sigma/Y)\}_{\dim 0}
\label{eqsgr}
\end{equation}
in $CH_0(Y_{\log}^\sigma)$.
In particular,
if $Y_{\log}^\sigma$
is a Cartier divisor $D_\sigma$ of $Y$,
we have
\begin{eqnarray}
(\Gamma_\sigma, \Delta_Y)^{\log}_{\mathbb Z}
&=&\{c(\Omega^1_{Y/F}(\log D))^*\cap
(1+D_\sigma)^{-1}\cap D_\sigma\}_{\dim 0}\\
&=&(-1)^{d-1}\{c(\Omega^1_{Y/F}(\log D))\cap
(1-D_\sigma)^{-1}\cap D_\sigma\}_{\dim 0}.\nonumber
\label{eqsgrd}
\end{eqnarray}
\end{lm}

\noindent{\it Proof.}
Clear from the definition of
the intersection product \cite{fulton} Proposition 6.1 (a)
and
$N_{Y/(Y\times Y)^\sim}
=\Omega^1_{Y/F}(\log D)$.
\qed

\begin{cor}\label{corsw}
Let $f:V\to U$
be a finite and \'etale morphism
of connected separated and smooth
scheme of finite type 
purely of dimension $d$ over 
a perfect field $F$
and let $\sigma$ be an automorphism of $V$
over $U$ of order $e$.
Let $Y$ be a smooth separated scheme
of finite type over $F$
containing $V$ as the complement of
a divisor $D$ with simple normal crossings.
If $\sigma$ is 
extended to an automorphism of $Y$ over $F$,
the following holds.

1. If $j$ is an integer prime to $e$,
we have $(\Gamma_\sigma, \Delta_Y)^{\log}_{\mathbb Z}
=(\Gamma_{\sigma^j}, \Delta_Y)^{\log}_{\mathbb Z}$
in $CH_0(Y\setminus V)$.

2.  
If $e$ is not a power of $p$,
we have
$(\Gamma_\sigma, \Delta_Y)^{\log}_{\mathbb Z}=0$
in $CH_0(Y\setminus V)$.
\end{cor}

\noindent{\it Proof.}
Let $g:Y'\to Y$ be the blow-up
associated to the subdivision by
baricenters.
Since $g_*:CH_0(Y'\setminus V)\to
CH_0(Y\setminus V)$ is an isomorphism,
by replacing $Y$ by $Y'$,
we may assume
that the action of $\sigma^j$
on $Y$ is admissible for each $j\in {\mathbb Z}$
by Lemma \ref{lmdec}.2.
Then it follows from Lemma \ref{lmsgr}
and Corollary \ref{corsig}.
\qed

\begin{lm}\label{lmord}
Conjecture \ref{cnord} is
true if $\dim U\le 2$.
\end{lm}

\noindent{\it Proof.}
It follows from
Lemma \ref{lmalter2}
and Corollary \ref{corsw}.
\qed

We consider the case of isolated fixed point.

\begin{lm}\label{lmiso}
Let $Y$ be a separated 
and smooth scheme of finite type
purely of dimension $d$ over $F$,
$y$ be a closed point of $Y$
and $\sigma$ be an automorphism
of $Y$ over a perfect field $F$.
Assume that the underlying set of
the fixed part $Y^\sigma$
is $\{y\}$.

Let $f:Y'\to Y$ be
the blow-up at $y$
and $D$ be the exceptional divisor.
Let $g:(Y'\times Y')'\to
(Y'\times Y')$
be the blow-up at $D\times D$.
Then the automorphism $\sigma'$ of $Y'$
induced by $\sigma$
is admissible.
Let $\Gamma'_{\sigma'} \subset 
(Y'\times Y')'$ denote the proper
transform of the graph
$\Gamma_{\sigma'} \subset 
Y'\times Y'$ of $\sigma$
and $\Delta_{Y'}\subset
(Y'\times Y')'$
be the log diagonal.
Then, we have 
\begin{equation}
f_*(\Gamma'_{\sigma'},
\Delta_{Y'})_{(Y'\times Y')'}=
[O_{Y^\sigma}]-[y]
\label{eqarsw}
\end{equation}
in $CH_0(y)={\mathbb Z}$.
\end{lm}

\noindent{\it Proof.}
We have
$[O_{Y^\sigma}]=
(\Gamma_{\sigma},\Delta_{Y})_{Y\times Y}$.
By the projection formula,
we have
$$(\Gamma_{\sigma},\Delta_{Y})_{Y\times Y}=
f_*(g^!(f\times f)^!\Gamma_{\sigma},
\Delta_{Y'})_{(Y'\times Y')'}.$$
Thus it is sufficient to show the equality
$$(g^!(f\times f)^!\Gamma_{\sigma}-
\Gamma'_{\sigma'},
\Delta_{Y'})_{(Y'\times Y')'}=
[y']$$
in $CH_0(D)$ for a $\kappa(y)$-rational point $y'\in D$.

We compute
$g^!(f\times f)^!\Gamma_{\sigma}$.
Since the irreducible components of
$(f\times f)^{-1}
(\Gamma_{\sigma})=
\Gamma_{\sigma}
\times_{Y\times Y}(Y'\times Y')$
are $\Gamma_{\sigma'}$
and $D\times D$,
we have
\begin{align*}
(f\times f)^![\Gamma_{\sigma}]=&\
[\Gamma_{\sigma'}]
+\{c(\Omega^1_{Y/F})^*s(D\times D/Y\times Y)\}_{\dim d}\\
=&\
[\Gamma_{\sigma'}]
+\{(1+D^{(1)})^{-1}(1+D^{(2)})^{-1}D^{(1)}\cdot D^{(2)}\}_{\dim d}.
\end{align*}
Here $D^{(1)}=D\times Y$ and
$D^{(2)}=Y\times D$.
The irreducible components of
$g^{-1}
(\Gamma_{\sigma'})=
\Gamma_{\sigma'}
\times_{Y'\times Y'}(Y'\times Y')'$
are $\Gamma'_{\sigma'}$
and the inverse image $E_D$ of 
the diagonal $D\subset D\times D$.
Hence we have
$$
g^! (f\times f)^![\Gamma_{\sigma}]=
[\Gamma'_{\sigma'}]
+[E_D]
+\{(1+g^*D^{(1)})^{-1}(1+g^*D^{(2)})^{-1}
g^*D^{(1)}\cdot g^*D^{(2)}\}_{\dim d}.$$
Thus we obtain
$$(g^!(f\times f)^!\Gamma_{\sigma}-
\Gamma'_{\sigma'},
\Delta_{Y'})_{(Y'\times Y')'}\\
=
(E_D,
\Delta_{Y'})_{(Y'\times Y')'}
+\{(1+D)^{-2}
D^2\}_{\dim 0}.$$
By
$$(E_D,
\Delta_{Y'})_{(Y'\times Y')'}
=(E_D,\Delta_D)_E
=(\Delta_D,\Delta_D)_{D\times D}
=d[y'],$$
$$\{(1+D)^{-2}
D^2\}_{\dim 0}=(-1)^d(d-1)D^d
=-(d-1)[y'],$$
the assertion follows.
\qed

\section{Swan class 
and Euler characteristic of a sheaf}

We keep the following notation
in this section.
Let $U$ be a connected,
separated and smooth scheme of finite type
purely of dimension $d$ over a perfect field $F$.
Let $\ell$ be a prime number 
different from the characteristic $p$ of $F$.

We consider a smooth $\bar {\mathbb F}_\ell$-sheaf 
${\cal F}$ on $U$
and a finite \'etale Galois covering
$f:V\to U$ trivializing ${\cal F}$.
We define and study
the Swan character class
in \S4.1.
Using it, we define
the Swan classes ${\rm Sw}_{V/U}({\cal F})
\in CH_0(\overline V\setminus V)
\otimes_{\mathbb Z}{\mathbb Q}$ 
and ${\rm Sw}({\cal F})
\in CH_0(\overline U\setminus U)
\otimes_{\mathbb Z}{\mathbb Q}$ 
in \S4.2.
We also prove the formula (\ref{eqgos}) in \S4.2.
In \S4.3, we state an integrality conjecture
(Conjecture \ref{cnint}) that is 
a generalization of the Hasse-Arf theorem
(Lemma \ref{lmdim1}).

\subsection{Swan character class}

We define the Swan character class for 
a ramified Galois covering
using the map
(\ref{eqmap})
$(\ ,\Delta_{\overline V})^{\log}:
CH_d(V\times_UV\setminus \Delta_V)
\to
CH_0(\overline V\setminus V)\otimes_{\mathbb Z}{\mathbb Q}$.

\begin{df}\label{dfSw}
Let $f:V\to U$
be a finite and \'etale Galois covering
of Galois group $G$
of connected separated and smooth 
schemes of finite type
purely of dimension $d$ over 
a perfect field $F$.
For
$\sigma\in G$,
we define the Swan character class 
$s_{V/U}(\sigma)\in 
CH_0(\overline V\setminus V)\otimes_{\mathbb Z}{\mathbb Q}$
by 
\begin{equation}
s_{V/U}(\sigma)=
\begin{cases}
D^{\log}_{V/U}
&\quad\text{ if }\sigma= 1\\
-(\Gamma_\sigma,\Delta_{\overline V})^{\log}
&\quad\text{ if }\sigma\neq 1.
\end{cases}
\label{eqsw}
\end{equation}
\end{df}

If ${\cal C}_V^{\rm sm,0}$
is cofinal in 
${\cal C}_V$,
the Swan character class
$s_{V/U}(\sigma)$
is defined in
$CH_0(\overline V\setminus V)$.

We show basic properties of
Swan character classes.

\begin{lm}\label{lmSwch}
Let the notation be as in Definition
\ref{dfSw}.

1. We have
\begin{equation}
\sum_{\sigma\in G}s_{V/U}(\sigma)=0.
\label{eqSwsum}
\end{equation}

2. If the order of $\sigma$
is not a power of
the characteristic of $F$,
we have
$s_{V/U}(\sigma)=0$.

3.
Let $H\subset G$
be a subgroup 
and $g:V\to U'$ be the 
corresponding intermediate
covering.
Then, for $\sigma\in H$,
we have
\begin{equation}
s_{V/U}(\sigma)=
\begin{cases}
s_{V/U'}(\sigma)&
\text{ if }\sigma\neq 1\\
s_{V/U'}(1)+
{\bar g}^*D_{U'/U}^{\log}&
\text{ if }\sigma= 1.
\end{cases}
\label{eqSwind}
\end{equation}

4. Let $N\subset G$ be a normal subgroup,
and $g:V\to V'$ be
the corresponding intermediate covering.
Then, we have
\begin{equation}
{\bar g}^*s_{V'/U}(\sigma)=
\sum_{\tilde \sigma\in G,\mapsto \sigma}
s_{V/U}(\tilde \sigma)
\label{eqSwchain}
\end{equation}
for $\sigma\in G/N$.
\end{lm}

\noindent{\it Proof.}
1. Clear from the definition and
$V\times_UV\setminus \Delta_V=
\coprod_{\sigma\neq 1}\Gamma_\sigma$.

2. Clear from Proposition \ref{prord}.

3. For $\sigma\neq 1$,
it is clear from Lemma \ref{lmfun}.1.
For $\sigma=1$,
it is nothing but (\ref{eqD}).

4. For $\sigma\neq 1$,
the equality (\ref{eqSwchain})
is clear from Lemma \ref{lmfun}.2.
For $\sigma=1$,
it follows from
the case $\sigma\neq 1$
and the equality (\ref{eqSwsum}).
\qed

\begin{rmk}\label{rmkSwch}
If Conjecture \ref{cnord} is true, we have
$s_{V/U}(\sigma)=s_{V/U}(\sigma^j)$
for an integer $j$ prime to
the order $e$ of $\sigma\in G$.
\end{rmk}

We have the following trace formula.

\begin{lm}\label{lmtrSw}
Let the notation be as in Definition
\ref{dfSw}.
Then, we have
\begin{equation}
{\rm deg}\ s_{V/U}(\sigma)=
\begin{cases}
[V:U]\chi_c(U_{\bar F},{\mathbb Q}_\ell)
-\chi_c(V_{\bar F},{\mathbb Q}_\ell)
\quad &\text{ if }\sigma=1\\
-{\rm Tr}(\sigma^*:H^*_c(V_{\bar F},{\mathbb Q}_\ell))
\quad &\text{ if }\sigma\neq1.
\end{cases}
\label{eqGal}
\end{equation}
\end{lm}

\noindent{\it Proof.}
Clear from the definition 
and Lemma \ref{lmtr2}.
\qed

\begin{cor}\label{cortrSw}
If $j$ is prime to
the order of $\sigma\in G$,
$${\rm deg}\ s_{V/U}(\sigma)=
{\rm deg}\ s_{V/U}(\sigma^j).$$
\end{cor}

\noindent{\it Proof.}
It suffices to consider the case $\sigma\neq 1$.
Since $j$ is prime to the order of $\sigma$,
${\rm Tr}(\sigma^{j*}:H^*_c(V_{\bar F},{\mathbb Q}_\ell))$
is a conjugate of
${\rm Tr}(\sigma^*:H^*_c(V_{\bar F},{\mathbb Q}_\ell))$
over ${\mathbb Q}$.
Hence, by the equality
${\rm deg}\ s_{V/U}(\sigma)=
-{\rm Tr}(\sigma^*:H^*_c(V_{\bar F},{\mathbb Q}_\ell))$,
the degree ${\rm deg}\ s_{V/U}(\sigma^j)$
is a conjugate of
${\rm deg}\ s_{V/U}(\sigma)$
over ${\mathbb Q}$.
Since
${\rm deg}\ s_{V/U}(\sigma)\in {\mathbb Q}$,
the assertion follows.
\qed

If $Y$ is a separated scheme of finite type
containing $U$ as a dense open subscheme,
let $s_{V/U,Y}(\sigma)
\in 
CH_0(Y\setminus V)\otimes_{\mathbb Z}{\mathbb Q}$
denote the image of
$s_{V/U}(\sigma)$.
Let $f:V\to U$
be a finite \'etale Galois covering
of separated smooth schemes
of finite type over $F$.
Let $G$ be the Galois group.
Let $X$ be a normal scheme
containing $U$ as a dense open subscheme
and $Y$ be the normalization of $X$ in $V$.
For a geometric point $\bar y$ of $Y\setminus V$,
let $I_{\bar y}\subset G$ be the inertia
group at $y$.
For a geometric point $\bar x$ of $X\setminus U$,
let $I_{\bar x}\subset G$ be the inertia
group $I_{\bar y}$ at a geometric point
$\bar y$ of $Y\setminus V$ lifting
$\bar x$, that is defined modulo conjugate.

\begin{lm}\label{lmvan}
Let
$f:V\to U$ be a finite \'etale Galois covering
of Galois group $G$
of connected, separated smooth
schemes of finite type 
purely of dimension $d$
over a perfect field $F$.
Let $X$ be a separated normal scheme 
of finite type containing $U$
as a dense open subscheme
and let $Y$ be the normalization of
$X$ in $V$.

Let $\sigma\in G$
be a non-trivial element
and $p$ be the characteristic of $F$.
Assume 
that
$\sigma$ is not in any conjugate of
any $p$-Sylow group of
the inertia subgroup
$I_{\bar x}\subset G$
for any geometric point $\bar x$
of $X\setminus U$.

Then, we have 
$s_{V/U,Y}(\sigma)=0$.
\end{lm}

\noindent{\it Proof.}
If the order of $\sigma$ 
is not a power of $p$,
it follows from
Lemma \ref{lmSwch}.2.
Thus, it suffices to show 
$s_{V/U,Y}(\sigma)=0$
assuming 
$\sigma$ is not in any conjugate of
the $I_{\bar x}\subset G$
for any geometric point $\bar x$
of $X\setminus U$.
The assumption means that
the $\sigma$-fixed part
$Y^{\sigma}$ is empty.
Hence, the assertion
follows from Corollary \ref{coret}.2.
\qed

For an isolated fixed point,
the following is a special case
of a conjecture of Serre.

\begin{cn}\label{cnart}
{\rm (Serre \cite{Sear} (1) p.418)}
Let $Y$ be a separated smooth scheme over 
a perfect field $F$
purely of dimension $d$
and 
$y$ be a closed point of $Y$.
Let $G$ be a finite group
of automorphisms of
$Y$ over $F$
such that,
for $\sigma\neq 1$,
the underlying set of
the fixed part $Y^\sigma$
is $\{y\}$.
Then, the function $a_G:G\to {\mathbb Z}$
defined by
\begin{equation}
a_G(\sigma)=
\begin{cases}
-{\rm length}\ O_{Y^\sigma,y}
&\quad\text{ if }\sigma\neq 1\\
-\sum_{\tau\in G\setminus\{1\}}
a_G(\tau)
&\quad\text{ if }\sigma= 1
\end{cases}
\label{eqart}
\end{equation}
is a character of the group $G$.
\end{cn}

Serre conjectures more precisely that
the character $a_G$
is rational over ${\mathbb Q}_\ell$
for all $\ell\neq p$
in loc.\ cit.\ (2).
Conjecture \ref{cnart} 
is proved in \cite{KSS}
assuming $\dim Y=2$.
In Corollary \ref{corint}.3,
we give a new proof
by deducing it from 
a generalization,
Conjecture \ref{cnint}.1,
assuming $\dim Y=2$.

We compare $a_G(\sigma)$
with the Swan character class
$s_{V/U,Y}(\sigma)$.

\begin{lm}\label{lmart}
Let $Y$ and $G$
be as in Conjecture \ref{cnart}.
We assume the quotient 
$\bar f:Y\to X=Y/G$ exists.
Let $x=\bar f(y)$ be the image of $y$.
Then,
the map
$f:V=Y\setminus\{ y\}\to 
U=X\setminus \{x\}$
is finite \'etale and
$V$ is a Galois covering of Galois group $G$.
Further, for $\sigma\in G$, we have
$$a_G(\sigma)=
\begin{cases}
s_{V/U,Y}(\sigma)-1
&\quad\text{ if }\sigma\neq 1\\
s_{V/U,Y}(\sigma)+|G|-1
&\quad\text{ if }\sigma= 1
\end{cases}$$
in $CH_0(y)={\mathbb Z}$.
\end{lm}

\noindent{\it Proof.}
We keep the notation in the proof of
Lemma \ref{lmiso}.
Then the natural map 
$X'=Y'/G\to X$
is an isomorphism on the complement
$U=X\setminus\{x\}$
and $U$ is a complement of a Cartier divisor of $X'$.
Hence the map
$(\ ,\Delta_Y)^{\log}:
CH_d(V\times_UV\setminus \Delta_V)\to 
CH_0(y)
={\mathbb Z}$
is induced by the intersection product
$(\ ,\Delta_{Y'})_{(Y'\times Y')'}$
and the assertion is clear from Lemma \ref{lmiso}.
\qed

\subsection{Swan class 
and Euler characteristic of a sheaf}

We define the Swan class of an $\bar {\mathbb F}_\ell$-sheaf
${\cal F}$
as a 0-cycle class on the boundary
of a covering trivializing ${\cal F}$.
For a finite group $G$
and a prime number $p$,
let $G_{(p)}\subset G$ be
the set of elements
of order a power of $p$.
If $p=0$,
we put $G_{(p)}=\emptyset$.
For a representation $M$
of $G$ and $\sigma\in G$,
let $M^\sigma$
denote the fixed part
$\{m\in M|\sigma(m)=m\}$.

\begin{df}\label{dfSwV}
Let $U$ be a smooth connected scheme
of dimension $d$ over 
a perfect field $F$ of characteristic $p$
and ${\cal F}$ be a smooth 
$\bar {\mathbb F}_\ell$-sheaf on $U$.
Let $f:V\to U$ 
be a finite \'etale Galois covering
of Galois group $G$ trivializing ${\cal F}$.
Let $M$ be the 
$\bar {\mathbb F}_\ell$-representation
of $G$ corresponding to 
${\cal F}$.

Then, we define the Swan class 
${\rm Sw}_{V/U}({\cal F})
\in 
CH_0(\overline V\setminus V)
\otimes_{\mathbb Z}{\mathbb Q}$
by
\begin{equation}
{\rm Sw}_{V/U}({\cal F})
=\sum_{\sigma \in G_{(p)}}\left( 
\dim_{{\mathbb F}_\ell} M^{\sigma}-
\frac
{\dim_{{\mathbb F}_\ell} M^{\sigma^p}/M^{\sigma}}{p-1}\right)
\cdot
s_{V/U}(\sigma).
\label{eqSwFC}
\end{equation}
\end{df}

Recall that we have $s_{V/U}(\sigma)=0$
if the order of $\sigma$ is not a power of $p$
by Lemma \ref{lmSwch}.2.
Thus we take the sum over $\sigma\in G_{(p)}$.
If $p=0$, we have 
${\rm Sw}_{V/U}({\cal F})=0$.

We define a variant of
the Swan class
expected to be the same as that
defined above.
For an $\bar {\mathbb F}_\ell$-automorphism $\sigma$
of an $\bar {\mathbb F}_\ell$-vector
space $M$ of dimension $m$,
the Brauer trace
${\rm Tr}^{Br}(\sigma:M)
\in {\mathbb Z}[\zeta_\infty]
\subset \bar {\mathbb Q}_\ell$
is defined as follows.
Let $\alpha_1,\ldots,\alpha_m$
be the eigenvalues of $\sigma$
counted with multiplicities
and let $\tilde \alpha_1,\ldots,
\tilde \alpha_m
\in {\mathbb Z}[\zeta_\infty]
\subset \bar {\mathbb Q}_\ell$
be the roots of unity
of order prime to $\ell$ lifting 
$\alpha_1,\ldots,\alpha_m$.
Then, we define
${\rm Tr}^{Br}(\sigma:M)=
\sum_{i=1}^m\tilde \alpha_i$.
If $\sigma$ is an automorphism
of order $p^e$ of $M$,
one can easily verify the equality
\begin{equation}
|({\mathbb Z}/p^e{\mathbb Z})^\times|
\cdot
(\dim_{{\mathbb F}_\ell} M^{\sigma}-
\frac
{\dim_{{\mathbb F}_\ell} M^{\sigma^p}/M^{\sigma}}{p-1})
=\sum_{i\in ({\mathbb Z}/p^e{\mathbb Z})^\times}
{\rm Tr}^{Br}(\sigma^i:M).
\label{eqBr}
\end{equation}

\begin{df}\label{dfSwVn}
Let the notation be as in 
Definition \ref{dfSwV}.
Then, we define the naive Swan class 
${\rm Sw}'_{V/U}({\cal F})
\in 
CH_0(\overline V\setminus V)
\otimes_{\mathbb Z}{\mathbb Q}(\zeta_{p^{\infty}})$
by
\begin{equation}
{\rm Sw}'_{V/U}({\cal F})
=\sum_{\sigma\in G_{(p)}}
s_{V/U}(\sigma)\otimes
{\rm Tr}^{Br}(\sigma:M).
\label{eqSwF}
\end{equation}
\end{df}

\begin{lm}\label{lmSwVn}
Let the notation be as in Definition
\ref{dfSwV}.

1. 
We have
\begin{equation}
{\rm deg}\ {\rm Sw}_{V/U}({\cal F})=
{\rm deg}\ {\rm Sw}'_{V/U}({\cal F}).
\label{eqSwdVn}
\end{equation}

2. If Conjecture \ref{cnord} holds,
we have
\begin{equation}
{\rm Sw}_{V/U}({\cal F})=
{\rm Sw}'_{V/U}({\cal F}).
\label{eqSwVn}
\end{equation}
\end{lm}

\noindent{\it Proof.}
1. It follows 
from the equality (\ref{eqBr})
for an element $\sigma\in G$
of order $p^e$
and Corollary \ref{cortrSw}.

2. It follows 
from the equality (\ref{eqBr})
for an element $\sigma\in G$
of order $p^e$.
\qed

\begin{lm}\label{lmSwV}
Let $f:V\to U$
be a finite and \'etale Galois covering
of connected separated and smooth 
schemes of finite type
purely of dimension $d$ over 
a perfect field $F$ of Galois group $G$.
Let $\ell$ be a prime number 
different form $p={\rm char}\ F$.

1. Let $0\to {\cal F}'\to {\cal F}
\to {\cal F}''\to 0$ be a
short exact sequence of
smooth $\bar {\mathbb F}_\ell$-sheaves on $U$
trivialized on $V$.
Then, 
we have
\begin{equation}
{\rm Sw}_{V/U}({\cal F})
={\rm Sw}_{V/U}({\cal F}')+
{\rm Sw}_{V/U}({\cal F}'').
\label{eqsum}
\end{equation}

2. Let $N\subset G$
be a normal subgroup 
and $g:V\to V'$ be the 
corresponding intermediate
covering.
Let ${\cal F}$ be a smooth 
$\bar {\mathbb F}_\ell$-sheaf on $U$
trivialized on $V'$.
Then, 
we have
\begin{equation}
{\rm Sw}_{V/U}({\cal F})
=g^*{\rm Sw}_{V'/U}({\cal F}).
\label{eqSwVq}
\end{equation}
\end{lm}

\noindent{\it Proof.}
1. Clear from the definition.

2. 
It is clear from Lemma \ref{lmSwch}.4.
\qed

\begin{cor}\label{lmSwan}
Let $U$ be a separated smooth scheme 
of finite type over $F$.
Let ${\cal F}$ be a smooth 
$\bar {\mathbb F}_\ell$-sheaf on $U$.
Let 
$f:V\to U$ 
be a finite \'etale Galois covering
of Galois group $G$
trivializing ${\cal F}$.

1. Then, 
$$\frac 1{|G|}f_*{\rm Sw}_{V/U}({\cal F})
\in 
CH_0(\overline U\setminus U)\otimes_{\mathbb Z}{\mathbb Q}$$
is independent of the choice of $V$.

2. We have
$${\rm Sw}_{V/U}({\cal F})=
\frac1{|G|}f^*f_*{\rm Sw}_{V/U}({\cal F}).$$
\end{cor}

{\it Proof.}
1. 
Clear from Lemma \ref{lmSwV}.2
and Corollary \ref{f*f!}.2.

2.
The Swan class ${\rm Sw}_{V/U}({\cal F})$
is invariant by the Galois group $G$.
Hence it follows from Corollary \ref{f*f!}.2.
\qed

Thus, we define the Swan class ${\rm Sw}({\cal F})$ in
$CH_0(\overline U\setminus U)\otimes_{\mathbb Z}{\mathbb Q}$
as follows.

\begin{df}\label{dfSwan}
Let $U$ be a separated smooth scheme 
of finite type over $F$.
Let ${\cal F}$ be a smooth 
$\bar {\mathbb F}_\ell$-sheaf on $U$.

We define the Swan class 
${\rm Sw}({\cal F})
\in 
CH_0(\overline U\setminus U)\otimes_{\mathbb Z}{\mathbb Q}$
by
\begin{equation}
{\rm Sw}({\cal F})=\frac 1{|G|}
f_*{\rm Sw}_{V/U}({\cal F})
\label{eqSwan}\end{equation}
that is independent of 
a finite \'etale Galois covering
$V\to U$ trivializing ${\cal F}$
by Corollary \ref{lmSwan}.

Similarly, we define the naive Swan class by
\begin{equation}
{\rm Sw}'({\cal F})=\frac 1{|G|}
f_*{\rm Sw}'_{V/U}({\cal F})
\label{eqSwann}\end{equation}
\end{df}

We also define the Swan class for
a smooth $\bar {\mathbb Q}_\ell$-sheaf.

\begin{lm}\label{lmSwan2}
Let $\ell$ be a prime number invertible in $F$.
Assume $U$ is connected.
Let ${\cal F}$ be 
a smooth $\bar {\mathbb Q}_\ell$-sheaf on $U$.
Then the class
${\rm Sw}({\cal F}_0\otimes_{\bar {\mathbb Z}_\ell}
\bar {\mathbb F}_\ell)
\in 
CH_0(\overline U\setminus U)\otimes_{\mathbb Z}{\mathbb Q}$
is indepenent of the choice of
a smooth $\bar {\mathbb Z}_\ell$-sheaf 
${\cal F}_0$ on $U$
satisfying 
${\cal F}={\cal F}_0\otimes_{\bar {\mathbb Z}_\ell}
\bar {\mathbb Q}_\ell$.
\end{lm}

\noindent{\it Proof.}
Clear from Lemma \ref{lmSwV}.1.
\qed

\begin{df}\label{dfSwan2}
Let $\ell$ be a prime number invertible in $F$.
Assume $U$ is connected.
For a smooth $\bar {\mathbb Q}_\ell$-sheaf ${\cal F}$ on $U$,
we define the Swan class 
${\rm Sw}({\cal F})
\in 
CH_0(\overline U\setminus U)\otimes_{\mathbb Z}{\mathbb Q}$
to be the class
${\rm Sw}_{V/U}({\cal F}_0\otimes_{\bar {\mathbb Z}_\ell}
\bar {\mathbb F}_\ell)$
in Lemma \ref{lmSwan2} that is independent of ${\cal F}_0$.
\end{df}

We prove the formula (\ref{eqgos})
for the Euler characteristic.
For a smooth $\bar {\mathbb Q}_\ell$-sheaf ${\cal F}$
on $U$,
we put
$$\chi_c(U_{\bar F},{\cal F})=
\sum_{q=0}^{2d}(-1)^q
\dim_{\bar {\mathbb Q}_\ell}
H^q_c(V_{\bar F},{\cal F}).$$
We define $\chi_c(U_{\bar F},{\cal F})$
similarly for 
a smooth $\bar {\mathbb F}_\ell$-sheaf ${\cal F}$
on $U$.

\begin{thm}\label{thmgos}
Let $U$ be a connected separated smooth scheme of dimension $d$
of finite type over $F$.
Let $\ell$ be a prime number invertible in $F$.
Let ${\cal F}$ be
a smooth $\bar {\mathbb F}_\ell$-sheaf or
a smooth $\bar {\mathbb Q}_\ell$-sheaf on $U$.
Then, we have
$$
\chi_c(U_{\bar F},{\cal F})
={\rm rank}\ {\cal F}\cdot 
\chi_c(U_{\bar F},{\mathbb Q}_\ell)-
{\rm deg}\ {\rm Sw}({\cal F}).
\leqno{(\ref{eqgos})}
$$
\end{thm}

\noindent
{\it Proof.}
It is sufficient to show the case where
${\cal F}$
is a smooth $\bar {\mathbb F}_\ell$-sheaf on $U$.
Let the notation be 
as in Definition \ref{dfSwV}.
Let $G_{\ell\text{-{\rm reg}}}$ be the subset of $G$
consisting of elements of order prime to $\ell$.
By Lemma 2.3 \cite{Il},
we have
$$
\chi_c(U_{\bar F},{\cal F})
=\frac 1{|G|}\sum_{\sigma\in G_{\ell\text{-{\rm reg}}}}
{\rm Tr}(\sigma^*:H^*_c(V_{\bar F},{\mathbb Q}_\ell))
\cdot
{\rm Tr}^{Br}(\sigma:M).$$
By 
Corollary \ref{corvan},
we may replace
$G_{\ell\text{-{\rm reg}}}$
in the summation by $G_{(p)}$.
Thus by Lemma \ref{lmtrSw},
we have
$$
\chi_c(U_{\bar F},{\cal F})
={\rm rank}\ {\cal F}\cdot 
\chi_c(U_{\bar F},{\mathbb Q}_\ell)-
{\rm deg}\ {\rm Sw}'({\cal F})
$$
where ${\rm Sw}'({\cal F})$
is the naive Swan class.
By Lemma \ref{lmSwVn}.1,
we have
${\rm deg}\ {\rm Sw}({\cal F})=
{\rm deg}\ {\rm Sw}'({\cal F})$
and the assertion follows.
\qed

\subsection{Properties of Swan classes}

We keep the notation that
$U$ denotes a connected
smooth scheme purely of dimension $d$ over 
a perfect field $F$
and $\ell$ is a prime number
different from the characteristic of $F$.

We define the wild discriminant and
show the induction formula for
Swan classes.

\begin{df}\label{dfdisc}
Let
$f:V\to U$
be a finite \'etale morphism
of connected, separated and smooth scheme
of finite type
purely of dimension $d$
over $F$.
Then we define the wild discriminant
$d^{\log}_{V/U}
\in 
CH_0(\overline U\setminus U)
\otimes_{\mathbb Z}{\mathbb Q}$
of $V$ over $U$ by 
\begin{equation}
d^{\log}_{V/U}=
f_*
D^{\log}_{V/U}.
\label{eqwdc}
\end{equation}
\end{df}

\begin{lm}\label{lmchaind}
Let
$V\to U'\to U$
be finite \'etale morphism
of separated and smooth schemes
of finite type
purely of dimension $d$
over $F$.
Assume $V\to U'$ is of constant degree
$[V:U']$
and let $h:U'\to U$
denote the map.
Then, we have
\begin{equation}
d^{\log}_{V/U}=
[V:U']\cdot d^{\log}_{U'/U}+
h_*d^{\log}_{V/U'}.
\label{eqdisc}
\end{equation}
\end{lm}

\noindent{\it Proof.}
Clear from Lemma \ref{lmchain}.
\qed

\begin{pr}\label{prind}
Let $f:V\to U$
be a finite and \'etale Galois covering of
connected separated schemes of 
of dimension $d$ of finite type over $F$.
Let $G$ be the Galois group 
and let $h:U'\to U$ be
the intermediate
covering corresponding to 
a subgroup 
$H\subset G$.

Let ${\cal F}$ be a smooth $\bar {\mathbb F}_\ell$-sheaf on $U'$.
Assume that the pull-back
$g^*{\cal F}$ by the map $g:V\to U'$ 
is constant.
Then, if $T\subset G$
is a complete set of representatives of
$G/H$,
we have
\begin{equation}
{\rm Sw}_{V/U}(h_*{\cal F})
=\sum_{\tau \in T}\tau^*({\rm Sw}_{V/U'}({\cal F})+
{\rm rank}\ {\cal F}\cdot g^*D_{U'/U}^{\log}).
\label{eqindV}
\end{equation}
In particular, we have
$$
{\rm Sw}_{V/U}(h_*{\mathbb F}_\ell)
=
\sum_{\tau \in T}\tau^*
g^*D_{U'/U}^{\log}.$$
\end{pr}

\noindent{\it Proof.}
As in Definition \ref{dfSwV},
let $p$ be the characteristic of $F$
and $G_{(p)}\subset G$
be the subset consisting
of elements of order a power of $p$.
Let $M$ be the
$\bar {\mathbb F}_\ell$-representation of $H$
corresponding to
${\cal F}$.
For $\sigma\in G$,
we have 
$$\dim ({\rm Ind}^G_H M)^{\sigma}=
\sum_{\tau\in T}
\frac
{\dim M^{
\langle \tau\sigma\tau^{-1}\rangle\cap H}}
{[\langle \tau\sigma\tau^{-1}\rangle:
\langle \tau\sigma\tau^{-1}\rangle\cap H]}.$$
Thus, we have
\begin{eqnarray*}
&&
\dim_{{\mathbb F}_\ell} ({\rm Ind}^G_H M)^{\sigma}-
\frac{\dim_{{\mathbb F}_\ell} 
({\rm Ind}^G_H M)^{\sigma^p}/
({\rm Ind}^G_H M)^{\sigma}
}{p-1}\\
&=&
\sum_{\tau\in T,\tau\sigma\tau^{-1}\in H_{(p)}}
(\dim_{{\mathbb F}_\ell} M^{\tau\sigma\tau^{-1}}-
\frac{
\dim_{{\mathbb F}_\ell} M^{\tau\sigma^p\tau^{-1}}
/M^{\tau\sigma\tau^{-1}}}{p-1}).
\end{eqnarray*}

Hence, the Swan class
$${\rm Sw}_{V/U}(h_*{\cal F})
=\sum_{\sigma \in G_{(p)}}
(\dim_{{\mathbb F}_\ell} ({\rm Ind}^G_H M)^{\sigma}-
\frac{\dim_{{\mathbb F}_\ell} 
({\rm Ind}^G_H M)^{\sigma^p}/
({\rm Ind}^G_H M)^{\sigma}
}{p-1})
\cdot
s_{V/U}(\sigma)$$
is equal to
\begin{eqnarray*}
&&
\sum_{\sigma \in G_{(p)}}\sum_{\tau\in T,\tau\sigma\tau^{-1}\in H_{(p)}}
(\dim_{{\mathbb F}_\ell} M^{\tau\sigma\tau^{-1}}-
\frac{
\dim_{{\mathbb F}_\ell} M^{\tau\sigma^p\tau^{-1}}
/M^{\tau\sigma\tau^{-1}}}{p-1})\cdot
s_{V/U}(\sigma)\\
&=&
\sum_{\sigma'\in H_{(p)}}\sum_{\tau\in T}
(\dim_{{\mathbb F}_\ell} M^{\sigma'}-
\frac{
\dim_{{\mathbb F}_\ell} M^{\sigma^{\prime p}}
/M^{\sigma'}}{p-1})\cdot
s_{V/U}(\tau^{-1}\sigma'\tau)
\\
&=&
\sum_{\tau\in T}\tau^*\left(\sum_{\sigma'\in H_{(p)}}
(\dim_{{\mathbb F}_\ell} M^{\sigma'}-
\frac{
\dim_{{\mathbb F}_\ell} M^{\sigma^{\prime p}}
/M^{\sigma'}}{p-1})\cdot
s_{V/U}(\sigma')\right).
\end{eqnarray*}
By Lemma \ref{lmSwch}.3,
the content of the big paranthese is equal to
\begin{eqnarray*}
&&
\sum_{\sigma'\in H_{(p)}}
(\dim_{{\mathbb F}_\ell} M^{\sigma'}-
\frac{
\dim_{{\mathbb F}_\ell} M^{\sigma^{\prime p}}
/M^{\sigma'}}{p-1})\cdot
s_{V/U'}(\sigma')
+
\dim M\cdot g^*D_{U'/U}^{\log}\\
&=&
{\rm Sw}_{V/U'}({\cal F})+
{\rm rank}\ {\cal F}\cdot g^*D_{U'/U}^{\log}.
\end{eqnarray*}
Thus the assertion follows.
\qed

\begin{cor}\label{corind}
Let $h:U'\to U$ be a finite and \'etale
morphism of
connected separated schemes of 
of dimension $d$ finite type over $F$.
Let ${\cal F}$ be a smooth $\bar {\mathbb F}_\ell$-sheaf on $U'$.

Then, we have
\begin{equation}
{\rm Sw}(h_*{\cal F})
=h_*{\rm Sw}({\cal F})+
{\rm rank}\ {\cal F}\cdot d_{U'/U}^{\log}.
\label{eqind}
\end{equation}
In particular, we have
$$
{\rm Sw}(h_*{\mathbb F}_\ell)
=
d_{U'/U}^{\log}.$$
\end{cor}

{\it Proof.}
Clear from Proposition \ref{prind}.
\qed

We study the integrality of Swan classes.
For a finite group $G$,
let $C_p(G)$
denote the set of cyclic subgroups $C\subset G$
of order a power of $p$.
For a cyclic subgroup $C\in C_p(G)$,
we put $C^p=\langle \sigma^p\rangle\in C_p(G)$
for a generator $\sigma$ of $C$
and $C^\times=\{\text{\rm generator of }C\}$.
Further, for an $\bar {\mathbb F}_\ell$-representation
$M$ of $G$, we put
$M^C
=\{m\in M|\sigma(m)=m\text{ \rm for all }\sigma \in C\}$.
It is clear that
the product $(\dim_{{\mathbb F}_\ell} M^C-
\frac
{\dim_{{\mathbb F}_\ell} M^{C^p}/M^C}{p-1})\cdot
|C^\times|$
is an integer.

\begin{df}\label{dfSwVZ}
Let $U$ be a smooth connected scheme
of dimension $d$ over 
a perfect field $F$
and ${\cal F}$ be a smooth 
$\bar {\mathbb F}_\ell$-sheaf on $U$.
Let $f:V\to U$ 
be a finite \'etale Galois covering
of Galois group $G$ trivializing ${\cal F}$.
Let $M$ be the 
$\bar {\mathbb F}_\ell$-representation
of $G$ corresponding to 
${\cal F}$.
We assume that ${\cal C}_V^{\rm sm,0}$
is cofinal in ${\cal C}_V$
and that Conjectures \ref{cnord} holds
for $\sigma\in G$.

Then, we define the integral Swan class 
${\rm Sw}_{V/U}({\cal F})_{\mathbb Z}
\in 
CH_0(\overline V\setminus V)$
by
\begin{equation}
{\rm Sw}_{V/U}({\cal F})_{\mathbb Z}
=\sum_{C \in C_p(G)}\left( 
\dim_{{\mathbb F}_\ell} M^C-
\frac
{\dim_{{\mathbb F}_\ell} M^{C^p}/M^C}{p-1}\right)
\cdot |C^\times|
\cdot
s_{V/U}(\sigma_C),
\label{eqSwFCZ}
\end{equation}
where $\sigma_C$ denotes an arbitrary generator of
$C\in C_p(G)$.
\end{df}

The assumptions that ${\cal C}_V^{\rm sm,0}$
is cofinal in ${\cal C}_V$
and that Conjectures \ref{cnord} holds
for $\sigma\in G$
are satisfied if $\dim U\le 2$.

We recall the classical theorem of Hasse-Arf for curves.

\begin{lm}\label{lmdim1}
Let $U$ be a smooth connected curve over 
a perfect field $F$
and ${\cal F}$ be a smooth 
$\bar {\mathbb F}_\ell$-sheaf on $U$
trivialized by a finite \'etale Galois covering
$f:V\to U$ 
of Galois group $G$.
Let $X$ be the proper smooth curve
containing $U$ as a dense open subscheme
and $\bar f:Y\to X$ be the 
normalization in $V$.
We identify
$CH_0(X\setminus U)=
CH_0(\overline U\setminus U)$
and
$CH_0(Y\setminus V)=
CH_0(\overline V\setminus V)$.

Then, 
the integral Swan class
${\rm Sw}_{V/U}({\cal F})_{\mathbb Z}
\in CH_0(Y\setminus V)
=\bigoplus_{y\in Y\setminus V}
{\mathbb Z}\cdot[y]$
is in the image of 
the injection
$\bar f^*:CH_0(X\setminus U)\to 
CH_0(Y\setminus V)$.
\end{lm}

\noindent{\it Proof.}
Since Conjecture \ref{cnord} holds in dimension 1,
the Swan class
${\rm Sw}_{V/U}({\cal F})$
is equal to the naive Swan class
${\rm Sw}'_{V/U}({\cal F})$
by Lemma \ref{lmSwVn}.2.
For $y\in Y\setminus V$,
let $I_y\subset G$
be the inertia group at $y$.
Let $M$ be the corresponding
$\bar {\mathbb F}_\ell$-representation
of $G$.
Then,
by Lemma \ref{lmcurve} and 
\cite{RG},
the Swan conductor
$$
{\rm Sw}_y({\cal F})
=
\frac1{|I_y|}
\sum_{\sigma\in I_y}
s_{V/U,y}(\sigma)
{\rm Tr}^{Br}(\sigma:M)$$
is in ${\mathbb N}$.
For $x\in X\setminus U$,
${\rm Sw}_y({\cal F})$
is independent of the inverse image
$y$ of $x$.
We put
${\rm Sw}_x({\cal F})=
{\rm Sw}_y({\cal F})$
for $x\in X\setminus U$
and 
${\rm Sw}({\cal F})=
\sum_{x\in X\setminus U}
{\rm Sw}_x({\cal F})\cdot [x]
\in CH_0(X\setminus U)$.
Then, we have
$${\rm Sw}'_{V/U}({\cal F})=
\sum_{y\in Y\setminus V}
|I_y|{\rm Sw}_y({\cal F})\cdot [y]=
\bar f^*\sum_{x\in X\setminus U}
{\rm Sw}_x({\cal F})\cdot [x]=
\bar f^*{\rm Sw}({\cal F})$$
and the assertion is proved.
\qed

We expect that Lemma \ref{lmdim1}
holds in higher dimension.

\begin{cn}\label{cnint}
Let $U$ be a smooth connected scheme
of dimension $d$ over 
a perfect field $F$
and ${\cal F}$ be a smooth 
$\bar {\mathbb F}_\ell$-sheaf on $U$.

1. 
The Swan class
${\rm Sw}({\cal F})\in
CH_0(\overline U\setminus U)
\otimes_{\mathbb Z}{\mathbb Q}$
is in the image of
$CH_0(\overline U\setminus U)$.

2. 
Let $f:V\to U$ 
be a finite \'etale Galois covering
trivializing ${\cal F}$.
Assume that 
${\cal C}_V^{\rm sm,0}$
is cofinal in ${\cal C}_V$
and that Conjecture \ref{cnord} holds
as in Definition \ref{dfSwVZ}.

Then, 
the integral Swan class
${\rm Sw}_{V/U}({\cal F})_{\mathbb Z}
\in CH_0(\overline V\setminus V)$
is in the image of 
$f^*:CH_0(\overline U\setminus U)\to 
CH_0(\overline V\setminus V)$.
\end{cn}

Conjecture \ref{cnint}.1
is equivalent to the assertion that
the Swan class
${\rm Sw}_{V/U}({\cal F})\in
CH_0(\overline V\setminus V)
\otimes_{\mathbb Z}{\mathbb Q}$
is in the image of
$f^*:CH_0(\overline U\setminus U)
\to CH_0(\overline V\setminus V)
\otimes_{\mathbb Z}{\mathbb Q}$
for a finite \'etale Galois covering $f:V\to U$
trivializing ${\cal F}$,
by Corollary \ref{lmSwan}.2.

By Lemma \ref{lmdim1},
Conjecture \ref{cnint} is true if $\dim U=1$.
We prove Conjecture \ref{cnint}.1
assuming $\dim U\le 2$
in Corollary \ref{corint}.1.
Conjecture \ref{cnint}.1
is reduced to the rank 1 case by
the induction formula as follows.

\begin{lm}\label{corrk1}
Let $f:V\to U$
be a finite \'etale Galois
covering of Galois group $G$.
We assume that
${\cal C}_{U'}^{\rm sm,0}$
is cofinal
in ${\cal C}_{U'}$
for every intermediate covering
$V\to U'\to U$.
We also assume that
${\rm Sw}\ {\cal G}
\in
CH_0(\overline{U'}\setminus U')
\otimes_{\mathbb Z}{\mathbb Q}$
is in the image of
$CH_0(\overline{U'}\setminus U')$
for every smooth $\bar{\mathbb F}_\ell$-sheaf
of rank 1 on an intermediate covering $U'$
trivialized on $V$.

Then, 
for every smooth $\bar{\mathbb F}_\ell$-sheaf ${\cal F}$
on $U$ trivialized on $V$,
the Swan class ${\rm Sw}\ {\cal F}
\in
CH_0(\overline U\setminus U)
\otimes_{\mathbb Z}{\mathbb Q}$
is in the image of
$CH_0(\overline U\setminus U)$
\end{lm}

\noindent{\it Proof.}
By Brauer's theorem \cite{RG},
we may assume ${\cal F}=h_*{\cal G}$
where $h:U'\to U$
is an intermediate covering
and ${\cal G}$
is a smooth $\bar {\mathbb F}_\ell$-sheaf of rank 1 on
$U'$.
Since 
${\cal C}_{U'}^{\rm sm,0}$
is assumed cofinal
in ${\cal C}_{U'}$,
the wild different
$D_{U'/U}^{\log}$
is defined in 
$CH_0(\overline {U'}\setminus U')$
by Proposition \ref{prdis}.
Hence,
the wild discriminant
$d_{U'/U}^{\log}$
is in the image of
$CH_0(\overline U\setminus U)$.
Thus it follows from
the assumption that
${\rm Sw}\ {\cal G}$
is in the image
of $CH_0(\overline {U'}\setminus U')$
and the induction formula
Corollary \ref{corind}.
\qed

If $X$ is a separated scheme of finite type
containing $U$ as a dense open subscheme,
let ${\rm Sw}_X({\cal F})
\in 
CH_0(X\setminus U)\otimes_{\mathbb Z}{\mathbb Q}$
denote the image of
${\rm Sw}({\cal F})$.
Similarly, if $Y$ is a separated scheme of finite type
containing $V$ as a dense open subscheme,
let ${\rm Sw}_{V/U,Y}({\cal F})
\in 
CH_0(Y\setminus V)\otimes_{\mathbb Z}{\mathbb Q}$
denote the image of
${\rm Sw}_{V/U}({\cal F})$.

\begin{lm}\label{lmartint}
Conjecture \ref{cnint}.1
implies Conjecture \ref{cnart}.
\end{lm}

\noindent
{\it Proof.}
Let the notation be 
as in Conjecture \ref{cnart}.
Since $|G|a_G$ is a character of $G$
by \cite{Sear} Proposition 7,
it is sufficient to show that
the Artin conductor
\begin{equation}
a_G(M)=\frac 1{|G|}\sum_{\sigma\in G}
a_G(\sigma){\rm Tr}(\sigma:M)
\label{eqarc}
\end{equation}
defined in ${\mathbb Q}$
is in ${\mathbb Z}$
for every $\bar {\mathbb Q}_\ell$-representation 
$M$ of $G$.
We may assume $Y$ is affine
and the quotient $X=Y/G$ exists.
Let $x\in X$ be the image of $y$
and 
${\cal F}$
be the smooth sheaf on $U=X
\setminus\{x\}$ corresponding
to the representation $M$.
Then, by Corollary \ref{lmart}
and Corollary \ref{corsw}.2, 
we have
$a_G(M)={\rm Sw}_X({\cal F})+\dim M-\dim M^G$
in $CH_0(x)\otimes_{\mathbb Z}{\mathbb Q}
={\mathbb Q}$.
Thus the assertion is proved.
\qed

We give a refinement of
Th\'eor\`eme 2.1 of \cite{Il}.

\begin{lm}\label{lmIl}
Let the notation be
as in Lemma \ref{lmvan}.
Let $p$ be the characteristic of 
a perfect field $F$.
Let ${\cal F}_1$ and ${\cal F}_2$ be
smooth $\bar {\mathbb F}_\ell$-sheaves on $U$
corresponding to
$\bar {\mathbb F}_\ell$-representations
$M_1$ and $M_2$ of $G$.
Assume that $X$ is normal
and that, for each geometric point
$\bar x$ of $X\setminus U$,
the restrictions of
$M_1$ and $M_2$ to 
a $p$-Sylow subgroup of
the inertia
subgroup $I_{\bar x}$
are isomorphic to each other.

Then, we have
$${\rm Sw}_{V/U,Y}({\cal F}_1)
={\rm Sw}_{V/U,Y}({\cal F}_2)$$
\end{lm}

\noindent{\it Proof.}
Clear from Lemma \ref{lmvan}
and Definition \ref{dfSwV}.
\qed

If the base field is finite,
we expect to have
the following refinement
of Theorem \ref{thmgos}.

\begin{cn}\label{cndet}
Let $U$ be a connected separated smooth scheme of dimension $d$
of finite type over a finite field $F$.
Let $X$ be a proper normal scheme
over $F$ containing $U$ as
a dense open subscheme.
Let $Fr_F\in {\rm Gal}(\bar F/F)$
be the geometric Frobenius
and let $\rho_X:CH_0(X)
\to \pi_1(X)^{\rm ab}$
be the reciprocity map sending
$[x]$ to the geometric Frobenius $Fr_x$
for closed points $x\in X$.

Let $\ell$ be a prime number invertible in $F$.
Let ${\cal F}$ be
a smooth $\bar {\mathbb F}_\ell$-sheaf or
a smooth $\bar {\mathbb Q}_\ell$-sheaf on $U$.
We assume
Conjecture \ref{cnint}.1 holds
and ${\rm Sw}_X({\cal F})
\in CH_0(X\setminus U)$
is defined.

Let ${\cal G}$ be a smooth
$\bar {\mathbb F}_\ell$-sheaf
or 
$\bar {\mathbb Q}_\ell$-sheaf
on $X$ and let
$\det{\cal G}:
\pi_1(X)^{\rm ab}
\to
\bar {\mathbb F}_\ell^\times$
or
$\pi_1(X)^{\rm ab}
\to
\bar {\mathbb Q}_\ell^\times$
be the character corresponding
to the smooth sheaf 
$\det{\cal G}$ of rank 1.
We put
$\det(-Fr_F:H^*_c(U_{\bar F},{\cal F}))
=
\prod_{q=0}^{2d}
\det(-Fr_F:H^q_c(U_{\bar F},{\cal F}))^{(-1)^q}.$

Then, we have
$$
\det(-Fr_F:H^*_c(U_{\bar F},{\cal F}\otimes {\cal G}))
=
\det(-Fr_F:H^*_c(U_{\bar F},{\cal F}))^{{\rm rank}\ {\cal G}}
\cdot \det{\cal G}(\rho_X({\rm Sw}_X({\cal F}))).$$
\end{cn}

If $\dim U=1$,
Conjecture \ref{cndet}
is a consequence of the product
formula for the constant term of
the functional equation of
$L$-functions
\cite{De}, \cite{La}.

\section{Computations of Swan classes}

We compare
the Swan classes
${\rm Sw}({\cal F})$
of sheaves of rank 1
with an invariant
defined in \cite{Kato}
in \S5.1.
Using the computation,
we prove the integrality conjecture
Conjecture \ref{cnint}.1
assuming $\dim U\le 2$.
We also compare the formula
(\ref{eqgos}) with Laumon's formula in \cite{L}.

We keep the notation that
$U$ denotes a connected
smooth scheme purely of dimension $d$ over 
a perfect field $F$
and $\ell$ is a prime number
different from the characteristic $p$ of $F$.

\subsection{Rank 1 case}

Let $X$ be a smooth separated scheme 
of finite type purely of dimension $d$
over $F$
and $U\subset X$ be the complement
of a divisor $D$ with simple normal crossings.
Let $\ell$ be a prime number
invertible in $F$.
We identify $\mu_p(\bar {\mathbb F}_\ell)=
{\mathbb Z}/p{\mathbb Z}$.

Let
${\cal F}$ be a smooth 
$\bar {\mathbb F}_{\ell}$-sheaf of rank 1
on $U$.
We briefly recall the definition of
the 0-cycle class 
$c_{\cal F}$ in \cite{Kato}.
Let $D_1,\ldots,D_m$
be the irreducible components of $D$. 
Let $\chi\in H^1(U,
\bar {\mathbb F}_{\ell}^\times)$
be the element
corresponding to ${\cal F}$.
In loc.\ cit.,
the Swan divisor $D_\chi=
\sum_{i=1}^m{\rm sw}_i({\chi})D_i\ge 0$
is defined.
Also the refined Swan character map
$${\rm rsw}_i(\chi):
O(-D_{\chi})|_{D_i}\to \Omega^1_{X/F}(\log D)|_{D_i}$$
is defined
for each irreducible component $D_i$
such that ${\rm sw}_i(\chi)>0$.

We put $E=\sum_{i;{\rm sw}_i(\chi)>0}D_i\subset D$.
If $\bar f:Y\to X$ is the normalization
in the cyclic \'etale covering $f:V\to U$
corresponding to $\chi$,
the closed subscheme
$E\subset X$ is the wild ramification locus of the covering
$Y\to X$.
The sheaf ${\cal F}$ is said to be clean
with respect to $X$
if the map
${\rm rsw}_i(\chi):
O(-D_{\chi})|_{D_i}\to \Omega^1_{X/F}(\log D)|_{D_i}$
is a locally splitting injection
for each component $D_i$ of $E$.
If ${\cal F}$ is clean
with respect to $X$,
the 0-cycle class $c_{\cal F}=c_\chi\in CH_0(E)$
is defined by
\begin{eqnarray}
c_{\cal F}=c_\chi&=&
\{c(\Omega^1_{X/F}(\log D))^*\cap
(1+D_\chi)^{-1}\cap D_\chi\}_{\dim 0}
\label{eqswchi}\\
&=&
(-1)^{d-1}\sum_{i=1}^m{\rm sw}_i(\chi)
c_{d-1}({\rm Coker}({\rm rsw}_i(\chi)))\cap [D_i].\nonumber
\end{eqnarray}
If one want to specify $X$,
we write $c_{{\cal F},X}$
for $c_{\cal F}$.

\begin{cn}\label{cncF}
Let $X$ be a separated
scheme of finite type
over a perfect field $F$
and $U\subset X$ be a dense open subscheme of $X$.
Assume $U$ is connected
and smooth purely of dimension $d$ over $F$.
Let $\ell$ be a prime number invertible in $F$
and ${\cal F}$ be a smooth $\bar {\mathbb F}_\ell$-sheaf 
of rank 1 on $U$.

1. 
Let
$$\begin{CD}
V@>{\subset}>> Y\\
@VfVV @VV{\bar f}V\\
U@>{\subset}>> X
\end{CD}\leqno{\rm (\ref{eqXY})}
$$
be a Cartesian diagram
of smooth separated schemes of finite type over $F$.
We assume $U\subset X$
and $V\subset Y$
are the complement of divisors
with simple normal crossings,
$f:V\to U$ is a connected finite \'etale Galois covering
of Galois group $G$
and ${\cal F}$ is constant on $V$.
If ${\cal F}$ is clean with respect to $X$,
we have
\begin{equation}
{\rm Sw}_{V/U,Y}({\cal F})
=\bar f^*c_{{\cal F},X}
\label{eqrk1}
\end{equation}
in $CH_0(E\times_XY)\otimes_{\mathbb Z}
{\mathbb Q}$.

2. There exists a 
Cartesian diagram
$$\begin{CD}
U@>{\subset}>> X'\\
@| @VV{\bar f}V\\
U@>{\subset}>> X
\end{CD}
$$
satisfying the following conditions:
the map $\bar f:X'\to X$ is proper,
$X'$ is smooth over $F$,
$U$ is the complement of
a divisor with simple normal crossings in $X'$
and ${\cal F}$ is clean with respect to $X'$.
\end{cn}

Conjecture \ref{cncF}.2
is proved if $\dim U\le 2$
in \cite{Kato} Theorem 4.1.
We prove Conjecture \ref{cncF}.1
assuming $\dim U\le 2$
later in Theorem \ref{thmint}.

\begin{lm}\label{lmrk1}
Conjecture \ref{cncF}
implies
Conjecture \ref{cnint}.1.
\end{lm}

\noindent{\it Proof.}
Note that Conjecture \ref{cncF}.2
is stronger than
the strong resolution of singularities.
Hence the assertion follows from Lemma \ref{corrk1}.
\qed

We prove Conjecture \ref{cncF}.1
in some cases.
We say
${\cal F}$ is $s$-clean
with respect to $X$
if it is clean and further if
the composition
$$\begin{CD}
O(-D_{\chi})|_{D_i}
@>{{\rm rsw}_i(\chi)}>>
\Omega^1_{X/F}(\log D)|_{D_i}
@>{{\rm res}_{D_i}}>> 
O_{D_i}\end{CD}$$
is either an isomorphism or
the 0-map 
for each component $D_i$ of $E$,
depending on $D_i$.
We recall results
in \cite{Sa}.

\begin{lm}\label{lmSa}
{\rm (\cite{Sa} Lemmas 1 and 2)}
Let $p>0$ be the characteristic of $F$
and $e\ge 1$ be an integer.
Let $X$ be a separated and smooth
scheme of finite type over $F$,
$U$ be the complement of
a divisor with simple normal crossings.
Let $f:V\to U$ be 
a finite \'etale connected cyclic
covering of degree $p^e$
and let $f_1:U_1\to U$ be
the intermediate covering
of degree $p$.
Let ${\cal F}$ 
and ${\cal G}$ be the smooth
$\bar {\mathbb F}_{\ell}$-sheaves of rank 1
corresponding to characters
$\chi,\theta:{\rm Gal}(V/U)\to 
\bar {\mathbb F}_{\ell}^\times$
of degree $p^e$ and $p$ respectively.
We assume that
the sheaf ${\cal G}$ is
$s$-clean with respect to $X$.
Let $E\subset X$ be the union of
irreducible components of $X\setminus U$
where ${\cal F}$ has wild ramification.

Then, there exists a Cartesian diagram
$$\begin{CD}
U_1@>{\subset}>> Y_1\\
@V{f_1}VV @VV{\bar f_1}V\\
U@>{\subset}>> X
\end{CD}
$$
of smooth separated scheme of finite
type satisfying the following condition:
\begin{enumerate}
\item[{\rm [\ref{lmSa}]}]
The map $\bar f_1:Y_1\to X$ is proper
and 
$U_1\subset Y_1$ is the complement of a
divisor with simple normal crossings.
If $\sigma$ is a generator of ${\rm Gal}(U_1/U)$, 
the action of $\sigma$ on $U_1$
is extended to an admissible action on $Y_1$
over $X$
and we have
\begin{equation}
p\cdot s_{U_1/U}(\sigma)
=-\bar f_1^* c_{\cal G}
\label{eqSa}
\end{equation}
in $CH_0(Y_1\times_XE)$.
Further if ${\cal F}$
is clean with respect to $X$
and if ${\cal F}_1=f_1^*{\cal F}$
is clean with respect to $Y_1$,
we have
\begin{equation}
\bar f_1^*c_{\cal F}=
c_{{\cal F}_1}+D^{\log}_{U_1/U}.
\label{eqSa2}
\end{equation}
\end{enumerate}
\end{lm}

\begin{pr}\label{prrk1}
Let the notation be as in Conjecture \ref{cncF}.1.
Let $\chi:G\to 
\bar {\mathbb F}_{\ell}^\times$
be the character 
corresponding to ${\cal F}$.
Let $n$ be the order of $\chi$ and
$e={\rm ord}_pn$ be
the $p$-adic valuation.
For $0\le i\le e$,
let $U_i$ be the intermediate \'etale covering
corresponding to the subgroup
$G_i\subset G$ of index $p^i$.

We assume that the diagram {\rm (\ref{eqXY})}
is inserted in a Cartesian diagram
\begin{equation}
\begin{CD}
V@>h>>U_e@>>> \cdots @>>> U_{i+1}@>>> U_i @= U_i@>>>\cdots @= U_0&=U\\
@V{\cap}VV@V{\cap}VV   @.    @V{\cap}VV @V{\cap}VV@V{\cap}VV @. @VVV\\
Y@>{\bar h}>>X_e@>{\bar g_e}>>\cdots 
@>{\bar g_{i+1}}>> Y_{i+1} @>{\bar f_{i+1}}>>X_i@>{\bar g_i}>> Y_i
@>{\bar f_i}>>\cdots @>{\bar g_0}>>Y_0&=X.
\end{CD}\label{eqAS}
\end{equation}
satisfying the following conditions
{\rm (\ref{eqAS}.1)-(\ref{eqAS}.3)}:

{\rm (\ref{eqAS}.1)} For $0\le i\le e$,
$X_i$ and $Y_i$ 
are separated and smooth over $F$ and 
contain $U_i$ as
the complement of divisors with simple normal crossings.
The pull-back ${\cal F}_i={\cal F}|_{U_i}$ 
is clean with respect to $X_i$
and to $Y_i$
and we have $\bar g_i^*(c_{{\cal F}_i,Y_i})=
c_{{\cal F}_i,X_i}$.

{\rm (\ref{eqAS}.2)}
For $0\le i< e$,
the smooth 
$\bar {\mathbb F}_{\ell}$-sheaf
${\cal G}_i$ on $U_i$ corresponding to
a non-trivial character
${\rm Gal}(U_{i+1}/U_i)
\to \bar {\mathbb F}_{\ell}^\times$
is $s$-clean with respect to $X_i$
and $\bar f_{i+1}:Y_{i+1}\to X_i$ satisfies
the condition {\rm [\ref{lmSa}]} in Lemma \ref{lmSa}.

{\rm (\ref{eqAS}.3)}
The actions of $G$ on $U_1,\ldots, U_e$ and on $V$
are extended to admissible actions on 
$X_1,\ldots, X_e$ and on $Y$.

\noindent
Then, we have
$${\rm Sw}_{V/U,Y}({\cal F})_{\mathbb Z}
=\bar f^*c_{{\cal F},X}
\leqno{\rm(\ref{eqrk1})}
$$
in $CH_0(E\times_XY)$.
\end{pr}

\noindent{\it Proof.}
First, we reduce it to the
case where $n=p^e$.
We decompose 
$G={\rm Gal}(V/U)
={\rm Gal}(U_e/U)
\times {\rm Gal}(U'/U)$
to the $p$-part 
${\rm Gal}(U_e/U)$
and the non-$p$-part
${\rm Gal}(U'/U)$.
Let $\chi'$ be the restriction to
the $p$-part ${\rm Gal}(U_e/U)$
and let ${\cal F}'$
be the corresponding sheaf on $U$.
By the definition
in \cite{Kato},
we have
$c_{\cal F}=
c_{{\cal F}'}$.
By Lemma \ref{lmSwch}.2,
we have
${\rm Sw}_{V/U,Y}({\cal F})_{\mathbb Z}
={\rm Sw}_{V/U,Y}({\cal F}')_{\mathbb Z}$.
By Lemma \ref{lmSwV}.2,
we have
${\rm Sw}_{V/U,Y}({\cal F}')_{\mathbb Z}
=\bar h^*
{\rm Sw}_{V_e/U,X_e}({\cal F}')_{\mathbb Z}$.
Thus the assertion is reduced
to the case where $n$ is a power of $p$.

We assume $n=p^e$ and
prove the assertion by induction on $e$.
We prove the case $n=p$.
By the condition (\ref{eqAS}.3) and 
Corollary \ref{corsw}.1,
we have ${\rm Sw}_{V/U,Y}({\cal F})_{\mathbb Z}
=-p\cdot s_{V/U}(\sigma)$
for a generator $\sigma$
of ${\rm Gal}(V/U)$.
Hence the assertion follows from
the equality (\ref{eqSa})
in Lemma \ref{lmSa} and
the assumption $\bar g_0^*c_{{\cal F},Y_0}=
c_{{\cal F},X_0}$
(\ref{eqAS}.1)
in the case $n=p$.

We assume $e\ge 2$.
By the induction hypothesis,
we may assume
${\rm Sw}_{V/U_1,Y}({\cal F}_1)_{\mathbb Z}=
\bar g^*c_{{\cal F}_1}$
where $\bar g:Y\to Y_1$
denotes the composition.
By the equality (\ref{eqSa2})
in Lemma \ref{lmSa} and
the assumption $\bar g_0^*c_{{\cal F},Y_0}=
c_{{\cal F},X_0}$
(\ref{eqAS}.1),
we have
$\bar f^*c_{\cal F}=
\bar g^*c_{{\cal F}_1}+\bar g^*D^{\log}_{U_1/U}$
in $CH_0(Y\setminus V)$.
By the condition (\ref{eqAS}.3), 
Corollary \ref{corsw}.1
and Lemma \ref{lmSwch}.3,
we have
${\rm Sw}_{V/U,Y}({\cal F})_{\mathbb Z}=
{\rm Sw}_{V/U_1,Y}({\cal F}_1)_{\mathbb Z}+
\bar g^*D^{\log}_{U_1/U}$.
Thus the assertion is proved.
\qed

\begin{thm}\label{thmint}
Conjecture \ref{cncF}.1
is true if $\dim U\le 2$.
More precisely,
we have
$${\rm Sw}_{V/U,Y}({\cal F})_{\mathbb Z}
=\bar f^*c_{{\cal F},X}
\leqno{\rm(\ref{eqrk1})}
$$
in $CH_0(E\times_XY)$.
\end{thm}

\noindent{\it Proof.}
Without loss of generality,
we may assume $X$ and $Y$
are proper over $F$,
since the strong resolution is known in
dimension $\le 2$.
If $\dim U=1$,
we obtain a diagram (\ref{eqAS})
satisfying the conditions 
{\rm (\ref{eqAS}.1)-(\ref{eqAS}.3)}
in Proposition \ref{prrk1}
by taking the normalizations
$X_i=Y_i$ of $X$ in $U_i$
and the assertion follows.

To prove the case $\dim U=2$,
first we recall some results
from \cite{Kato}.

\begin{lm}\label{lmint2}
Let $X$ and $X'$
be smooth surfaces of finite type over $F$
containing $U$ as the complement
of divisors with simple normal crossings
and $g:X'\to X$ be a morphism over $F$
inducing the identity on $U$.
Let ${\cal F}$ be a smooth $\bar {\mathbb F}_\ell$-sheaf
of rank 1 on $U$
clean with respect to $X$.

1.  The sheaf ${\cal F}$ is also clean
with respect to $X'$ and 
we have $g^*c_{{\cal F},X}=c_{{\cal F},X'}$.

2. Assume ${\cal F}$ corresponds to 
a character of order $p$.
Then ${\cal F}$ is $s$-clean 
with respect to the complement of
at most finitely many closed points of
$X\setminus U$.
If $g:X'\to X$ is the blow-up at
the points where ${\cal F}$ is 
not $s$-clean,
then ${\cal F}$ is $s$-clean with respect to $X'$.
\end{lm}

\noindent{\it Proof.}
1. It is sufficient to consider the case
where $g:X'\to X$ is the blow-up
at a closed point
of the complement of $U$.
Then, 
${\cal F}$ is clean 
with respect to $X'$
by \cite{Kato} Remark 4.13.
Further, 
we have $c_{{\cal F},X}=g_*c_{{\cal F},X'}$
by \cite{Kato} Theorem 5.2.
Hence by Lemma \ref{lmCV},
we have $g^*c_{{\cal F},X}=c_{{\cal F},X'}$.

2.
The first assertion is clear from the definition
of $s$-cleanness. 
We show the second assertion.
We may assume 
${\cal F}$ is $s$-clean with respect to
$X\setminus \{x\}$ where
$x\in X\setminus U$
is a closed point.
Then, the characterization given in 
\cite{Kato} (3.6) shows that
${\cal F}$
is defined by an Artin-Schreier equation
$T^p-T=s/t^n$ where
$(s,t)$ is a local coordinate at $x$
and $n$ is prime to $p={\rm char}\ F$,
on an \'etale neighborhood of $x$.
(In \cite{Kato} p. 773, $h=gf$ in line 7
should read $h=g^{-1}f$
and $\pi_{1\le i\le r}\pi_i$ in line 12
should read $\prod_{1\le i\le r}\pi_i$.)
Then the assertion is easily checked.
\qed

We go back to the proof of Theorem \ref{thmint}
in the case $\dim U=2$.
By Lemma \ref{lmCV},
we may replace $Y$
by a successive blow-up
$Y'\to Y$
at closed points in the complement of $V$.
By Lemma \ref{lmint2},
we may also replace $X$
by a successive blow-up
$X'\to X$
at closed points in the complement of $U$.
By Proposition \ref{prrk1},
it is sufficient to 
construct a diagram (\ref{eqAS})
satisfying the conditions 
{\rm (\ref{eqAS}.1)-(\ref{eqAS}.3)}
after possibly replacing $X$ and $Y$ by
successive blow-ups 
at closed points in the complements.

For $0\le i< e$,
there exist a proper smooth surface
$X'_i$ containing $U_i$
as the complement of a
divisor with simple normal crossing
such that ${\cal F}_i$
and ${\cal G}_i$
are clean with respect to $X'_i$
and that the map $U_i\to U$
is extended to a map $X'_i\to X$
for $0\le i<e$,
by \cite{Kato} Theorem 4.1.
By Lemma \ref{lmalter2},
there exist a proper smooth surface
$Y'$ containing $V$
as the complement of a
divisor with simple normal crossing
such that 
the maps $V\to U_i$
are extended to maps $Y'\to X'_i$
and that the action of $G$ on $V$
is extended to an admissible action on $Y'$
over $X$.

We define a diagram 
(\ref{eqAS})
satisfying 
the conditions 
{\rm (\ref{eqAS}.1)-(\ref{eqAS}.3)}
in Proposition \ref{prrk1}
inductively
after possibly replacing $X$ and $Y$
by a successive blow-up.
Applying Lemma \ref{lmalter2}
to the quotient $Y'/G$,
we obtain a proper smooth surface
$Y_0$ containing $U_0$
as the complement of a divisor
with simple normal crossing
with a map 
$Y_0\to Y'/G$
extending the identity of $U$.
Since the identity of $U=U_0$
is extended to a map
$Y'/G\to X'_0$,
the identity of $U_0$
is extended to a map
$Y_0\to X'_0$.
By replacing $X$ by $Y_0$,
we put $X=Y_0$.

We define $Y_{i+1}$ and $X_i$
inductively by
assuming that $Y_i$ is already
defined and that 
the identity of $U_i$
is extended to a map
$Y_i\to X'_i$.
Applying Lemma \ref{lmalter2}
to $Y_i$,
we obtain a proper smooth scheme $Y'_i$
that contains $U_i$
as the complement of a divisor with simple normal crossings
and that the action of $G$ on $U_i$
is extended to an admissible action on $Y'_i$ over $X$.
Since $Y'_i$ dominates $X'_i$,
the sheaves ${\cal F}_i$
and ${\cal G}_i$
are clean with respect to $Y'_i$,
by Lemma \ref{lmint2}.1.
Let $X_i\to Y'_i$ 
be the blowing-up 
at the closed points where
${\cal G}_i$ are not $s$-clean
and $\bar g_i:X_i\to Y_i$
be the composition.
Then 
the sheaf ${\cal G}_i$
is $s$-clean with respect to $X_i$
by Lemma \ref{lmint2}.2.
Further ${\cal F}_i$
is clean with respect to $X_i$
and the condition $\bar g_i^*c_{{\cal F}_i,Y_i}=
c_{{\cal F}_i,X_i}$ is satisfied
by Lemma \ref{lmint2}.1.
Applying Lemma \ref{lmSa},
we obtain $Y_{i+1}\to X_i$.

By the construction,
we see that $Y_{i+1}$ dominates $X'_{i+1}$.
Repeating this construction inductively,
we obtain a diagram
(\ref{eqAS})
except the map $\bar h:Y\to X_e$.
We define $Y''$ by applying
the construction in Lemma \ref{lmalter2}
to the normalization of $X_e$ in $V$.
Then the action of $G$ on $V$ is
extended to an admissible action on $Y''$ over $X$.
Replacing $Y$ by $Y''$,
we obtain a diagram
(\ref{eqAS})
satisfying the conditions
{\rm (\ref{eqAS}.1)-(\ref{eqAS}.3)} in
Proposition \ref{prrk1}.
Thus the assertion is proved.
\qed

\begin{cor}\label{corint}
1. Conjecture \ref{cnint}.1
is true if $\dim U\le 2$.

2. {\rm (\cite{KSS})} Conjecture \ref{cnart}
is true if $\dim Y\le 2$.
\end{cor}

\noindent{\it Proof.}
Clear from Lemmas \ref{lmrk1}
and \ref{lmartint}
respectively.
\qed

\subsection{Comparison with
Laumon's formula}

In \cite{L}, Laumon proves a
generalization of the
Grothendieck-Ogg-Shafarevich formula for surfaces
under the assumption (NF) below on ramification. 
We compare the formula
(\ref{eqgos}) with Laumon's formula in \cite{L}.

For simplicity, we assume $F$
is an algebraically closed field.
Let $X$ be a proper
normal connected surface over $F$ and
$U$ be a smooth dense open subscheme. Let
${\cal F}$ be a smooth ${\bar {\mathbb
F}}_{\ell}$-sheaf on $U$. 
Let $B_1,\ldots,B_m$ be the irreducible components
of dimension 1
of the complement
$B=X\setminus U$,
let $\xi_i$ be the generic point of
$B_i$, and let $K_i$ be the 
field of fractions of
the completion of $O_{X,\xi_i}$. We assume 
the following condition.

\medskip

(NF) For each $i$, the finite Galois
extension of $K_i$ that trivializes ${\cal F}$
has {\it separable} residue extension.

\medskip
\noindent
By this assumption, the Swan conductor
$\text{Sw}_i({\cal F})\in {\mathbb N}$ of
$\cal F$ for the local field $K_i$ is defined
by the classical ramification theory,
as in the proof of Lemma \ref{lmdim1}.
In \cite{L},
a smooth dense open subscheme
$B_i^\circ\subset B_i$ 
for each $i$
and an integer
${\rm Sw}_x({\cal F})\in {\mathbb Z}$ for each closed point 
$x\in \Sigma=B\setminus \bigcup_{i=1}^mB_i^\circ$
are defined
and the formula
\begin{equation}
\chi_c(U, {\cal F})= 
{\rm rank}\ {\cal F}\cdot 
\chi_c(U,  {\mathbb Q}_\ell) - 
\sum_{i=1}^m {\rm Sw}_i({\cal F}) 
\cdot 
\chi_c(B_i^\circ,{\mathbb Q}_\ell) 
+ 
\sum_{x\in \Sigma} 
\text{Sw}_x({\cal F})
\label{eqlaumon}
\end{equation}
is proved.

To compare the formula
(\ref{eqlaumon})
with (\ref{eqgos}),
we give a slight reformulation.
Let $\pi_i:{\overline B}_i\to B_i$ be the normalization
for each irreducible component of dimension 1.
For each closed point $x\in \Sigma$,
we put
$$S_x({\cal F})= - {\rm Sw}_x({\cal F}) +
\sum_{i=1}^m{\rm Sw}_i({\cal F})\cdot |\pi_i^{-1}(x)|.$$
Then, the formula
(\ref{eqlaumon}) is equivalent to
\begin{equation}
\chi_c(U, {\cal F})= 
{\rm rank} \,{\cal F}\cdot\chi_c(U,{\mathbb Q}_\ell) - 
\left(\sum_{i=1}^m
{\rm Sw}_i({\cal F})
\chi({\overline B}_i, {\mathbb Q}_\ell) +
 \sum_{x\in \Sigma} S_x({\cal F})\right).
\label{eqlaumon2}
\end{equation}

We compute the Swan class ${\rm Sw}({\cal F})$
assuming the condition (NF)
and give a relation with $S_x({\cal F})$.
By Lemma \ref{lmalter2},
there exist a finite \'etale Galois covering $V
\to U$ that trivializes ${\cal F}$ and a
Cartesian diagram
$$\begin{CD}
V@>{\subset}>> Y\\
@VfVV @VV{\bar f}V\\
U@>{\subset}>> X
\end{CD}\leqno{\rm (\ref{eqXY})}
$$
such that $Y$ is smooth,
$Y\to X$ is proper,
$V\subset Y$ is the complement of
a divisor with simple normal crossings
and that the action of
$G={\rm Gal}(V/U)$
is extended to an admissible action on $Y$.
We may further assume that there exist
a proper scheme $X'$ containing $U$
as the complement of a Cartier divisor
and that $f:V\to U$
is extended to a morphism $Y\to X'$.
Furthermore, by the assumption (NF), we may assume
the following condition (NF$'$)
is satisfied
where
$\{\eta_{i1},\ldots,\eta_{ik_i}\}$ denotes the
inverse image of $\xi_i$ in $Y$ for
$i=1,\ldots,m$.

\medskip
(NF$'$)
For each $i,j$,
the extension
$\kappa(\eta_{ij})$
is {\it separable} over
$\kappa(\xi_i)$.

\medskip
Let $\sigma\neq 1$ be an element
of the Galois group $G={\rm Gal}(V/U)$.
For a generic point $\eta_{ij}$ as above,
we put
$m_{ij}(\sigma)=
{\rm length}\ {O_{Y^\sigma_{\log},\eta_{ij}}}$.
We define a divisor $D_\sigma$ of $Y$
by
$D_\sigma=\sum_{i,j}m_{ij}(\sigma)D_{ij}$
where
$D_{ij}$ is the closure of $\{\eta_{ij}\}$.
The Cartier divisor $D_\sigma$
is a closed subscheme of $Y^\sigma_{\log}$.
We define the residual subscheme
$R_\sigma\subset Y^\sigma_{\log}$ to
be the closed subscheme of $Y$
satisfying
$I_{Y^\sigma_{\log}}=I_{D_\sigma}I_{R_\sigma}$
where $I_Z$ denotes the ideal sheaf
of $O_Y$ defining a closed subscheme $Z\subset Y$.
Then, by the residual intersection formula
\cite{fulton} Theorem 9.2,
we have
\begin{align*}
(\Gamma_\sigma,\Delta_Y)^{\log}_{\mathbb Z}
=&\
-s_{V/U}(\sigma)\\
=&\
\{c(\Omega^1_{Y/F}(\log D))^*\cap 
(1+D_\sigma)^{-1}\cap D_\sigma\}_{\dim 0}
+{\mathbb R}_\sigma\\
=&\
-\sum_{i,j}m_{ij}(\sigma)
(c_1(\Omega^1_{Y/F}(\log D)) 
+D_\sigma)\cap D_{ij}
+{\mathbb R}_\sigma
\end{align*}
where
${\mathbb R}_\sigma
=\{c(\Omega^1_{Y/F}(\log D)\otimes O_Y(-D_\sigma))^*\cap 
s(R_\sigma/Y)\}_{\dim 0}$
is a 0-cycle class supported on
the inverse image of finitely many closed points of $B$.

To compute the first term in the right hand side,
we define a complex ${\cal K}_{ij}$
of $O_{D_{ij}}$-modules by
\begin{equation}
\begin{CD}
{\cal K}_{ij}=[
\varphi_{ij}^*(\Omega^1_{{\bar
B}_i/F})@>>>
 \Omega^1_{Y/F}(\log D)|_{D_{ij}}
@>\alpha>>  O_Y(-D_\sigma)|_{D_{ij}}].
\end{CD}
\label{eqNF}
\end{equation}
The sheaf
$\Omega^1_{Y/F}(\log D)|_{D_{ij}}$
is put on degree 0,
the map $\varphi_{ij}:D_{ij}\to {\overline B}_i$ is the natural one
and the map $\alpha$ is defined by
$da\mapsto \sigma(a)-a$
and $d\log b \mapsto \sigma(b)/b-1$.
By the assumption (NF$'$), the cohomology sheaves
${\cal H}^q({\cal K}_{ij})$
are 0 except for
$q=0,1$
and are supported
on finitely many closed points
for $q=0,1$.
Thus we have
$$(c_1(\Omega^1_{Y/F}(\log D)) 
+D_\sigma)\cap D_{ij}
=\varphi_{ij}^*c_1(\Omega^1_{\overline B_i/F})
+[{\cal H}^*({\cal K}_{ij})]$$
where
$[{\cal H}^*({\cal K}_{ij})]=
[{\cal H}^0({\cal K}_{ij})]-
[{\cal H}^1({\cal K}_{ij})]$.
Let $Z_0(B)$
denote the free abelian group
generated by
the classes of
the closed points in $B$.
We define a 0-cycle
$S_\sigma\in Z_0(B)$ by
$$S_\sigma=
{\bar f}_*(\sum_{i,j}m_{ij}(\sigma)
[{\cal H}^*({\cal K}_{ij})]
-{\mathbb R}_\sigma)$$
and put
$m_i(\sigma)
=\sum_jm_{ij}(\sigma)
[\kappa(\eta_{ij})
:\kappa(\xi_i)]$.
Then, we obtain
$$
{\bar f}_*s_{V/U}(\sigma)
=\sum_{i=1}^mm_i(\sigma)
g_{i*}(c_1(\Omega^1_{\overline B_i/F})
\cap [\overline B_i])
+S_\sigma$$
where $g_i:\overline B_i\to X$
is the natural map.
We define a 0-cycle $S_{\cal F}\in 
Z_0(B)\otimes{\mathbb Q}$ by
$$\frac 1{|G|}\sum_{\sigma \in G_{(p)}\setminus\{1\}}
\left(\dim_{{\mathbb F}_\ell} M^{\sigma}-
\frac{\dim_{{\mathbb F}_\ell} M^{\sigma^p}/M^{\sigma}}{p-1}
-\dim M\right)\cdot S_\sigma .$$
By ${\rm Sw}_i({\cal F})=
-\frac 1{|G|}
\sum_{\sigma\in G_{(p)}\setminus\{1\}}
m_i(\sigma)(
\dim_{{\mathbb F}_\ell} M^{\sigma}-
\dim_{{\mathbb F}_\ell}( M^{\sigma^p}/M^{\sigma})/(p-1)
-\dim M)$
and Lemma \ref{lmSwch}.1,
we have
\begin{eqnarray*}
{\rm Sw}({\cal F})
&=&\frac 1{|G|}\sum_{\sigma\in G_{(p)}\setminus\{1\}}
(
\dim_{{\mathbb F}_\ell} M^{\sigma}-
\frac{\dim_{{\mathbb F}_\ell} M^{\sigma^p}/M^{\sigma}}{p-1}
-\dim_{{\mathbb F}_\ell} M)
\cdot \bar f_*s_{V/U}(\sigma)\\
&=&
-\sum_{i=1}^m{\rm Sw}_i({\cal F})
g_{i*}(c_1(\Omega^1_{\overline B_i/F})\cap [\overline B_i])
+S_{\cal F}.
\end{eqnarray*}
Since
$\chi(\overline B_i,{\mathbb Q}_\ell)=-
{\rm deg}
(c_1(\Omega^1_{\overline B_i/F})\cap [\overline B_i])$,
the formula (\ref{eqgos})
together with the following proposition
will imply the formula
(\ref{eqlaumon2}).

\begin{pr}\label{prlaumon}
Under the notation above,
we have an equality
\begin{equation}
S_{\cal F}=
\sum_{x\in \Sigma} S_x({\cal F})[x]
\label{eqSwx}
\end{equation}
in $Z_0(B)$.
\end{pr}

In \cite{KSS} Theorem (6.7),
the invariant ${\rm Sw}_x({\cal F})$
is shown to be equal to
another invariant that is defined in \cite{SS}
using intersection classes
without introducing log products.
A similar computation 
gives a proof of
Proposition \ref{prlaumon}
but we leave the detail to the reader.

\bigskip

\begin{tabbing}
xxxx\=
xxxxxxxxxxxxxxxxxxxxxxxxxxxxxxx\= xxxxxxxxxxxxxxxxxxxxxxxxxxxxxxxxxxxxxx\kill
\> Kazuya Kato\>                     Takeshi Saito\\
\> Department of Mathematics,\>      Department of Mathematical Sciences,\\
\> Kyoto University,\>               University of Tokyo,\\
\> Kyoto 606-8502 Japan\>            Tokyo 153-8914 Japan\\
\> kazuya@math.kyoto-u.ac.jp\>       t-saito@ms.u-tokyo.ac.jp
\end{tabbing}

\end{document}